\documentclass[11pt]{amsart}
\usepackage{mathrsfs}
\usepackage{amssymb,latexsym}
\usepackage{pstricks,color}

 \numberwithin{equation}{section}

\newtheorem{theorem}{Theorem}[section]
\newtheorem{proposition}[theorem]{Proposition}
\newtheorem{corollary}[theorem]{Corollary}

\newtheorem{remark}[theorem]{Remark}

\newtheorem{lemma}[theorem]{Lemma}

\def\qed{\hfill $\Box$}
\def\pf{\noindent {\it Proof.} }

\title{Some statistics on generalized Motzkin paths with vertical steps }

\begin{document}
\maketitle
\begin{center}
Yidong Sun\footnote{Corresponding author: Yidong Sun.}, Di Zhao$^{2}$, Wenle Shi$^{3}$ and Weichen Wang$^{4}$

School of Science, Dalian Maritime University, 116026 Dalian, P.R. China\\[5pt]

{\it Emails: $^{1}$sydmath@dlmu.edu.cn, $^{2}$zd1129@dlmu.edu.cn, $^{3}$wlshi@dlmu.edu.cn, $^{4}$weichenw@dlmu.edu.cn }

\end{center}\vskip0.2cm

\subsection*{Abstract} Recently, several authors have considered lattice paths with various steps, including vertical steps permitted. In this paper, we consider a kind of generalized Motzkin paths, called {\it G-Motzkin paths} for short, that is lattice paths from $(0, 0)$ to $(n, 0)$ in the first quadrant of the $XOY$-plane that consist of up steps $\mathbf{u}=(1, 1)$, down steps $\mathbf{d}=(1, -1)$, horizontal steps $\mathbf{h}=(1, 0)$ and vertical steps $\mathbf{v}=(0, -1)$. We mainly count the number of G-Motzkin paths of length $n$ with given number of $\mathbf{z}$-steps for $\mathbf{z}\in \{\mathbf{u}, \mathbf{h}, \mathbf{v}, \mathbf{d}\}$, and enumerate the statistics ``number of $\mathbf{z}$-steps" at given level in G-Motzkin paths for $\mathbf{z}\in \{\mathbf{u}, \mathbf{h}, \mathbf{v}, \mathbf{d}\}$, some explicit formulas and combinatorial identities are given by bijective and algebraic methods, some enumerative results are linked with Riordan arrays according to the structure decompositions of G-Motzkin paths. We also discuss the statistics ``number of $\mathbf{z}_1\mathbf{z}_2$-steps" in G-Motzkin paths for $\mathbf{z}_1, \mathbf{z}_2\in \{\mathbf{u}, \mathbf{h}, \mathbf{v}, \mathbf{d}\}$, the exact counting formulas except for $\mathbf{z}_1\mathbf{z}_2=\mathbf{dd}$ are obtained by the Lagrange inversion formula and their generating functions.

\medskip

{\bf Keywords}: Dyck path, G-Motzkin path, Catalan number, Riordan array.

\noindent {\sc 2010 Mathematics Subject Classification}: Primary 05A15; Secondary 05A05, 05A19.

{\bf \section{ Introduction } }

Lattice paths have been studied by many mathematicians and have produced numerous interesting and important results. Research in this area has resulted in well known classes of lattice paths such as those named after Dyck \cite{Deutsch99,Stanley,StanleyEC2}, Motzkin \cite{Aigner,BarPinSpr,DonShap}, Schr\"{o}der \cite{Comtet} and Delannoy \cite{BandSchw}. They are used in physics \cite{Rensburg}, computer science \cite{Knuth,Viennot}, biology \cite{Choi,DosSvrVel,HofSchusSta,Robeva,SheStr,ZukSan} and probability theory \cite{Defant,Koshy,Mand,Mohanty,Narayana,Takacs}. We refer the reader to the wonderful survey by Humphreys \cite{Humphreys} for additional historical information.

A {\it Dyck path} of length $2n$ is a lattice path from $(0, 0)$ to $(2n, 0)$ in the first quadrant of the XOY-plane that consists of up steps $\mathbf{u}=(1, 1)$ and down steps $\mathbf{d}=(1, -1)$. Let $\mathcal{C}_n$ be the set of Dyck paths of length $2n$. It is well known \cite{Deutsch99,Stanley,StanleyEC2} that $|\mathcal{C}_n|=C_n=\frac{1}{n+1}\binom{2n}{n}$, the $n$th Catalan number \cite[A000108]{Sloane}, has the generating function
\begin{eqnarray}\label{eqn 1.1}
C(x)=\sum_{n\geq 0}C_nx^n=\frac{1-\sqrt{1-4x}}{2x}
\end{eqnarray}
with the relation $C(x)=1+xC(x)^2=\frac{1}{1-xC(x)}$.

A {\it Motzkin path} of length $n$ is a lattice path from $(0, 0)$ to $(n, 0)$ in the first quadrant of the XOY-plane that consists of up steps $\mathbf{u}=(1, 1)$, down steps $\mathbf{d}=(1, -1)$ and horizontal steps $\mathbf{h}=(1, 0)$. Let $\mathcal{M}_n$ be the set of Motzkin paths of length $n$. It is well known \cite{BarPinSpr,DonShap,StanleyEC2} that $|\mathcal{M}_n|=M_n$, the $n$th Motzkin number \cite[A001006]{Sloane}, has the generating function
\begin{eqnarray*}
M(x)=\sum_{n\geq 0}M_nx^n=\frac{1-x-\sqrt{1-2x-3x^2}}{2x^2}
\end{eqnarray*}
with the relation $M(x)=1+xM(x)+x^2M(x)^2=\frac{1}{1-x}C(\frac{x^2}{(1-x)^2})$. This implies the following relation between the Catalan numbers $C_n$ and the Motzkin numbers $M_n$ \cite{DonShap}, i.e.,
\begin{eqnarray*}
M_n=\sum_{k=0}^{[\frac{n}{2}]}\binom{n}{2k}C_k.
\end{eqnarray*}
In fact, the sequence $M_{n,k}=\binom{n}{2k}C_k$ counts the number of Motzkin paths of length $n$ with $k$ $\mathbf{d}$-steps. The first values of $M_{n,k}=\binom{n}{2k}C_k$ are illustrated in Table 1.1.

\begin{center}
\begin{eqnarray*}
\begin{array}{c|cccccc}\hline
n/k & 0   & 1   & 2    & 3    & 4    & 5       \\\hline
  0 & 1   &     &      &      &      &        \\
  1 & 1   &     &      &      &      &        \\
  2 & 1   & 1   &      &      &      &         \\
  3 & 1   & 3   &      &      &      &         \\
  4 & 1   & 6   & 2    &      &      &         \\
  5 & 1   & 10  & 10   &      &      &        \\
  6 & 1   & 15  & 30   & 5    &      &        \\
  7 & 1   & 21  & 70   & 35   &      &        \\
  8 & 1   & 28  & 140  & 140  &  14  &        \\
  9 & 1   & 36  & 252  & 420  &  126 &        \\
 10 & 1   & 45  & 420  & 1050 &  630 &  42    \\ \hline
\end{array}
\end{eqnarray*}
Table 1.1. The first values of $M_{n,k}$.
\end{center}

A {\it Schr\"{o}der path} of length $2n$ is a path from $(0, 0)$ to $(2n, 0)$ in the first quadrant of the XOY-plane that consists of up steps $\mathbf{u}=(1, 1)$, down steps $\mathbf{d}=(1, -1)$ and horizontal steps $\mathbf{H}=(2, 0)$. Let $\mathcal{S}_n$ be the set of Schr\"{o}der paths of length $2n$. It is well known \cite{StanleyEC2} that $|\mathcal{S}_n|=R_n$, the $n$th large Schr\"{o}der number \cite[A006318]{Sloane}, has the generating function
\begin{eqnarray*}
R(x)=\sum_{n\geq 0}R_nx^n=\frac{1-x-\sqrt{1-6x+x^2}}{2x}
\end{eqnarray*}
with the relation $R(x)=1+xR(x)+xR(x)^2=\frac{1}{1-x}C(\frac{x}{(1-x)^2})$. This shows the following relation between the Catalan numbers $C_n$ and the large Schr\"{o}der numbers $R_n$, i.e.,
\begin{eqnarray*}
R_n=\sum_{k\geq 0}M_{n+k,k}=\sum_{k=0}^{n}\binom{n+k}{2k}C_k.
\end{eqnarray*}
In fact, the sequence $R_{n,k}=\binom{n+k}{2k}C_k$ counts the number of Schr\"{o}der paths of length $2n$ with $k$ $\mathbf{d}$-steps. Also, the $n$th little Schr\"{o}der number $r_n$ counts the number of Schr\"{o}der paths of length $2n$ with no $\mathbf{H}$-steps at $X$-axis \cite[A001003]{Sloane}.

Recently, several authors \cite{Dziem-A,Dziem-B,Dziem-C,IrvMelRusk,IrvRusk,YanZhang,YangYang} have considered lattice paths with various steps, including vertical steps permitted. In this paper, we consider a kind of generalized Motzkin paths, called {\it G-Motzkin paths} for short. That is, a G-Motzkin path of length $n$ is a lattice path from $(0, 0)$ to $(n, 0)$ in the first quadrant of the XOY-plane that consists of up steps $\mathbf{u}=(1, 1)$, down steps $\mathbf{d}=(1, -1)$, horizontal steps $\mathbf{h}=(1, 0)$ and vertical steps $\mathbf{v}=(0, -1)$. See Figure 1 for a G-Motzkin path of length $24$.

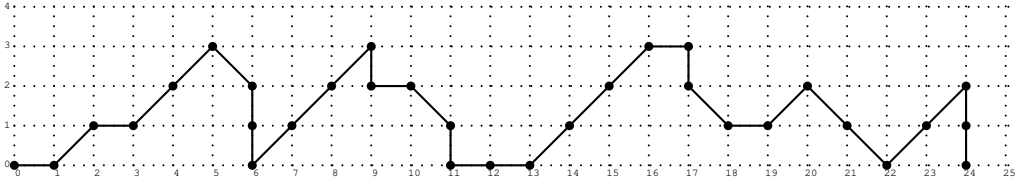
\begin{figure}[h] \setlength{\unitlength}{0.5mm}

\begin{center}
\begin{pspicture}(13,2.2)
\psset{xunit=15pt,yunit=15pt}\psgrid[subgriddiv=1,griddots=4,
gridlabels=4pt](0,0)(25,4)

\psline(0,0)(1,0)(2,1)(3,1)(5,3)(6,2)(6,1)(6,0)(9,3)(9,2)(10,2)(11,1)(11,0)(13,0)(16,3)(17,3)(17,2)(18,1)(19,1)
\psline(19,1)(20,2)(22,0)(24,2)(24,0)

\pscircle*(0,0){0.06}\pscircle*(1,0){0.06}\pscircle*(2,1){0.06}
\pscircle*(3,1){0.06}\pscircle*(4,2){0.06}\pscircle*(5,3){0.06}
\pscircle*(6,2){0.06}\pscircle*(6,1){0.06}\pscircle*(6,0){0.06}

\pscircle*(7,1){0.06}\pscircle*(8,2){0.06}\pscircle*(9,2){0.06}
\pscircle*(9,3){0.06}\pscircle*(10,2){0.06}\pscircle*(11,1){0.06}\pscircle*(11,0){0.06}
\pscircle*(12,0){0.06}\pscircle*(13,0){0.06}\pscircle*(14,1){0.06}
\pscircle*(15,2){0.06}\pscircle*(16,3){0.06}\pscircle*(17,3){0.06}\pscircle*(17,2){0.06}
\pscircle*(18,1){0.06}\pscircle*(19,1){0.06}\pscircle*(20,2){0.06}
\pscircle*(21,1){0.06}\pscircle*(22,0){0.06}\pscircle*(23,1){0.06}\pscircle*(24,2){0.06}
\pscircle*(24,1){0.06}\pscircle*(24,0){0.06}

\end{pspicture}
\end{center}\vskip0.2cm

\caption{\small A G-Motzkin path of length $24$.}

\end{figure}

Let $\varepsilon$ be the empty path, that is a dot path. If $\mathbf{P}_1$ and $\mathbf{P}_2$ are G-Motzkin paths, then we define $\mathbf{P}_1\mathbf{P}_2$ as the concatenation of $\mathbf{P}_1$ and $\mathbf{P}_2$. For example, $\mathbf{P}_1=\mathbf{uhuduuvvdhh}$ and $\mathbf{P}_2=\mathbf{uhuhvvuudd}$, then $\mathbf{P}_1\mathbf{P}_2=\mathbf{uhuduuvvdhh}\mathbf{uhuhvvuudd}$.

A point of a G-Motzkin path with ordinate $\ell$ is said to be at {\it level} $\ell$. A step of a G-Motzkin path is said to be at level $\ell$ if the ordinate of its endpoint is $\ell$. A {\it $\mathbf{ud}$-peak} ({\it $\mathbf{uv}$-peak}) in a G-Motzkin path is an occurrence of $\mathbf{ud}$ ($\mathbf{uv}$). A {\it $\mathbf{du}$-valley} ({\it $\mathbf{vu}$-valley}) in a G-Motzkin path is an occurrence of $\mathbf{du}$ ($\mathbf{vu}$). A {\it peak (valley)} in a G-Motzkin path is a $\mathbf{ud}$-peak or $\mathbf{uv}$-peak ($\mathbf{du}$-valley or $\mathbf{vu}$-valley). By the {\it hight of a peak} ({\it valley}) we mean the level of the intersection point of its two steps. By a {\it return step} we mean a $\mathbf{d}$-step or $\mathbf{v}$-step at level $0$. A {\it matching step} of a $\mathbf{u}$-step at level $k\geq 1$ in a G-Motzkin path is the leftmost step among all $\mathbf{d}$-steps or $\mathbf{v}$-steps at level $k-1$ right to the $\mathbf{u}$-step. A G-Motzkin path $\mathbf{P}$ is said to be {\it primitive} if $\mathbf{P}=\mathbf{u}\mathbf{P}'\mathbf{d}$ or $\mathbf{P}=\mathbf{u}\mathbf{P}'\mathbf{v}$ for certain G-Motzkin path $\mathbf{P}'$.

In the present paper, we concentrate on several statistics in G-Motzkin paths. Precisely, the next section mainly counts the number of G-Motzkin paths of length $n$ with given number of $\mathbf{z}$-steps for $\mathbf{z}\in \{\mathbf{u}, \mathbf{h}, \mathbf{v}, \mathbf{d}\}$, some explicit formulas and combinatorial identities are given by bijective and algebraic methods. The third section mainly focuses on the enumeration of statistics ``number of $\mathbf{z}$-steps" at given level in G-Motzkin paths for $\mathbf{z}\in \{\mathbf{u}, \mathbf{h}, \mathbf{v}, \mathbf{d}\}$, the enumerative results are linked with Riordan arrays according to the structure decompositions of G-Motzkin paths. The last section discusses
the statistics ``number of $\mathbf{z}_1\mathbf{z}_2$-steps" in G-Motzkin paths for $\mathbf{z}_1, \mathbf{z}_2\in \{\mathbf{u}, \mathbf{h}, \mathbf{v}, \mathbf{d}\}$, the exact counting formulas except for $\mathbf{z}_1\mathbf{z}_2=\mathbf{dd}$ are provided according to the method of the first return decomposition of G-Motzkin paths and the Lagrange inversion formula.

\section{The statistics ``number of $\mathbf{z}$-steps" in G-Motzkin paths }

In this section, we first consider the weighted G-Motzkin paths. The {\it weight} of each step of a G-Motzkin path $\mathbf{P}$ is assigned as follows. The $\mathbf{u}$-step, $\mathbf{h}$-step, $\mathbf{v}$-step and $\mathbf{d}$-step are weighted respectively by $1, a, b$ and $c$. The {\it weight} of $\mathbf{P}$, denoted by $w(\mathbf{P})$, is the product of the weight of each step of $\mathbf{P}$. For example, $w(\mathbf{uhuduuvvdhh})=a^3b^2c^2$. The {\em weight} of a subset $\mathcal{A}$ of the set $\mathcal{G}(a, b, c)$ of all weighted G-Motzkin paths, denoted by $w(\mathcal{A})$, is the sum of the total weights of all paths in $\mathcal{A}$.
Denoted by $w(\mathcal{G}_{n}(a, b, c))=G_n(a, b, c)$ the weight of the set $\mathcal{G}_{n}(a, b, c)$ of all weighted G-Motzkin paths of length $n$. When $a=b=c=1$, we write
$\mathcal{G}=\mathcal{G}(1, 1, 1), \mathcal{G}_{n}=\mathcal{G}_{n}(1, 1, 1), G_n=G_n(1,1,1)$ for short.

Let $G(a,b,c; x)=\sum_{n=0}^{\infty}G_n(a, b, c)x^n$ be the generating function. According to the
first return decomposition, a G-Motzkin path $\mathbf{P}$ can be decomposed as one of the following four forms:
$$\mathbf{P}=\varepsilon, \ \mathbf{P}=\mathbf{h}\mathbf{Q}_1,\ \mathbf{P}=\mathbf{u}\mathbf{Q}_1\mathbf{v}\mathbf{Q}_2 \ \mbox{or}\  \mathbf{P}=\mathbf{u}\mathbf{Q}_1\mathbf{d}\mathbf{Q}_2, $$
where $\mathbf{Q}_1$ and $\mathbf{Q}_2 $ are (possibly empty) G-Motzkin paths. Then we get the relation
\begin{eqnarray}\label{eqn 2.1}
G(a,b,c; x)=1+axG(a,b,c; x)+bxG(a,b,c; x)^2+cx^2G(a,b,c; x)^2.
\end{eqnarray}
Solve this, we have
\begin{eqnarray}\label{eqn 2.2}
G(a,b,c; x)=\frac{1-ax-\sqrt{(1-ax)^2-4x(b+cx)}}{2x(b+cx)}=\frac{1}{1-ax}C\Big(\frac{x(b+cx)}{(1-ax)^2}\Big).
\end{eqnarray}
By (\ref{eqn 1.1}), taking the coefficient of $x^n$ in $G(a,b,c; x)$, we derive that
\begin{proposition}
For any integer $n\geq 0$, there holds
\begin{eqnarray*}
G_n(a, b, c) \hskip-.22cm &=&\hskip-.22cm \sum_{k=0}^{n}\sum_{j=0}^{n-k}\binom{k}{j}\binom{n+k-j}{2k}C_ka^{n-k-j}b^{k-j}c^{j}  \\
             \hskip-.22cm &=&\hskip-.22cm \sum_{k=0}^{n}\sum_{j=0}^{k}\binom{k}{j}\binom{n+j}{2k}C_ka^{n-2k+j}b^{j}c^{k-j}  \\
             \hskip-.22cm &=&\hskip-.22cm \sum_{k=0}^{n}\sum_{j=0}^{n-k}\binom{2k+j}{j}\binom{k}{n-k-j}C_ka^{j}b^{2k+j-n}c^{n-k-j}.
\end{eqnarray*}
\end{proposition}

Set $T=xG(a,b,c; x)$, (\ref{eqn 2.1}) produces
\begin{eqnarray*}
T=x\frac{1+aT+cT^2}{1-bT},
\end{eqnarray*}
using the Lagrange inversion formula \cite{Gessel}, taking the coefficient of $x^{n+1}$ in $T$ in three different ways, we derive that
\begin{proposition}
For any integer $n\geq 0$, there holds
\begin{eqnarray}
G_n(a, b, c) \hskip-.22cm &=&\hskip-.22cm \frac{1}{n+1}\sum_{k=0}^{[\frac{n}{2}]}\sum_{j=0}^{n-2k}\binom{n+1}{k}\binom{n+1-k}{j}\binom{2n-2k-j}{n-2k-j} a^{j}b^{n-2k-j}c^{k}  \nonumber\\
             \hskip-.22cm &=&\hskip-.22cm \frac{1}{n+1}\sum_{k=0}^{n}\sum_{j=0}^{[\frac{n-k}{2}]}\binom{n+1}{k}\binom{n+1-k}{j}\binom{2n-k-2j}{n-k-2j} a^{k}b^{n-k-2j}c^{j}  \nonumber\\
             \hskip-.22cm &=&\hskip-.22cm \frac{1}{n+1}\sum_{k=0}^{n}\sum_{j=0}^{n-k}\binom{n+1}{k}\binom{k}{j}\binom{2n-k-j}{n-k-j} a^{k-j}b^{n-k-j}c^{j}.  \label{eqn 2.2b}
\end{eqnarray}
\end{proposition}

Clearly, $G_n=G_n(1,1,1)$ is the number of G-Motzkin paths of length $n$ with the generating function
\begin{eqnarray}\label{eqn 2.3}
G(x)=\sum_{n=0}^{\infty}G_nx^n=\frac{1-x-\sqrt{1-6x-3x^2}}{2x(1+x)}=\frac{1}{1-x}C\Big(\frac{x(1+x)}{(1-x)^2}\Big).
\end{eqnarray}
The explicit form of $G(x)$ was given by Drake \cite{Drake} by counting lattice paths without regard to area and by Dziemia\'{n}czuk \cite{Dziem-A} by counting special lattice paths with four types of steps. The sequence
$$(G_n)_{n\geq 0} = (1, 2, 7, 29, 133, 650, 3319, 17498, 94525, 520508, 2910895, . . .)$$
is denoted by A064641 in OEIS \cite{Sloane}, and has the formula \cite{Dziem-A}
\begin{eqnarray*}
G_{n}=\frac{1}{n+1}\sum_{k=0}^{n}\sum_{j=0}^{k}\binom{n+1}{j}\binom{j}{k-j}\binom{2n-k}{n}
\end{eqnarray*}
which obeys the recurrence relation \cite{Sloane}
$$(n+1)G_n=(5n-4)G_{n-1}+9(n-1)G_{n-2}+3(n-2)G_{n-3}. $$

If setting $T=xG(1,1,1; x)$, (\ref{eqn 2.1}) produces
\begin{eqnarray*}
T=x\frac{1+T+T^2}{1-T}=x\frac{1-T^3}{(1-T)^2},
\end{eqnarray*}
using the Lagrange inversion formula, taking the coefficient of $x^{n+1}$ in $T$, one has another simple formula for $G_n$, namely,
\begin{eqnarray*}
G_{n}=\frac{1}{n+1}\sum_{k=0}^{[\frac{n}{3}]}(-1)^k\binom{n+1}{k}\binom{3n-3k+1}{n-3k}.
\end{eqnarray*}

Let $G_{n,k}=\frac{1}{n+1}\binom{n+1}{k}\binom{3n-3k+1}{n-3k}$, the first values of $G_{n,k}$ are illustrated in Table 2.0.
\begin{center}
\begin{eqnarray*}
\begin{array}{c|ccccccc}\hline
n/k & 0      & 1      & 2       & 3     & 4    & 5        \\\hline
  0 & 1      &        &         &       &      &         \\
  1 & 2      &        &         &       &      &         \\
  2 & 7      &        &         &       &      &          \\
  3 & 30     & 1      &         &       &      &          \\
  4 & 143    & 10     &         &       &      &          \\
  5 & 728    & 78     &         &       &      &         \\
  6 & 3876   & 560    & 3       &       &      &         \\
  7 & 21318  & 3876   & 56      &       &      &         \\
  8 & 120175 & 26334  & 684     &       &      &         \\
  9 & 690690 & 177100 & 6930    & 12    &      &         \\\hline
\end{array}
\end{eqnarray*}
Table 2.0. The first values of $G_{n,k}$.
\end{center}

The following is an interesting identity related to $G_{n,k}$, that is,
\begin{eqnarray*}
\sum_{k=0}^{[\frac{n}{2}]}G_{n+k, k}=\frac{1}{n+1}\sum_{k=0}^{[\frac{n}{2}]}\binom{n+k}{k}\binom{3n+1}{n-2k}=2^nC_n,
\end{eqnarray*}
which can be proved as follows,
\begin{eqnarray*}
\frac{1}{n+1}\sum_{k=0}^{[\frac{n}{2}]}\binom{n+k}{k}\binom{3n+1}{n-2k}
\hskip-.22cm &=&\hskip-.22cm \frac{1}{n+1}\sum_{k=0}^{[\frac{n}{2}]}\binom{n+k}{k}[x^{n-2k}](1+x)^{3n+1} \\
\hskip-.22cm &=&\hskip-.22cm \frac{1}{n+1}[x^{n}]\sum_{k=0}^{\infty}\binom{n+k}{k}x^{2k}(1+x)^{3n+1}     \\
\hskip-.22cm &=&\hskip-.22cm \frac{1}{n+1}[x^{n}]\frac{(1+x)^{3n+1}}{(1-x^2)^{n+1}} =\frac{1}{n+1}[x^{n}]\frac{(1+x)^{2n}}{(1-x)^{n+1}}            \\
\hskip-.22cm &=&\hskip-.22cm \frac{1}{n+1}\sum_{k=0}^{n}\binom{2n}{n-k}\binom{n+k}{k} = \frac{1}{n+1}\sum_{k=0}^{n}\binom{n}{k}\binom{2n}{n}         \\
\hskip-.22cm &=&\hskip-.22cm 2^nC_n,
\end{eqnarray*}
where $[x^n]f(x)$ denotes the coefficient of $x^{n}$ in $f(x)$.

\subsection{ The statistics ``number of $\mathbf{v}$-steps" }

Let $V_{n,i}$ denote the number of G-Motzkin paths of length $n$ with $i$ $\mathbf{v}$-steps, the first values of $V_{n,i}$ are illustrated in Table 2.1.

\begin{center}
\begin{eqnarray*}
\begin{array}{c|ccccccc}\hline
n/i & 0   & 1   & 2    & 3    & 4    & 5    & 6     \\\hline
  0 & 1   &     &      &      &      &      &      \\
  1 & 1   & 1   &      &      &      &      &      \\
  2 & 2   & 3   & 2    &      &      &      &       \\
  3 & 4   & 10  & 10   & 5    &      &      &       \\
  4 & 9   & 30  & 45   & 35   &  14  &      &       \\
  5 & 21  & 90  & 175  & 196  & 126  & 42   &        \\
  6 & 51  & 266 & 644  & 924  & 840  & 462  & 132     \\\hline
\end{array}
\end{eqnarray*}
Table 2.1. The first values of $V_{n,i}$.
\end{center}

This shows that there is a close relation between $V_{n,i}$ and $C_n$. Exactly, any Dyck paths $\mathbf{P}$ of length $2k$ can be extended to G-Motzkin paths of length $n$ with $i$ $\mathbf{v}$-steps for $i\leq k\leq n$. Note that there exist $2k+1$ points and $k$ $\mathbf{d}$-steps in $\mathbf{P}$, so there are $\binom{k}{i}$ ways to replace $i$ $\mathbf{d}$-steps by $\mathbf{v}$-steps and there are $\big(\binom{2k+1}{n-2k+i}\big)=\binom{n+i}{2k}$ ways to insert repeatedly $n-2k+i$ $\mathbf{h}$-steps into $2k+1$ points of $\mathbf{P}$ to form G-Motzkin paths of length $n$ with $i$ $\mathbf{v}$-steps. Summarizing over $k$, we have
\begin{theorem}\label{theom 2.1.1}
For any integers $n\geq i\geq 0$, there holds
\begin{eqnarray*}
V_{n,i}=\sum_{k=i}^{n}\binom{k}{i}\binom{n+i}{2k}C_k.
\end{eqnarray*}
\end{theorem}
Note that $V_{n-i,i}$ also counts the G-Motzkin paths of length $n-i$ with $n$ steps and $V_{n,i}$ is also the coefficient of $b^i$ in $G_n(1, b, 1)$ in (\ref{eqn 2.2b}) which has another expression
\begin{eqnarray*}
V_{n,i}=\frac{1}{n+1}\binom{n+i}{i}\sum_{k=0}^{n-i}\binom{n+1}{k}\binom{k}{n-k-i}.
\end{eqnarray*}

\subsection{ The statistics ``number of $\mathbf{h}$-steps" }

Let $H_{n,i}$ denote the number of G-Motzkin paths of length $n$ with $i$ $\mathbf{h}$-steps, the first values of $H_{n,i}$ are illustrated in Table 2.2.

\begin{center}
\begin{eqnarray*}
\begin{array}{c|ccccccc}\hline
n/i & 0   & 1   & 2    & 3    & 4    & 5    & 6     \\\hline
  0 & 1   &     &      &      &      &      &      \\
  1 & 1   & 1   &      &      &      &      &      \\
  2 & 3   & 3   & 1    &      &      &      &       \\
  3 & 9   & 13  & 6    & 1    &      &      &       \\
  4 & 31  & 55  & 36   & 10   & 1    &      &       \\
  5 & 113 & 241 & 200  & 80   & 15   & 1    &        \\
  6 & 431 & 1071& 1080 & 560  & 155  & 21   & 1     \\\hline
\end{array}
\end{eqnarray*}
Table 2.2. The first values of $H_{n,i}$.
\end{center}

Similarly, any Dyck paths $\mathbf{P}$ of length $2k$ can be extended to G-Motzkin paths of length $n$ with $i$ $\mathbf{h}$-steps for $[\frac{n-i}{2}]\leq k\leq n-i$. Note that there are $\binom{k}{n-i-k}$ ways to replace $n-i-k$ $\mathbf{d}$-steps by $\mathbf{v}$-steps and there are $\big(\binom{2k+1}{i}\big)=\binom{2k+i}{2k}$ ways to insert repeatedly $i$ $\mathbf{h}$-steps into $2k+1$ points of $\mathbf{P}$ to form G-Motzkin paths of length $n$ with $i$ $\mathbf{h}$-steps. Summarizing over $k$, we have
\begin{theorem}\label{theom 2.2.1}
For any integers $n\geq i\geq 0$, there holds
\begin{eqnarray*}
H_{n,i}=\sum_{k=[\frac{n-i}{2}]}^{n-i}\binom{2k+i}{2k}\binom{k}{n-i-k}C_k.
\end{eqnarray*}
\end{theorem}
Note that $H_{n,i}$ is also the coefficient of $a^i$ in $G_n(a, 1, 1)$ in (\ref{eqn 2.2b}) which has another expression
\begin{eqnarray*}
H_{n,i}=\frac{1}{n+1}\binom{n+1}{i}\sum_{j=0}^{[\frac{n-i}{2}]}\binom{n+1-i}{j}\binom{2n-i-2j}{n-i-2j},
\end{eqnarray*}
and the special case $G(-2, 1, 1; x)=\frac{1}{1+x}$ by (\ref{eqn 2.2}) deduces the following identity whose combinatorial proof is also provided.
\begin{theorem}\label{theom 2.2.2}
For any integer $n\geq 0$, there holds
\begin{eqnarray*}
\sum_{i=0}^{n}(-2)^{i}H_{n,i}=(-1)^{n}.
\end{eqnarray*}
\end{theorem}
\pf Let $\mathcal{H}_n^{e}$ ($\mathcal{H}_n^{o}$) denote the set of weighted G-Motzkin paths of length $n$ with even (odd) number of $\mathbf{h}$-steps such that each $\mathbf{h}$-step is weighted by $2$ (regarded as $\mathbf{h}_1$ and $\mathbf{h}_2$ for convenience) and other steps are weighted by $1$.
Clearly,
$$w(\mathcal{H}_n^{e})=\sum_{i\ even}2^{i}H_{n,i}\ \ \mbox{and}\ \ w(\mathcal{H}_n^{o})=\sum_{i\ odd}2^{i}H_{n,i}. $$
So it is sufficient to give a bijection $\phi$ between $\mathcal{H}_n^{e}/\{\mathbf{h}_1^n\}$ and $\mathcal{H}_n^{o}$ for $n$ even and between $\mathcal{H}_n^{e}$ and $\mathcal{H}_n^{o}/\{\mathbf{h}_1^n\}$ for $n$ odd. When $n\geq 2$, any $\mathbf{P}\in \mathcal{H}_n^{e}/\{\mathbf{h}_1^n\}$ for $n$ even or $\mathbf{P}\in \mathcal{H}_n^{e}$ for $n$ odd has at least one of the four subpaths, $\mathbf{d}$, $\mathbf{h}_1\mathbf{v}$, $\mathbf{h}_2$ and $\mathbf{uv}$, find the last one, say $\mathbf{x}$, $\mathbf{P}$ can be partitioned uniquely into $\mathbf{P}=\mathbf{P}_1\mathbf{x}\mathbf{P}_2$, where $\mathbf{P}_2=\mathbf{v}^k$ for certain $0\leq k<n$. Then define $\phi(\mathbf{P})=\mathbf{P}_1\mathbf{x}'\mathbf{P}_2$, where
\begin{eqnarray*}
\mathbf{x}'=\left\{
\begin{array}{rl}
 \mathbf{h}_1\mathbf{v},   &  \mbox{if}\ \mathbf{x}=\mathbf{d}, \\
 \mathbf{d},               &  \mbox{if}\ \mathbf{x}=\mathbf{h}_1\mathbf{v}, \\
 \mathbf{uv},              &  \mbox{if}\ \mathbf{x}=\mathbf{h}_2, \\
 \mathbf{h}_2,             &  \mbox{if}\ \mathbf{x}=\mathbf{uv}.
\end{array}\right.
\end{eqnarray*}
This way ensures that the number of $\mathbf{h}$-steps in $\phi(\mathbf{P})$ is one more or less than that in $\mathbf{P}\in \mathcal{H}_n^{e}/\{\mathbf{h}_1^n\}$ for $n$ even or in $\mathbf{P}\in \mathcal{H}_n^{e}$ for $n$ odd, so $\phi(\mathbf{P})\in \mathcal{H}_n^{o}$ for $n$ even and $\phi(\mathbf{P})\in \mathcal{H}_n^{o}/\{\mathbf{h}_1^n\}$ for $n$ odd. Moreover, $\mathbf{x}'$ in $\phi(\mathbf{P})$ is also the last one of the four subpaths, $\mathbf{d}$, $\mathbf{h}_1\mathbf{v}$, $\mathbf{h}_2$ and $\mathbf{uv}$. The reverse procedure can be handled similarly.
Hence, $\phi$ is a bijection (an involution in fact) between $\mathcal{H}_n^{e}/\{\mathbf{h}_1^n\}$ and $\mathcal{H}_n^{o}$ for $n$ even and between $\mathcal{H}_n^{e}$ and $\mathcal{H}_n^{o}/\{\mathbf{h}_1^n\}$ for $n$ odd. This completes the proof.    \qed\vskip0.1cm

\begin{theorem}\label{theom 2.2.3}
For any integers $n, m\geq 0$, there holds
\begin{eqnarray*}
\sum_{i=0}^{2n}(-1)^{i}\binom{2n}{i}H_{n+m+i,m+i}=C_n.
\end{eqnarray*}
\end{theorem}
\pf By Theorem \ref{theom 2.2.1}, we have
\begin{eqnarray*}
\sum_{i=0}^{2n}(-1)^{i}\binom{2n}{i}H_{n+m+i,m+i} \hskip-.22cm &=&\hskip-.22cm \sum_{i=0}^{2n}(-1)^{i}\binom{2n}{i}\sum_{k=[\frac{n}{2}]}^{n}\binom{2k+m+i}{2k}\binom{k}{n-k}C_k \\
\hskip-.22cm &=&\hskip-.22cm \sum_{k=[\frac{n}{2}]}^{n}\binom{k}{n-k}C_k\sum_{i=0}^{2n}(-1)^{i}\binom{2n}{i}\binom{2k+m+i}{2k}\\
\hskip-.22cm &=&\hskip-.22cm (-1)^{m}\sum_{k=[\frac{n}{2}]}^{n}\binom{k}{n-k}C_k\sum_{i=0}^{2n}\binom{2n}{2n-i}\binom{-2k-1}{m+i}\\
\hskip-.22cm &=&\hskip-.22cm (-1)^{m}\sum_{k=[\frac{n}{2}]}^{n}\binom{k}{n-k}C_k\binom{2n-2k-1}{2n+m}\\
\hskip-.22cm &=&\hskip-.22cm (-1)^{m}\binom{n}{0}C_n\binom{-1}{2n+m}   \ \ \ \ \ \  (k=n)     \\
\hskip-.22cm &=&\hskip-.22cm C_n,
\end{eqnarray*}
where the fourth equality follows from the Chu-Vandemonde identity. This completes the proof.    \qed\vskip0.2cm

\subsection{ The statistics ``number of $\mathbf{d}$-steps" }
Let $D_{n,i}$ denote the number of G-Motzkin paths of length $n$ with $i$ $\mathbf{d}$-steps, the first values of $D_{n,i}$ are illustrated in Table 2.3.

\begin{center}
\begin{eqnarray*}
\begin{array}{c|cccccc}\hline
n/i & 0       & 1       & 2       & 3       & 4    & 5       \\\hline
  0 & 1       &         &         &         &      &        \\
  1 & 2       &         &         &         &      &        \\
  2 & 6       & 1       &         &         &      &         \\
  3 & 22      & 7       &         &         &      &         \\
  4 & 90      & 41      & 2       &         &      &         \\
  5 & 394     & 231     & 25      &         &      &        \\
  6 & 1806    & 1289    & 219     & 5       &      &        \\
  7 & 8558    & 7183    & 1666    & 91      &      &        \\
  8 & 41586   & 40081   & 11780   & 1064    & 14   &        \\\hline
\end{array}
\end{eqnarray*}
Table 2.3. The first values of $D_{n,i}$.
\end{center}

Similarly, any Dyck paths $\mathbf{P}$ of length $2k$ can be extended to G-Motzkin paths of length $n$ with $i$ $\mathbf{d}$-steps for $i\leq k\leq n-i$. Note that there are $\binom{k}{i}$ ways to replace $k-i$ $\mathbf{d}$-steps by $\mathbf{v}$-steps and there are $\big(\binom{2k+1}{n-k-i}\big)=\binom{n+k-i}{2k}$ ways to insert repeatedly $n-k-i$ $\mathbf{h}$-steps into $2k+1$ points of $\mathbf{P}$ to form G-Motzkin paths of length $n$ with $i$ $\mathbf{d}$-steps. Summarizing over $k$, we have
\begin{theorem}\label{theom 2.3.1}
For any integers $n\geq i\geq 0$, there holds
\begin{eqnarray*}
D_{n,i}=\sum_{k=i}^{n-i}\binom{k}{i}\binom{n-i+k}{2k}C_k.
\end{eqnarray*}
\end{theorem}
Note that $D_{n,i}$ is also the coefficient of $c^i$ in $G_n(1, 1, c)$ in (\ref{eqn 2.2b}) which has another expression
\begin{eqnarray*}
D_{n,i}=\frac{1}{n+1}\binom{n+1}{i}\sum_{k=i}^{n-i}\binom{n+1-i}{k-i}\binom{2n-i-k}{n-i-k},
\end{eqnarray*}
$D_{2n,n}$ is the $n$-th Catalan number $C_n$ and $D_{n,0}$ is the $n$-th large Schr\"{o}der number $R_n$.
Since any G-Motzkin path of length $n$ with no $\mathbf{d}$-steps can generate a Schr\"{o}der path of length $2n$ by replacing each $\mathbf{h}$-step by an $\mathbf{H}$-step and each $\mathbf{v}$-step by a $\mathbf{d}$-step, and vice versa. The special cases $G(1, 1, -2; x)=\frac{1}{1-2x}$ and $G(1, 1, -1; x)=\frac{1}{1-x}C\big(\frac{x}{1-x}\big)$ by (\ref{eqn 2.2}) deduce the following identities whose combinatorial proofs are also provided.
\begin{theorem}\label{theom 2.3.2}
For any integer $n\geq 0$, there holds
\begin{eqnarray}
\sum_{i=0}^{n}(-2)^{i}D_{n,i} \hskip-.22cm &=&\hskip-.22cm 2^{n}, \label{eqn D2.1} \\
\sum_{i=0}^{n}(-1)^{i}D_{n,i} \hskip-.22cm &=&\hskip-.22cm \sum_{k=0}^{n}\binom{n}{k}C_k. \label{eqn D2.2}
\end{eqnarray}
\end{theorem}
\pf Let $\mathcal{D}_n^{e}$ ($\mathcal{D}_n^{o}$) denote the set of weighted G-Motzkin paths of length $n$ with even (odd) number of $\mathbf{d}$-steps such that each $\mathbf{d}$-step is weighted by $2$ (regarded as $\mathbf{d}_1$ and $\mathbf{d}_2$ for convenience) and other steps are weighted by $1$.
Let $\mathcal{D}_n^{*}$ be the subset of $\mathcal{D}_n^{e}$ such that each path in $\mathcal{D}_n^{*}$ has no $\mathbf{d}$-steps and only consists of $\mathbf{h}$-steps and $\mathbf{uv}$-peaks. Clearly,
$$w(\mathcal{D}_n^{e})=\sum_{i\ even}2^{i}D_{n,i}, \ \ w(\mathcal{D}_n^{o})=\sum_{i\ odd}2^{i}D_{n,i}\ \ \mbox{and} \ \ w(\mathcal{D}_n^{*})=2^n. $$
To prove (\ref{eqn D2.1}), it is sufficient to give a bijection $\tau$ between $\mathcal{D}_n^{e}/\mathcal{D}_n^{*}$ and $\mathcal{D}_n^{o}$. It is trivial for $n=0, 1$. For $n\geq 2$, any $\mathbf{P}\in \mathcal{D}_n^{e}/\mathcal{D}_n^{*}$ has at least one of the four subpaths, $\mathbf{d}_1$, $\mathbf{d}_2$, $\mathbf{uvv}$ and $\mathbf{hv}$, find the last one, say $\mathbf{z}$, $\mathbf{P}$ can be partitioned uniquely into $\mathbf{P}=\mathbf{P}_1\mathbf{z}\mathbf{P}_2$, where $\mathbf{P}_2=\mathbf{v}^k$ for certain $0\leq k<n$. Then define $\tau(\mathbf{P})=\mathbf{P}_1\mathbf{z}'\mathbf{P}_2$, where
\begin{eqnarray*}
\mathbf{z}'=\left\{
\begin{array}{rl}
 \mathbf{uvv}, &  \mbox{if}\ \mathbf{z}=\mathbf{d}_1, \\
 \mathbf{hv},  &  \mbox{if}\ \mathbf{z}=\mathbf{d}_2, \\
 \mathbf{d}_1, &  \mbox{if}\ \mathbf{z}=\mathbf{uvv}, \\
 \mathbf{d}_2, &  \mbox{if}\ \mathbf{z}=\mathbf{hv}.
\end{array}\right.
\end{eqnarray*}
This way ensures that the number of $\mathbf{d}$-steps in $\tau(\mathbf{P})$ is one more or less than that in $\mathbf{P}\in \mathcal{D}_n^{e}/\mathcal{D}_n^{*}$, so $\tau(\mathbf{P})\in \mathcal{D}_n^{o}$ and
$\mathbf{z}'$ in $\tau(\mathbf{P})$ is also the last one of the four subpaths, $\mathbf{d}_1$, $\mathbf{d}_2$, $\mathbf{uvv}$ and $\mathbf{hv}$.
The reverse procedure can be handled similarly. Hence, $\tau$ is a bijection (an involution in fact) between $\mathcal{D}_n^{e}/\mathcal{D}_n^{*}$ and $\mathcal{D}_n^{o}$. This completes the proof of (\ref{eqn D2.1}). \vskip0.2cm

Let $\mathcal{\bar{D}}_n^{e}$ ($\mathcal{\bar{D}}_n^{o}$) denote the set of weighted G-Motzkin paths of length $n$ with even (odd) number of $\mathbf{d}$-steps such that each step is weighted by $1$. Clearly,
$$w(\mathcal{\bar{D}}_n^{e})=\sum_{i\ even}D_{n,i} \ \mbox{and}\ w(\mathcal{\bar{D}}_n^{o})=\sum_{i\ odd}D_{n,i}. $$

Let $\mathcal{\bar{D}}_n^{*}$ be the subset of $\mathcal{\bar{D}}_n^{e}$ such that each path $\mathbf{Q}\in \mathcal{\bar{D}}_n^{*}$ has no $\mathbf{d}$-steps and no $\mathbf{hv}$-steps. Note that any $\mathbf{Q}\in \mathcal{\bar{D}}_n^{*}$ with $k$ $\mathbf{u}$-steps (with $k$ $\mathbf{v}$-steps naturally) and $n-k$ $\mathbf{h}$-steps can be obtained from Dyck paths $\mathbf{Q}'$ of length $2k$ for $0\leq k\leq n$ as follows. First replace each $\mathbf{d}$-step of $\mathbf{Q}'$ by a $\mathbf{v}$-step to get a G-Motzkin path $\mathbf{Q}''$ with no $\mathbf{h}$-steps and no $\mathbf{d}$-steps, and there are $\big(\binom{k+1}{n-k}\big)=\binom{n}{k}$ ways to insert $n-k$ $\mathbf{h}$-steps repeatedly into the $k+1$ positions exactly before $k$ $\mathbf{u}$-steps and at the endpoint of the path $\mathbf{Q}''$ to get $\mathbf{Q}$. This way can not produce $\mathbf{d}$-steps and $\mathbf{hv}$-steps in $\mathbf{Q}$. Summarizing over $k$, one has
$$ w(\mathcal{\bar{D}}_n^{*})=\sum_{k=0}^{n}\binom{n}{k}C_k. $$

To prove (\ref{eqn D2.2}), it is sufficient to give a bijection $\bar{\tau}$ between $\mathcal{\bar{D}}_n^{e}/\mathcal{\bar{D}}_n^{*}$ and $\mathcal{\bar{D}}_n^{o}$. It is trivial for $n=0$. For $n\geq 1$, any $\mathbf{Q}\in \mathcal{\bar{D}}_n^{e}/\mathcal{\bar{D}}_n^{*}$ has at least one of the two subpaths, $\mathbf{d}$ and $\mathbf{hv}$, find the last one, say $\mathbf{z}$, $\mathbf{Q}$ can be partitioned uniquely into $\mathbf{Q}=\mathbf{Q}_1\mathbf{z}\mathbf{Q}_2$. Then define $\bar{\tau}(\mathbf{Q})=\mathbf{Q}_1\mathbf{z}'\mathbf{Q}_2$, where
\begin{eqnarray*}
\mathbf{z}'=\left\{
\begin{array}{rl}
 \mathbf{hv},  &  \mbox{if}\ \mathbf{z}=\mathbf{d}, \\
 \mathbf{d},   &  \mbox{if}\ \mathbf{z}=\mathbf{hv}.
\end{array}\right.
\end{eqnarray*}
This way ensures that the number of $\mathbf{d}$-steps in $\bar{\tau}(\mathbf{Q})$ is one more or less than that in $\mathbf{Q}\in \mathcal{\bar{D}}_n^{e}/\mathcal{\bar{D}}_n^{*}$, so $\bar{\tau}(\mathbf{Q})\in \mathcal{\bar{D}}_n^{o}$ and
$\mathbf{z}'$ in $\bar{\tau}(\mathbf{Q})$ is also the last one of the two subpaths, $\mathbf{d}$ and $\mathbf{hv}$.
The reverse procedure can be handled similarly. Hence, $\bar{\tau}$ is a bijection (an involution in fact) between $\mathcal{\bar{D}}_n^{e}/\mathcal{\bar{D}}_n^{*}$ and $\mathcal{\bar{D}}_n^{o}$. This completes the proof (\ref{eqn D2.2}).    \qed\vskip0.2cm

\begin{theorem}\label{theom 2.3.4}
For any integer $n\geq 0$, there hold
\begin{eqnarray}
\sum_{i=0}^{n}(-1)^{i}D_{n+i,i} \hskip-.22cm &=&\hskip-.22cm 1, \label{eqn D2.3}\\
\sum_{i=0}^{n}(-2)^{i}D_{n+i,i} \hskip-.22cm &=&\hskip-.22cm \left\{
\begin{array}{ccc}
 1, &  \mbox{if}\ n=0, \\
 0, &  \mbox{otherwise}.
\end{array}\right.\label{eqn D2.4}
\end{eqnarray}
\end{theorem}
\pf Let $\mathcal{D}_{n+i, i}$ denote the set of weighted G-Motzkin paths of length $n+i$ with $i$ $\mathbf{d}$-steps such that all steps are weighted by $1$. Set
$$\mathcal{A}_n^{e}=\bigcup_{i=0, i\ even}^n\mathcal{D}_{n+i, i},\  \mathcal{A}_n^{o}=\bigcup_{i=0, i\ odd}^n\mathcal{D}_{n+i, i}.$$
Clearly,
$$w(\mathcal{A}_n^{e})=\sum_{i\ even}D_{n+i,i} \ \mbox{and} \  w(\mathcal{A}_n^{o})=\sum_{i\ odd}D_{n+i,i}. $$
To prove (\ref{eqn D2.3}), it is sufficient to give a bijection $\varphi$ between $\mathcal{A}_n^{e}/\{\mathbf{h}^{n}\}$ and $\mathcal{A}_n^{o}$. It is trivial for $n=0$. For $n\geq 1$, any $\mathbf{P}\in \mathcal{A}_n^{e}/\{\mathbf{h}^{n}\}$ has at least a $\mathbf{u}$-step, so there exist $\mathbf{d}$-steps or $\mathbf{v}$-steps in $\mathbf{P}$. Find the last return step $\mathbf{z}$, $\mathbf{P}$ can be partitioned uniquely into $\mathbf{P}=\mathbf{P}_1\mathbf{z}\mathbf{h}^k$ for certain $0\leq k<n$. Then define $\varphi(\mathbf{P})$ as follows:
\begin{eqnarray*}
\varphi(\mathbf{P})=\left\{
\begin{array}{rl}
 \mathbf{P}_1\mathbf{d}\mathbf{h}^k,   &  \mbox{if}\ \mathbf{z}=\mathbf{v}, \\
 \mathbf{P}_1\mathbf{v}\mathbf{h}^k,   &  \mbox{if}\ \mathbf{z}=\mathbf{d}.
\end{array}\right.
\end{eqnarray*}
This way ensures that the number of $\mathbf{d}$-steps in $\varphi(\mathbf{P})$ is one more or less than that in $\mathbf{P}\in \mathcal{A}_n^{e}/\{\mathbf{h}^{n}\}$, so $\varphi(\mathbf{P})\in \mathcal{A}_n^{o}$.
The reverse procedure can be handled similarly. Hence, $\varphi$ is a bijection (an involution in fact) between $\mathcal{A}_n^{e}/\{\mathbf{h}^{n}\}$ and $\mathcal{A}_n^{o}$. This completes the proof of (\ref{eqn D2.3}).

Let $\mathcal{\bar{D}}_{n+i, i}$ denote the set of weighted G-Motzkin paths of length $n+i$ with $i$ $\mathbf{d}$-steps such that each $\mathbf{d}$-step is weighted by $2$ (regarded as $\mathbf{d}_1$ and $\mathbf{d}_2$ for convenience) and other steps are weighted by $1$. Set
$$\mathcal{\bar{A}}_n^{e}=\bigcup_{i=0, i\ even}^n\mathcal{\bar{D}}_{n+i, i},\  \mathcal{\bar{A}}_n^{o}=\bigcup_{i=0, i\ odd}^n\mathcal{\bar{D}}_{n+i, i}.$$
Clearly,
$$w(\mathcal{\bar{A}}_n^{e})=\sum_{i\ even}2^iD_{n+i,i} \ \mbox{and} \  w(\mathcal{\bar{A}}_n^{o})=\sum_{i\ odd}2^iD_{n+i,i}. $$
It is trivial for $n=0$ in (\ref{eqn D2.4}). To prove (\ref{eqn D2.4}), it is sufficient to give a bijection $\bar{\varphi}$ between $\mathcal{\bar{A}}_n^{e}$ and $\mathcal{\bar{A}}_n^{o}$ for $n\geq 1$. Note that any $\mathbf{Q}\in \mathcal{\bar{A}}_n^{e}$ for $n\geq 1$ has at least one of the four subpaths, $\mathbf{h}$, $\mathbf{ud}_1$, $\mathbf{d}_2$ and $\mathbf{v}$, find the last one, say $\mathbf{z}$, $\mathbf{Q}$ can be partitioned uniquely into $\mathbf{Q}=\mathbf{Q}_1\mathbf{z}\mathbf{d}_1^k$ for certain $0\leq k<n$. Then define $\bar{\varphi}(\mathbf{Q})=\mathbf{Q}_1\mathbf{z}'\mathbf{d}_1^k$, where
\begin{eqnarray*}
\mathbf{z}'=\left\{
\begin{array}{rl}
 \mathbf{ud}_1,   &  \mbox{if}\ \mathbf{z}=\mathbf{h}, \\
 \mathbf{h},      &  \mbox{if}\ \mathbf{z}=\mathbf{ud}_1, \\
 \mathbf{v},      &  \mbox{if}\ \mathbf{z}=\mathbf{d}_2, \\
 \mathbf{d}_2,    &  \mbox{if}\ \mathbf{z}=\mathbf{v}.
\end{array}\right.
\end{eqnarray*}
This way ensures that the number of $\mathbf{d}$-steps in $\bar{\varphi}(\mathbf{Q})$ is one more or less than that in $\mathbf{Q}\in \mathcal{\bar{A}}_n^{e}$, so $\bar{\varphi}(\mathbf{Q})\in \mathcal{\bar{A}}_n^{o}$ and
$\mathbf{z}'$ in $\bar{\varphi}(\mathbf{Q})$ is also the last one of the four subpaths, $\mathbf{h}$, $\mathbf{ud}_1$, $\mathbf{d}_2$ and $\mathbf{v}$.
The reverse procedure can be handled similarly. Hence, $\bar{\varphi}$ is a bijection (an involution in fact) between $\mathcal{\bar{A}}_n^{e}$ and $\mathcal{\bar{A}}_n^{o}$. This completes the proof of (\ref{eqn D2.4}).          \qed\vskip0.2cm

\begin{theorem}\label{theom 2.3.5}
For any integer $n\geq 0$, there holds
\begin{eqnarray}
\sum_{i=0}^{n}y^{i}D_{n+i,i} \hskip-.22cm &=&\hskip-.22cm (y+1)^nN_n\Big(\frac{y+2}{y+1}\Big), \label{eqn D2.5}
\end{eqnarray}
where $N_n(y)=\sum_{k=1}^{n}\frac{1}{n}\binom{n}{k-1}\binom{n}{k}y^{k}=y^{n+1}N_n(y^{-1})$ with $N_0(y)=1$ is the Narayana polynomial \cite{ManSun}.
\end{theorem}
\pf Let $\mathcal{\hat{D}}_{n+i, i}$ denote the set of weighted G-Motzkin paths of length $n+i$ with $i$ $\mathbf{d}$-steps such that each $\mathbf{d}$-step is weighted by $y$ (regarded as $\mathbf{d}_y$ for convenience) and other steps are weighted by $1$. Let $\mathcal{\hat{C}}_{n, k}$ denote the set of weighted Dyck paths of length $2n$ with $k$ $\mathbf{ud}$-peaks such that each $\mathbf{d}$-step in $\mathbf{ud}$-peak is weighted by $y+2$ (regarded as $\mathbf{d}_y$, $\mathbf{d}_1$ and $\mathbf{d}_2$ for convenience) and other $\mathbf{d}$-steps are weighted by $y+1$ (regarded as $\mathbf{d}_y$ and $\mathbf{d}_1$ for convenience).  Set $\mathcal{\hat{D}}_{n}=\bigcup_{i=0}^n\mathcal{\hat{D}}_{n+i, i}$ and $\mathcal{\hat{C}}_{n}=\bigcup_{k=0}^n\mathcal{\hat{C}}_{n, k}$.
Since the Narayana number $N_{n,k}=\frac{1}{n}\binom{n}{k-1}\binom{n}{k}$ counts the number of Dyck paths of length $2n$ with $k$ $\mathbf{ud}$-peaks \cite{Deutsch99}, it is clear that
$$|\mathcal{\hat{C}}_{n, k}|=N_{n,k}(y+2)^{k}(y+1)^{n-k}. $$

So it is sufficient to give a bijection $\hat{\varphi}$ between $\mathcal{\hat{C}}_n$ and $\mathcal{\hat{D}}_n$ for $n\geq 1$. For any $\mathbf{Q}\in \mathcal{\hat{C}}_{n}$ for $n\geq 1$ with $k_0$ $\mathbf{d}_y$-steps, $k_1$ $\mathbf{d}_1$-steps and $k_2$ $\mathbf{d}_2$-steps, note that $k_0+k_1+k_2=n$ and each $\mathbf{d}_2$-step must in a $\mathbf{ud}_2$-peak, replace each $\mathbf{d}_1$-step by a $\mathbf{v}$-step and each $\mathbf{ud}_2$-peak by an $\mathbf{h}$-step, we obtain a weighted G-Motzkin path $\mathbf{Q}'\in \mathcal{\hat{D}}_{n+k_0, k_0}$. Then define $\hat{\varphi}(\mathbf{Q})=\mathbf{Q}'.$ It is not difficult to verify that
$\hat{\varphi}$ is a bijection $\mathcal{\hat{C}}_n$ and $\mathcal{\hat{D}}_n$. This completes the proof of (\ref{eqn D2.5}).          \qed\vskip0.2cm

In order to give a more intuitive view on the bijection $\hat{\varphi}$, a pictorial description of
$\hat{\varphi}$ is presented for $\mathbf{Q}=\mathbf{u}\mathbf{d}_y\mathbf{uu}\mathbf{d}_2
\mathbf{uuu}\mathbf{d}_y\mathbf{d}_1\mathbf{d}_y\mathbf{uuuu}\mathbf{d}_2\mathbf{u}\mathbf{d}_1\mathbf{d}_y\mathbf{d}_y\mathbf{d}_1\mathbf{d}_1\mathbf{u}\mathbf{d}_2\mathbf{uuu}\mathbf{d}_y\mathbf{d}_1\mathbf{d}_1$, we have
$$\hat{\varphi}(\mathbf{Q})=\mathbf{u}\mathbf{d}_y\mathbf{u}\mathbf{h}\mathbf{uuu}\mathbf{d}_y\mathbf{v}\mathbf{d}_y\mathbf{uuu}\mathbf{h}
\mathbf{u}\mathbf{v}\mathbf{d}_y\mathbf{d}_y\mathbf{v}\mathbf{v}\mathbf{h}\mathbf{uuu}\mathbf{d}_y\mathbf{v}\mathbf{v}. $$
See Figure 2 for detailed illustrations.

\begin{figure}[h] \setlength{\unitlength}{0.5mm}

\begin{center}
\begin{pspicture}(18,3.5)
\psset{xunit=15pt,yunit=15pt}\psgrid[subgriddiv=1,griddots=4,
gridlabels=4pt](0,0)(30,7)

\psline(0,0)(1,1)(2,0)(4,2)(5,1)(8,4)(11,1)(14,4)(15,5)(16,4)(17,5)(22,0)(23,1)(24,0)(27,3)(30,0)

\pscircle*(0,0){0.06}\pscircle*(1,1){0.06}\pscircle*(2,0){0.06}
\pscircle*(3,1){0.06}\pscircle*(4,2){0.06}\pscircle*(5,1){0.06}
\pscircle*(6,2){0.06}\pscircle*(7,3){0.06}\pscircle*(8,4){0.06}
\pscircle*(9,3){0.06}\pscircle*(10,2){0.06}\pscircle*(11,1){0.06}
\pscircle*(12,2){0.06}\pscircle*(13,3){0.06}\pscircle*(14,4){0.06}\pscircle*(15,5){0.06}
\pscircle*(16,4){0.06}\pscircle*(17,5){0.06}\pscircle*(18,4){0.06}
\pscircle*(19,3){0.06}\pscircle*(20,2){0.06}\pscircle*(21,1){0.06}\pscircle*(22,0){0.06}
\pscircle*(23,1){0.06}\pscircle*(24,0){0.06}\pscircle*(25,1){0.06}
\pscircle*(26,2){0.06}\pscircle*(27,3){0.06}\pscircle*(28,2){0.06}\pscircle*(29,1){0.06}\pscircle*(30,0){0.06}

\put(2.7,3.3){$\mathbf{Q}=\mathbf{u}\mathbf{d}_y\mathbf{uu}\mathbf{d}_2
\mathbf{uuu}\mathbf{d}_y\mathbf{d}_1\mathbf{d}_y\mathbf{uuuu}\mathbf{d}_2
\mathbf{u}\mathbf{d}_1\mathbf{d}_y\mathbf{d}_y\mathbf{d}_1\mathbf{d}_1\mathbf{u}\mathbf{d}_2\mathbf{uuu}\mathbf{d}_y\mathbf{d}_1\mathbf{d}_1$}

\put(.75,.35){$\mathbf{d}_y$}\put(2.3,.9){$\mathbf{d}_2$}\put(4.45,1.9){$\mathbf{d}_y$}\put(5,1.4){$\mathbf{d}_1$}\put(5.45,.9){$\mathbf{d}_y$}
\put(8.2,2.4){$\mathbf{d}_2$}\put(9.25,2.4){$\mathbf{d}_1$}\put(9.75,1.9){$\mathbf{d}_y$}\put(10.25,1.4){$\mathbf{d}_y$}\put(10.8,.85){$\mathbf{d}_1$}\put(11.3,.3){$\mathbf{d}_1$}
\put(12.4,.3){$\mathbf{d}_2$}\put(14.45,1.4){$\mathbf{d}_y$}\put(15,.85){$\mathbf{d}_1$}\put(15.5,.3){$\mathbf{d}_1$}

\end{pspicture}
\end{center}\vskip0.5cm

$\Updownarrow \hat{\varphi}$
\vskip0.5cm

\begin{center}
\begin{pspicture}(12,3.6)
\psset{xunit=15pt,yunit=15pt}\psgrid[subgriddiv=1,griddots=4,
gridlabels=4pt](0,0)(22,7)

\psline(0,0)(1,1)(2,0)(3,1)(4,1)(7,4)(8,3)(8,2)(9,1)(12,4)(13,4)(14,5)(14,4)(16,2)(16,0)(17,0)(20,3)(21,2)(21,0)

\pscircle*(0,0){0.06}\pscircle*(1,1){0.06}\pscircle*(2,0){0.06}
\pscircle*(3,1){0.06}\pscircle*(4,1){0.06}\pscircle*(5,2){0.06}\pscircle*(6,3){0.06}
\pscircle*(7,4){0.06} \pscircle*(8,3){0.06}\pscircle*(8,2){0.06}
\pscircle*(9,1){0.06}\pscircle*(10,2){0.06}\pscircle*(11,3){0.06}\pscircle*(12,4){0.06}
\pscircle*(13,4){0.06}\pscircle*(14,5){0.06}\pscircle*(14,4){0.06}
\pscircle*(15,3){0.06}\pscircle*(16,2){0.06}\pscircle*(16,1){0.06}\pscircle*(16,0){0.06}
\pscircle*(17,0){0.06}\pscircle*(18,1){0.06}\pscircle*(19,2){0.06}\pscircle*(20,3){0.06}
\pscircle*(21,2){0.06}\pscircle*(21,1){0.06}\pscircle*(21,0){0.06}

\put(.75,.35){$\mathbf{d}_y$}\put(4.,1.9){$\mathbf{d}_y$}\put(4.5,.85){$\mathbf{d}_y$}
\put(7.75,1.85){$\mathbf{d}_y$}\put(8.25,1.35){$\mathbf{d}_y$}\put(10.85,1.35){$\mathbf{d}_y$}

\put(1.2,3.3){$\hat{\varphi}(\mathbf{Q})=\mathbf{u}\mathbf{d}_y\mathbf{u}\mathbf{h}
\mathbf{uuu}\mathbf{d}_y\mathbf{v}\mathbf{d}_y\mathbf{uuu}\mathbf{h}
\mathbf{u}\mathbf{v}\mathbf{d}_y\mathbf{d}_y\mathbf{v}\mathbf{v}\mathbf{h}\mathbf{uuu}\mathbf{d}_y\mathbf{v}\mathbf{v}$}

\end{pspicture}
\end{center}\vskip0.5cm

\caption{\small An example of the bijection $\hat{\varphi}$ described in the proof of Theorem \ref{theom 2.3.5}. }

\end{figure}
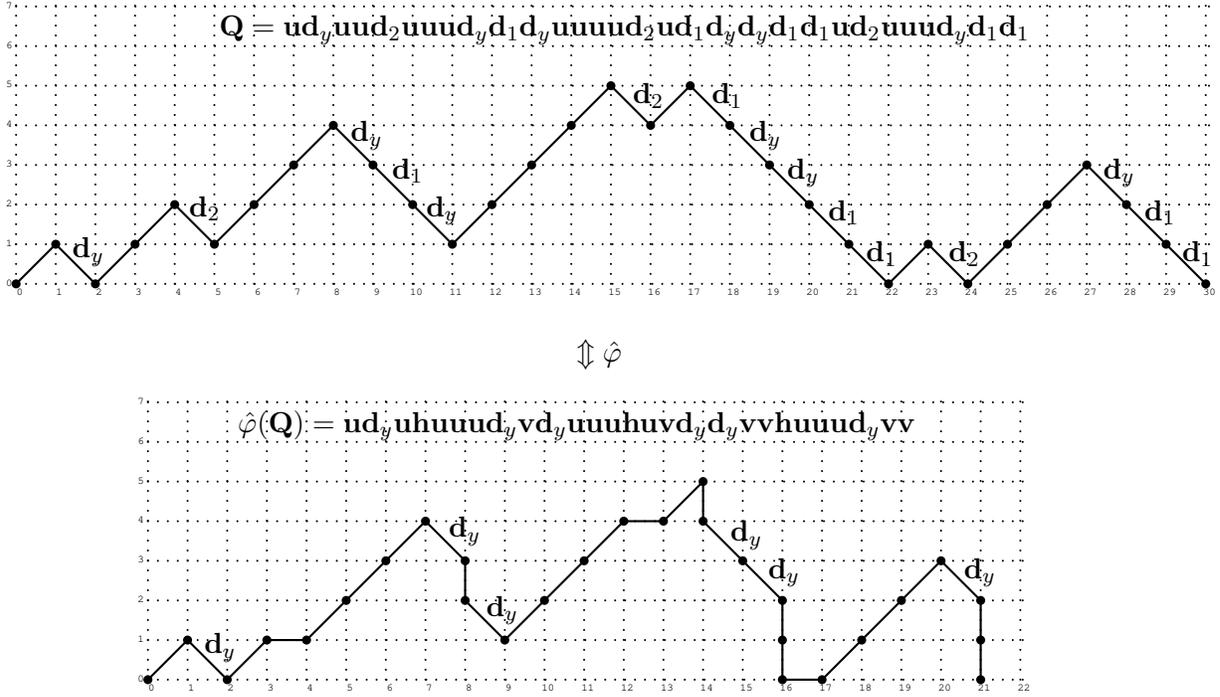

When $y=-3$ in (\ref{eqn D2.5}), by the well-known relation $R_n=N_n(2)=2r_n$ for $n\geq 1$, one can derive that the following identity which is asked for
a direct combinatorial proof similar to that of Theorem \ref{theom 2.3.4}.
\begin{corollary}\label{coro 2.3.6}
For any integer $n\geq 0$, there holds
\begin{eqnarray*}
\sum_{i=0}^{n}(-3)^{i}D_{n+i,i} \hskip-.22cm &=&\hskip-.22cm (-1)^nr_n,
\end{eqnarray*}
where $r_n$ is the little Schr\"{o}der number.
\end{corollary}

\subsection{ The statistics ``number of $\mathbf{u}$-steps" }
Let $U_{n,i}$ denote the number of G-Motzkin paths of length $n$ with $i$ $\mathbf{u}$-steps, the first values of $U_{n,i}$ are illustrated in Table 2.4.

\begin{center}
\begin{eqnarray*}
\begin{array}{c|ccccccc}\hline
n/i & 0   & 1   & 2    & 3    & 4    & 5    & 6     \\\hline
  0 & 1   &     &      &      &      &      &      \\
  1 & 1   & 1   &      &      &      &      &      \\
  2 & 1   & 4   & 2    &      &      &      &       \\
  3 & 1   & 9   & 14   & 5    &      &      &       \\
  4 & 1   & 16  & 52   & 50   &  14  &      &       \\
  5 & 1   & 25  & 140  & 260  & 182  & 42   &        \\
  6 & 1   & 36  & 310  & 950  & 1218 & 462  & 132     \\\hline
\end{array}
\end{eqnarray*}
Table 2.4. The first values of $U_{n,i}$.
\end{center}

Any Dyck paths $\mathbf{P}$ of length $2i$ can be extended to G-Motzkin paths of length $n$ with $i$ $\mathbf{u}$-steps for $0\leq i\leq n$. Note that for $0\leq k\leq i$ there are $\binom{i}{k}$ ways to replace $k$ $\mathbf{d}$-steps by $\mathbf{v}$-steps and $\big(\binom{2i+1}{n-2i+k}\big)=\binom{n+k}{2i}$ ways to insert repeatedly $n-2i+k$ $\mathbf{h}$-steps into $2i+1$ points of $\mathbf{P}$ to form G-Motzkin paths of length $n$ with $i$ $\mathbf{u}$-steps. Summarizing over $k$, we have
\begin{theorem}\label{theom 2.4.1}
For any integers $n\geq i\geq 0$, there holds
\begin{eqnarray*}
U_{n,i}=\sum_{k=0}^{i}\binom{i}{k}\binom{n+k}{2i}C_i.
\end{eqnarray*}
\end{theorem}
Note that $U_{n,i}$ is also the coefficient of $b^i$ in $G_n(1, b, b)$ in (\ref{eqn 2.2b}) which has another expression
\begin{eqnarray*}
U_{n,i}=\frac{1}{n+1}\binom{n+1}{i+1}\sum_{j=0}^{i}\binom{n-i}{j}\binom{n+i-j}{i-j},
\end{eqnarray*}
and $U_{n,n}$ is the $n$-th Catalan number $C_n$. By (\ref{eqn 2.2}), the special case $G(1, -2, -2; x)=\frac{1}{1+x}$ deduces the following identity. We provide a combinatorial proof for this identity.

\begin{theorem}\label{theom 2.4.2}
For any integer $n\geq 0$, there holds
\begin{eqnarray*}
\sum_{i=0}^{n}(-2)^{i}U_{n,i}=(-1)^{n}.
\end{eqnarray*}
\end{theorem}

\pf Let $\mathcal{U}_n^{e}$ ($\mathcal{U}_n^{o}$) denote the set of weighted G-Motzkin paths of length $n$ with even (odd) number of $\mathbf{u}$-steps such that the matching step of each $\mathbf{u}$-step, that is each of $\mathbf{d}$-steps and $\mathbf{v}$-steps is weighted by $2$ (regarded respectively as $\mathbf{d}_1$ and $\mathbf{d}_2$, $\mathbf{v}_1$ and $\mathbf{v}_2$ for convenience) and other steps are weighted by $1$. Clearly,
$$w(\mathcal{U}_n^{e})=\sum_{i\ even}2^{i}U_{n,i}\ \ \mbox{and}\ \ w(\mathcal{U}_n^{o})=\sum_{i\ odd}2^{i}U_{n,i}. $$
So it is sufficient to give a bijection $\theta$ between $\mathcal{U}_n^{e}/\{(\mathbf{uv}_1)^n\}$ and $\mathcal{U}_n^{o}$ for $n$ even and between $\mathcal{U}_n^{e}$ and $\mathcal{U}_n^{o}/\{(\mathbf{uv}_1)^n\}$ for $n$ odd. When $n\geq 1$, any $\mathbf{P}\in \mathcal{U}_n^{e}/\{(\mathbf{uv}_1)^n\}$ for $n$ even or $\mathbf{P}\in \mathcal{U}_n^{e}$ for $n$ odd has at least one of the six subpaths, $\mathbf{h}$, $\mathbf{u}\mathbf{v}_2$, $\mathbf{d}_1$, $\mathbf{u}\mathbf{v}_1^2$, $\mathbf{d}_2$ and $\mathbf{u}\mathbf{v}_1\mathbf{v}_2$, find the last one, say $\mathbf{s}$, $\mathbf{P}$ can be partitioned uniquely into $\mathbf{P}=\mathbf{P}_1\mathbf{s}\mathbf{P}_2$, where $\mathbf{P}_2=\mathbf{t}_1\cdots \mathbf{t}_{k}$ and $\mathbf{t}_{j}\in \{\mathbf{v}_1, \mathbf{v}_2 \}$ for $0\leq j\leq k<n$. Then define $\theta(\mathbf{P})=\mathbf{P}_1\mathbf{s}'\mathbf{P}_2$, where
\begin{eqnarray*}
\mathbf{s}'=\left\{
\begin{array}{rl}
 \mathbf{u}\mathbf{v}_2,               &  \mbox{if}\ \mathbf{s}=\mathbf{h}, \\
 \mathbf{h},                           &  \mbox{if}\ \mathbf{s}=\mathbf{u}\mathbf{v}_2, \\
 \mathbf{u}\mathbf{v}_1^2,             &  \mbox{if}\ \mathbf{s}=\mathbf{d}_1, \\
 \mathbf{d}_1,                         &  \mbox{if}\ \mathbf{s}=\mathbf{u}\mathbf{v}_1^2, \\
 \mathbf{u}\mathbf{v}_1\mathbf{v}_2,   &  \mbox{if}\ \mathbf{s}=\mathbf{d}_2, \\
 \mathbf{d}_2,                         &  \mbox{if}\ \mathbf{s}=\mathbf{u}\mathbf{v}_1\mathbf{v}_2.
\end{array}\right.
\end{eqnarray*}
This way ensures that the number of $\mathbf{u}$-steps in $\theta(\mathbf{P})$ is one more or less than that in $\mathbf{P}\in \mathcal{U}_n^{e}/\{(\mathbf{uv}_1)^n\}$ for $n$ even or in $\mathbf{P}\in \mathcal{U}_n^{e}$ for $n$ odd, so $\theta(\mathbf{P})\in \mathcal{U}_n^{o}$ for $n$ even and $\theta(\mathbf{P})\in \mathcal{U}_n^{o}/\{(\mathbf{uv}_1)^n\}$ for $n$ odd. Moreover, $\mathbf{s}'$ is also the last one of the six subpaths, $\mathbf{h}$, $\mathbf{u}\mathbf{v}_2$, $\mathbf{d}_1$, $\mathbf{u}\mathbf{v}_1^2$, $\mathbf{d}_2$ and $\mathbf{u}\mathbf{v}_1\mathbf{v}_2$. The reverse procedure can be handled similarly.
Hence, $\theta$ is a bijection (an involution in fact) between $\mathcal{U}_n^{e}/\{(\mathbf{uv}_1)^n\}$ and $\mathcal{U}_n^{o}$ for $n$ even and between $\mathcal{U}_n^{e}$ and $\mathcal{U}_n^{o}/\{(\mathbf{uv}_1)^n\}$ for $n$ odd. This completes the proof.            \qed\vskip0.2cm

\section{The statistics ``number of $\mathbf{z}$-steps" at given level in G-Motzkin paths }

In the literature, there are several statistics considered at given level in Dyck paths, such as ``number of $\mathbf{u}$-steps", ``number of $\mathbf{ud}$-peaks" and ``number of points" \cite{Cheng}. Precisely, let $B_{n,i}$ be the number of $\mathbf{u}$-steps" at level $i+1$ in all Dyck paths in $\mathcal{C}_{n+1}$ and $C_{n,i}$ be the number of points at level $i$ in all Dyck paths in $\mathcal{C}_{n}$. It is known that $B_{n,i}=\frac{2i+3}{2n+3}\binom{2n+3}{n-i}$ is the bollat number and $C_{n,i}=\frac{i+1}{n+1}\binom{2n+1}{n-i}$
also enumerates the number of $\mathbf{ud}$-peaks at level $i+1$ in all Dyck paths in $\mathcal{C}_{n+1}$. The first values of $B_{n,i}$ and $C_{n,i}$ are illustrated respectively in Table 3.1 and in Table 3.2.

\begin{center}
\begin{eqnarray*}
\begin{array}{c|cccccccc}\hline
n/i & 0   & 1   & 2    & 3    & 4    & 5    & 6       \\\hline
  0 & 1   &     &      &      &      &      &        \\
  1 & 3   & 1   &      &      &      &      &        \\
  2 & 9   & 5   & 1    &      &      &      &         \\
  3 & 28  & 20  & 7    & 1    &      &      &         \\
  4 & 90  & 75  & 35   & 9    &  1   &      &         \\
  5 & 297 & 275 & 154  & 54   &  11  &  1   &        \\
  6 & 1001 &1001 &637  & 273  &  77  &  13  &  1     \\\hline
\end{array}
\end{eqnarray*}
Table 3.1. The first values of $B_{n,i}$.

\begin{eqnarray*}
\begin{array}{c|cccccccc}\hline
n/i & 0   & 1   & 2    & 3    & 4    & 5    & 6       \\\hline
  0 & 1   &     &      &      &      &      &        \\
  1 & 2   & 1   &      &      &      &      &        \\
  2 & 5   & 4   & 1    &      &      &      &         \\
  3 & 14  & 14  & 6    & 1    &      &      &         \\
  4 & 42  & 48  & 27   & 8    &  1   &      &         \\
  5 & 132 & 165 & 110  & 44   &  10  &  1   &        \\
  6 & 429 & 572 & 429  & 208  &  65  &  12  &  1     \\\hline
\end{array}
\end{eqnarray*}
Table 3.2. The first values of $C_{n,i}$.

\end{center}

Actually, the matrices $\big(C_{n,i}\big)_{n\geq i\geq 0}$ and $\big(B_{n,i}\big)_{n\geq i\geq 0}$ form Riordan arrays $(C(x)^2, xC(x)^2)$ and $(C(x)^3, xC(x)^2)$ respectively. Recall that a {\it Riordan array} \cite{ShapB, ShapGet, Sprug} is an infinite lower triangular matrix $\mathscr{D}=(d_{n,i})_{n,i \in \mathbb{N}}$ such that its $i$-th column has generating function $d(x)h(x)^i$, where $d(x)$ and $h(x)$ are formal power series with $d(0)=1$ and $h(0)=0$. That is, the general term of $\mathscr{D}$ is $d_{n,i}=[x^n]d(x)h(x)^i$, where $[x^n]$ is the coefficient operator. The matrix $\mathscr{D}$ corresponding to the pair $d(x)$ and $h(x)$ is denoted by $(d(x),h(x))$.
The product of a Riordan array $(d(x),h(x))$ and a formal power series $A(x)=\sum_{n\geq 0}a_nx^n$ is given by $(d(x),h(x))A(x)=d(x)A(h(x))$, this implies that
\begin{eqnarray}\label{eqn 3.1}
\sum_{i=0}^{n}d_{n,i}a_i=[x^n]d(x)A(h(x)).
\end{eqnarray}


In this section, we focus on the enumeration of statistics ``number of $\mathbf{z}$-steps" for $\mathbf{z}\in \{\mathbf{u}, \mathbf{h}, \mathbf{v}, \mathbf{d}\}$, ``number of return steps" and ``number of points" at given level in G-Motzkin paths. Some counting results are linked with Riordan arrays.

\subsection{ The statistics ``number of $\mathbf{u}$-steps" at level $i+1$ }
Let $\alpha_{n, i}$ denote the number of $\mathbf{u}$-steps at level $i+1$ in all G-Motzkin paths of length $n+1$, the first values of $\alpha_{n, i}$ are illustrated in Table 3.3.

\begin{center}
\begin{eqnarray*}
\begin{array}{c|ccccccc}\hline
n/i & 0      & 1      & 2       & 3       & 4    & 5    & 6    \\\hline
  0 & 1      &        &         &         &      &      &     \\
  1 & 7      & 1      &         &         &      &      &     \\
  2 & 39     & 12     & 1       &         &      &      &      \\
  3 & 212    & 96     & 17      & 1       &      &      &      \\
  4 & 1157   & 665    & 178     & 22      & 1    &      &      \\
  5 & 6384   & 4320   & 1513    & 285     & 27   & 1    &      \\
  6 & 35647  & 27177  & 11522   & 2881    & 417  & 32   & 1      \\\hline
\end{array}
\end{eqnarray*}
Table 3.3. The first values of $\alpha_{n, i}$.
\end{center}

\begin{theorem}\label{theom 3.1.1}
For any integers $n\geq i\geq 0$, there holds
\begin{eqnarray*}
\alpha_{n,i}=\sum_{j=i}^{n}B_{j,i}\sum_{k=0}^{n-j}\binom{j+1}{k}\binom{n+j-k+2}{n-j-k},
\end{eqnarray*}
Moreover, $\alpha_{n,i}$ is the $(n,i)$-entry of the
Riordan array
$$\Big(\frac{1+x}{(1-x)^3}C\Big(\frac{x(1+x)}{(1-x)^2}\Big)^3, \frac{x(1+x)}{(1-x)^2}C\Big(\frac{x(1+x)}{(1-x)^2}\Big)^2\Big). $$
\end{theorem}
\pf For each Dyck path $\mathbf{P}$ of length $2j+2$, it can be extended to G-Motzkin paths $\mathbf{Q}$ of length $n+1$ for $i\leq j\leq n$ such that $\mathbf{P}$ and $\mathbf{Q}$ have the same number of $\mathbf{u}$-steps at level $i+1$. First replace $j+1-k$ $\mathbf{d}$-steps in $\mathbf{P}$ by $\mathbf{v}$-steps to get $\mathbf{P}'$, there are $\binom{j+1}{k}$ ways, and insert repeatedly $n-j-k$ $\mathbf{h}$-steps into $2j+3$ points of $\mathbf{P}'$ to form G-Motzkin paths $\mathbf{Q}$ of length $n+1$, there $\big(\binom{2j+3}{n-j-k}\big)=\binom{n+j-k+2}{n-j-k}$ ways. Note that there are totally $B_{j, i}$ $\mathbf{u}$-steps at level $i+1$ in all $\mathbf{P}\in \mathcal{C}_{j+1}$, summarizing over $k$ and $j$, we obtain the desired result.


On the other hand, for any G-Motzkin path $\mathbf{G}\in \mathcal{G}$ with at least one $\mathbf{u}$-step at level $i+1$, given such a $\mathbf{u}$-step, marked as $\mathbf{u}^{*}$, $\mathbf{G}$ can be partitioned uniquely into
$$\mathbf{G}=\mathbf{G}_0\mathbf{u}\mathbf{G}_1\dots \mathbf{u}\mathbf{G}_{i}\mathbf{u}^{*}\mathbf{G}_{i+1}\mathbf{z}_1\mathbf{\bar{G}}_{1}\mathbf{z}_2\mathbf{\bar{G}}_{2}\dots \mathbf{z}_{i+1}\mathbf{\bar{G}}_{i+1},$$
where $\mathbf{G}_0, \dots, \mathbf{G}_{i+1}, \mathbf{\bar{G}}_0, \dots, \mathbf{\bar{G}}_{i+1} \in \mathcal{G}$ and $\mathbf{z}_1, \dots, \mathbf{z}_{i+1} \in \{\mathbf{d}, \mathbf{v}\}$. Since each of $\mathbf{G}_k$ and $\mathbf{\bar{G}}_k$ has the generating function $G(x)$, each of $\mathbf{u}$ and $\mathbf{d}$ produces an $x$ and each $\mathbf{v}$ leads to a $1$, this makes $\mathbf{z}_1\mathbf{z}_2\dots \mathbf{z}_{i+1}$ generate $(1+x)^{i+1}$, so according to the length of $\mathbf{G}$, all $\mathbf{G}\in \mathcal{G}$ produce the generating function
$x^{i+1}(1+x)^{i+1}G(x)^{2i+3}$. Hence, the total number $\alpha_{n,i}$ of $\mathbf{u}^{*}$-steps in all G-Motzkin paths $\mathbf{G}\in \mathcal{G}_{n+1}$ is the coefficient of $x^{n+1}$ in $x^{i+1}(1+x)^{i+1}G(x)^{2i+3}$, namely,
$$\alpha_{n,i}=[x^{n+1}]x^{i+1}(1+x)^{i+1}G(x)^{2i+3}=[x^n](1+x)G(x)^3\Big(x(1+x)G(x)^2\Big)^{i}. $$
By (\ref{eqn 2.3}), $\alpha_{n,i}$ is the $(n,i)$-entry of the Riordan array
$$\Big((1+x)G(x)^3, x(1+x)G(x)^2\Big)=\Big(\frac{1+x}{(1-x)^3}C\Big(\frac{x(1+x)}{(1-x)^2}\Big)^3, \frac{x(1+x)}{(1-x)^2}C\Big(\frac{x(1+x)}{(1-x)^2}\Big)^2\Big). $$
This completes the proof of Theorem 3.1.        \qed\vskip0.2cm

\begin{theorem}\label{theom 3.1.2}
For any integers $n, m\geq 0$, there holds
\begin{eqnarray}\label{eqn 3.1.2}
\sum_{i=0}^{n}(-1)^{n-i}\binom{n}{i}\alpha_{n+m+i,m+i}=5^{n}.
\end{eqnarray}
\end{theorem}
\pf By Theorem \ref{theom 3.1.1}, we have
\begin{eqnarray*}
\lefteqn{\sum_{i=0}^{n}(-1)^{n-i}\binom{n}{i}\alpha_{n+m+i,m+i} } \\
\hskip-.22cm &=&\hskip-.22cm \sum_{i=0}^{n}(-1)^{n-i}\binom{n}{i}[x^{n+m+i}]\frac{1+x}{(1-x)^3}C\Big(\frac{x(1+x)}{(1-x)^2}\Big)^3 \Big(\frac{x(1+x)}{(1-x)^2}C\Big(\frac{x(1+x)}{(1-x)^2}\Big)^2\Big)^{m+i} \\
\hskip-.22cm &=&\hskip-.22cm [x^n]\frac{(1+x)^{m+1}}{(1-x)^{2m+3}}C\Big(\frac{x(1+x)}{(1-x)^2}\Big)^{2m+3}\sum_{i=0}^{n}(-1)^{n-i}\binom{n}{i} \Big(\frac{1+x}{(1-x)^2}C\Big(\frac{x(1+x)}{(1-x)^2}\Big)^2\Big)^{i} \\
\hskip-.22cm &=&\hskip-.22cm [x^n]\frac{(1+x)^{m+1}}{(1-x)^{2m+3}}C\Big(\frac{x(1+x)}{(1-x)^2}\Big)^{2m+3}\Big(\frac{1+x}{(1-x)^2}C\Big(\frac{x(1+x)}{(1-x)^2}\Big)^2-1\Big)^{n} \\
\hskip-.22cm &=&\hskip-.22cm [x^n]\frac{(1+x)^{m+1}}{(1-x)^{2m+3}}C\Big(\frac{x(1+x)}{(1-x)^2}\Big)^{2m+3}\Big(\frac{1+x}{(1-x)^2}-1+\frac{2x(1+x)^2}{(1-x)^4}+ \sum_{k=2}^{\infty}C_k\frac{x^k(1+x)^{k+1}}{(1-x)^{2k+2}}\Big)^{n} \\
\hskip-.22cm &=&\hskip-.22cm [x^n]\frac{(1+x)^{m+1}}{(1-x)^{2m+3}}C\Big(\frac{x(1+x)}{(1-x)^2}\Big)^{2m+3}\Big(\frac{x(5-x+3x^2+x^3)}{(1-x)^4}+ \sum_{k=2}^{\infty}C_k\frac{x^k(1+x)^{k+1}}{(1-x)^{2k+2}}\Big)^{n} \\
\hskip-.22cm &=&\hskip-.22cm [x^0]\frac{(1+x)^{m+1}}{(1-x)^{2m+3}}C\Big(\frac{x(1+x)}{(1-x)^2}\Big)^{2m+3}\Big(\frac{5-x+3x^2+x^3}{(1-x)^4}+ \sum_{k=2}^{\infty}C_k\frac{x^{k-1}(1+x)^{k+1}}{(1-x)^{2k+2}}\Big)^{n} \\
\hskip-.22cm &=&\hskip-.22cm  \frac{(1+x)^{m+1}}{(1-x)^{2m+3}}C\Big(\frac{x(1+x)}{(1-x)^2}\Big)^{2m+3}\Big(\frac{5-x+3x^2+x^3}{(1-x)^4}+ \sum_{k=2}^{\infty}C_k\frac{x^{k-1}(1+x)^{k+1}}{(1-x)^{2k+2}}\Big)^{n}\Big|_{x=0} \\
\hskip-.22cm &=&\hskip-.22cm 5^n.
\end{eqnarray*}
This completes the proof.  \qed\vskip0.2cm

\begin{theorem}\label{theom 3.1.3}
For any integer $n\geq 0$, there holds
\begin{eqnarray}\label{eqn 3.1.3}
\sum_{i=0}^{n}(-1)^{i}\binom{i+2}{2}\alpha_{n,i}=(n+1)^{2}.
\end{eqnarray}
\end{theorem}
\pf By (\ref{eqn 3.1}) and Theorem \ref{theom 3.1.2}, together with the relation $C(x)=1+xC(x)^2$, we have
\begin{eqnarray*}
\lefteqn{\sum_{i=0}^{n}(-1)^{i}\binom{i+2}{2}\alpha_{n,i}} \\
\hskip-.22cm &=&\hskip-.22cm
[x^n]\Big(\frac{1+x}{(1-x)^3}C\Big(\frac{x(1+x)}{(1-x)^2}\Big)^3, \frac{x(1+x)}{(1-x)^2}C\Big(\frac{x(1+x)}{(1-x)^2}\Big)^2\Big)\frac{1}{(1+x)^3}  \\
\hskip-.22cm &=&\hskip-.22cm  [x^n]\frac{1+x}{(1-x)^3} =(n+1)^2.
\end{eqnarray*}
This completes the proof.  \qed\vskip0.2cm
One can be asked for combinatorial proofs of these two identities (\ref{eqn 3.1.2}) and (\ref{eqn 3.1.3}).

\subsection{ The statistics ``number of $\mathbf{v}$-steps" and ``number of $\mathbf{d}$-steps" at level $i$ }
Let $\beta_{n, i}$ denote the number of $\mathbf{v}$-steps at level $i$ in all G-Motzkin paths of length $n+1$ and let $\gamma_{n, i}$ denote the number of $\mathbf{d}$-steps at level $i$ in all G-Motzkin paths of length $n+2$, the first values of $\beta_{n, i}$ and $\gamma_{n, i}$ are illustrated in Table 3.4.

\begin{center}
\begin{eqnarray*}
\begin{array}{c|ccccccc}\hline
n/i & 0     & 1     & 2      & 3      & 4      & 5    & 6   \\\hline
  0 & 1     &       &        &        &        &      &     \\
  1 & 6     & 1     &        &        &        &      &     \\
  2 & 33    & 11    & 1      &        &        &      &      \\
  3 & 179   & 85    & 16     & 1      &        &      &      \\
  4 & 978   & 580   & 162    & 21     & 1      &      &      \\
  5 & 5406  & 3740  & 1351   & 264    & 26     & 1    &      \\
  6 & 30241 & 23437 & 10171  & 2617   & 391    & 31  & 1      \\\hline
\end{array}
\end{eqnarray*}
Table 3.4. The first values of $\beta_{n, i}$ and $\gamma_{n, i}$.
\end{center}

\begin{lemma}\label{lemma 3.2.1}
For any integers $n\geq i\geq 0$, there holds
\begin{eqnarray*}
\beta_{n, i}=\gamma_{n,i}.
\end{eqnarray*}
\end{lemma}
\pf Given a $\mathbf{v}$-step counted at level $i$ in a G-Motzkin path $\mathbf{P}$ of length $n+1$, replace it by a $\mathbf{d}$-step, one get a G-Motzkin path $\mathbf{P}'$ of length $n+2$ with a $\mathbf{d}$-step counted at level $i$ and vice versa. This implies that $\beta_{n, i}=\gamma_{n,i}$.  \qed\vskip0.2cm

\begin{theorem}\label{theom 3.2.2}
For any integers $n\geq i\geq 0$, there holds
\begin{eqnarray*}
\beta_{n, i}=\gamma_{n,i}=\sum_{j=i}^{n}B_{j,i}\sum_{k=0}^{j}\binom{j}{k}\binom{n+j+2-k}{n-j-k}.
\end{eqnarray*}
Moreover, $\beta_{n,i}$ is the $(n,i)$-entry of the
Riordan array
$$\Big(\frac{1}{(1-x)^3}C\Big(\frac{x(1+x)}{(1-x)^2}\Big)^3, \frac{x(1+x)}{(1-x)^2}C\Big(\frac{x(1+x)}{(1-x)^2}\Big)^2\Big). $$
\end{theorem}
\pf For any G-Motzkin path $\mathbf{G}\in \mathcal{G}$ with at least one $\mathbf{v}$-step at level $i$, given such a $\mathbf{v}$-step, marked as $\mathbf{v}^{*}$, $\mathbf{G}$ can be partitioned uniquely into
$$\mathbf{G}=\mathbf{G}_0\mathbf{u}\mathbf{G}_1\dots \mathbf{u}\mathbf{G}_{i}\mathbf{u}\mathbf{G}_{i+1}\mathbf{v}^{*}\mathbf{\bar{G}}_{0}\mathbf{z}_1\mathbf{\bar{G}}_{1}\dots \mathbf{z}_i\mathbf{\bar{G}}_{i}, $$
where $\mathbf{G}_0, \dots, \mathbf{G}_{i+1}, \mathbf{\bar{G}}_0, \dots, \mathbf{\bar{G}}_{i} \in \mathcal{G}$ and $\mathbf{z}_1, \dots, \mathbf{z}_{i} \in \{\mathbf{d}, \mathbf{v}\}$. Similar to the proof of Theorem \ref{theom 3.1.1}, according to the length of $\mathbf{G}$, all $\mathbf{G}\in \mathcal{G}$ produce the generating function
$x^{i+1}(1+x)^{i}G(x)^{2i+3}$. Hence, the total number $\beta_{n,i}$ of $\mathbf{v}^{*}$-steps in all G-Motzkin paths $\mathbf{G}\in \mathcal{G}_{n+1}$ is the coefficient of $x^{n+1}$ in $x^{i+1}(1+x)^{i}G(x)^{2i+3}$, namely,
$$\beta_{n,i}=[x^{n+1}]x^{i+1}(1+x)^{i}G(x)^{2i+3}=[x^n]G(x)^3\Big(x(1+x)G(x)^2\Big)^{i}. $$
By (\ref{eqn 2.3}), $\beta_{n,i}$ is the $(n,i)$-entry of the Riordan array
$$\Big(G(x)^3,  x(1+x)G(x)^2\Big)=\Big(\frac{1}{(1-x)^3}C\Big(\frac{x(1+x)}{(1-x)^2}\Big)^3, \frac{x(1+x)}{(1-x)^2}C\Big(\frac{x(1+x)}{(1-x)^2}\Big)^2\Big). $$
By the relation $B_{j,i}=[x^j]C(x)^3(xC(x)^2)^{i}$ and by Lemma \ref{lemma 3.2.1}, we have
\begin{eqnarray*}
\gamma_{n,i}\hskip-.22cm &=&\hskip-.22cm \beta_{n, i}=[x^n]\frac{1}{(1-x)^3}C\Big(\frac{x(1+x)}{(1-x)^2}\Big)^3\Big(\frac{x(1+x)}{(1-x)^2}C\Big(\frac{x(1+x)}{(1-x)^2}\Big)^2\Big)^{i} \\
\hskip-.22cm &=&\hskip-.22cm [x^n]\frac{1}{(1-x)^3}\sum_{j=i}^{\infty}B_{j,i}\frac{x^j(1+x)^j}{(1-x)^{2j}}=[x^n]\sum_{j=i}^{\infty}B_{j,i}\frac{x^j(1+x)^{j}}{(1-x)^{2j+3}}\\
\hskip-.22cm &=&\hskip-.22cm \sum_{j=i}^{n}B_{j,i}\sum_{k=0}^{n-j}\binom{j}{k}\binom{n+j-k+2}{n-j-k}.
\end{eqnarray*}
This completes the proof of Theorem \ref{theom 3.2.2}.  \qed\vskip0.2cm

By Theorem \ref{theom 3.2.2}, similar to the proofs of Theorem \ref{theom 3.1.2} and Theorem \ref{theom 3.1.3}, we have
\begin{theorem}\label{theom 3.2.3}
For any integers $n, m\geq 0$, there holds
\begin{eqnarray}
\sum_{i=0}^{n}(-1)^{n-i}\binom{n}{i}\beta_{n+m+i,m+i}  \hskip-.22cm &=&\hskip-.22cm  5^{n}, \label{eqn 3.2.1} \\
\sum_{i=0}^{n}(-1)^{i}\binom{i+2}{2}\beta_{n,i}        \hskip-.22cm &=&\hskip-.22cm  \binom{n+2}{2}. \label{eqn 3.2.2}
\end{eqnarray}
\end{theorem}
One can be asked for combinatorial proofs for these two identities (\ref{eqn 3.2.1}) and (\ref{eqn 3.2.2}).

Note that any $\mathbf{u}$-step at level $i+1$ in a G-Motzkin path $\mathbf{P}$ of length $n+1$ has a matching step, it is a $\mathbf{v}$-step or $\mathbf{d}$-step at level $i$, together with Lemma \ref{lemma 3.2.1}, we have
\begin{corollary}\label{coro 3.2.4}
For any integers $n\geq i\geq 0$, there holds
\begin{eqnarray*}
\alpha_{n,i}=\beta_{n,i}+\beta_{n-1,i}.
\end{eqnarray*}
\end{corollary}

\subsection{ The statistics ``number of $\mathbf{h}$-steps" and ``number of points" at level $i$ }

Let $\mu_{n, i}$ denote the number of $\mathbf{h}$-steps at level $i$ in all G-Motzkin paths of length $n+1$ and let $\lambda_{n, i}$ denote the number of points at level $i$ in all G-Motzkin paths of length $n$, the first values of $\mu_{n, i}$ and $\lambda_{n, i}$ are illustrated in Table 3.5.

\begin{center}
\begin{eqnarray*}
\begin{array}{c|ccccccc}\hline
n/i & 0     & 1     & 2      & 3      & 4    & 5    & 6       \\\hline
  0 & 1     &       &        &        &      &      &        \\
  1 & 4     & 1     &        &        &      &      &        \\
  2 & 18    & 9     & 1      &        &      &      &         \\
  3 & 86    & 60    & 14     & 1      &      &      &         \\
  4 & 431   & 368   & 127    & 19     & 1    &      &         \\
  5 & 2238  & 2190  & 970    & 219    & 24   & 1    &         \\
  6 & 11941 & 12894 & 6803   & 2017   & 336  & 29   & 1      \\\hline
\end{array}
\end{eqnarray*}
Table 3.5. The first values of $\mu_{n, i}$ and $\lambda_{n, i}$.
\end{center}

\begin{lemma}\label{lemma 3.3.1}
For any integers $n\geq i\geq 0$, there holds
\begin{eqnarray*}
\lambda_{n, i}=\mu_{n,i}.
\end{eqnarray*}
\end{lemma}
\pf Given a point at level $i$ in a G-Motzkin path $\mathbf{P}$ of length $n$, insert an $\mathbf{h}$-step into the point, one get a G-Motzkin path $\mathbf{P}'$ of length $n+1$ with an $\mathbf{h}$-step counted at level $i$. Conversely, given an $\mathbf{h}$-step at level $i$ in a G-Motzkin path $\mathbf{P}'$ of length $n+1$, remove the $\mathbf{h}$-step, one get a G-Motzkin path $\mathbf{P}$ of length $n$ with a point counted at level $i$. This one-to-one mapping implies that $\lambda_{n, i}=\mu_{n,i}$.  \qed\vskip0.2cm

\begin{theorem}\label{theom 3.3.2}
For any integers $n\geq i\geq 0$, there holds
\begin{eqnarray*}
\mu_{n,i}=\lambda_{n, i}=\sum_{j=i}^{n}C_{j,i}\sum_{k=0}^{n-j}\binom{j}{k}\binom{n+j-k+1}{n-j-k}.
\end{eqnarray*}
Moreover, $\mu_{n,i}$ is the $(n,i)$-entry of the
Riordan array
$$\Big(\frac{1}{(1-x)^2}C\Big(\frac{x(1+x)}{(1-x)^2}\Big)^2, \frac{x(1+x)}{(1-x)^2}C\Big(\frac{x(1+x)}{(1-x)^2}\Big)^2\Big).$$
\end{theorem}
\pf For any G-Motzkin path $\mathbf{G}\in \mathcal{G}$ with at least one $\mathbf{h}$-step at level $i$, given such an $\mathbf{h}$-step, marked as $\mathbf{h}^{*}$, $\mathbf{G}$ can be partitioned uniquely into
$$\mathbf{G}=\mathbf{G}_0\mathbf{u}\mathbf{G}_1\dots \mathbf{u}\mathbf{G}_{i}\mathbf{h}^{*}\mathbf{\bar{G}}_{0}\mathbf{z}_1\mathbf{\bar{G}}_{1}\dots \mathbf{z}_i\mathbf{\bar{G}}_{i}, $$
where $\mathbf{G}_0, \dots, \mathbf{G}_{i}, \mathbf{\bar{G}}_0, \dots, \mathbf{\bar{G}}_{i} \in \mathcal{G}$ and $\mathbf{z}_1, \dots, \mathbf{z}_{i} \in \{\mathbf{d}, \mathbf{v}\}$. Similar to the proof of Theorem \ref{theom 3.1.1}, according to the length of $\mathbf{G}$, all $\mathbf{G}\in \mathcal{G}$ produce the generating function
$x^{i+1}(1+x)^{i}G(x)^{2i+2}$. Hence, the total number $\mu_{n,i}$ of $\mathbf{h}^{*}$-steps in all G-Motzkin paths $\mathbf{G}\in \mathcal{G}_{n+1}$ is the coefficient of $x^{n+1}$ in $x^{i+1}(1+x)^{i}G(x)^{2i+2}$, namely,
$$\mu_{n,i}=[x^{n+1}]x^{i+1}(1+x)^{i}G(x)^{2i+2}=[x^n]G(x)^2\Big(x(1+x)G(x)^2\Big)^{i}. $$
By (\ref{eqn 2.3}), $\mu_{n,i}$ is the $(n,i)$-entry of the Riordan array
$$\Big(G(x)^2,  x(1+x)G(x)^2\Big)=\Big(\frac{1}{(1-x)^2}C\Big(\frac{x(1+x)}{(1-x)^2}\Big)^2, \frac{x(1+x)}{(1-x)^2}C\Big(\frac{x(1+x)}{(1-x)^2}\Big)^2\Big). $$
By the relation $C_{j,i}=[x^j]C(x)^2(xC(x)^2)^{i}$ and by Lemma \ref{lemma 3.3.1}, we have
\begin{eqnarray*}
\lambda_{n, i}\hskip-.22cm &=&\hskip-.22cm \mu_{n,i}=[x^n]\frac{1}{(1-x)^2}C\Big(\frac{x(1+x)}{(1-x)^2}\Big)^2\Big(\frac{x(1+x)}{(1-x)^2}C\Big(\frac{x(1+x)}{(1-x)^2}\Big)^2\Big)^{i} \\
\hskip-.22cm &=&\hskip-.22cm [x^n]\frac{1}{(1-x)^2}\sum_{j=i}^{\infty}C_{j,i}\frac{x^j(1+x)^j}{(1-x)^{2j}}=[x^n]\sum_{j=i}^{\infty}C_{j,i}\frac{x^j(1+x)^{j}}{(1-x)^{2j+2}}\\
\hskip-.22cm &=&\hskip-.22cm \sum_{j=i}^{n}C_{j,i}\sum_{k=0}^{n-j}\binom{j}{k}\binom{n+j-k+1}{n-j-k}.
\end{eqnarray*}
This completes the proof of Theorem \ref{theom 3.3.2}.    \qed\vskip0.2cm

By Theorem \ref{theom 3.3.2}, similar to the proofs of Theorem \ref{theom 3.1.2} and Theorem \ref{theom 3.1.3}, we have
\begin{theorem}\label{theom 3.3.3}
For any integers $n, m\geq 0$, there holds
\begin{eqnarray}
\sum_{i=0}^{n}(-1)^{n-i}\binom{n}{i}\mu_{n+m+i,m+i} \hskip-.22cm &=&\hskip-.22cm  5^{n}, \label{eqn 3.3.1}\\
\sum_{i=0}^{n}(-1)^{i}(i+1)\mu_{n,i}                \hskip-.22cm &=&\hskip-.22cm  n+1. \label{eqn 3.3.2}
\end{eqnarray}
\end{theorem}

One can be asked for combinatorial proofs of these two identities (\ref{eqn 3.3.1}) and (\ref{eqn 3.3.2}).

\begin{remark}
Note that if replacing the marked $\mathbf{h}$-step $\mathbf{h}^{*}$ by marked $\mathbf{uv}$-peak, one can derive that the number of $\mathbf{uv}$-peaks at level $i+1$ in all G-Motzkin paths of length $n+1$ is $\mu_{n, i}$. Similarly, if replacing the marked $\mathbf{h}$-step $\mathbf{h}^{*}$ by marked $\mathbf{ud}$-peak, one can derive that the number of $\mathbf{ud}$-peaks at level $i+1$ in all G-Motzkin paths of length $n+2$ is also $\mu_{n, i}$.
\end{remark}

\subsection{ The statistics ``number of return steps" }

Let $r_{n, i}$ denote the number of G-Motzkin paths of length $n$ with $i$ return steps, the first values of $r_{n, i}$ are illustrated in Table 3.6.

\begin{center}
\begin{eqnarray*}
\begin{array}{c|ccccccc}\hline
n/i & 0   & 1     & 2      & 3      & 4    & 5    & 6    \\\hline
  0 & 1   &       &        &        &      &      &     \\
  1 & 1   & 1     &        &        &      &      &     \\
  2 & 1   & 5     & 1      &        &      &      &      \\
  3 & 1   & 18    & 9      & 1      &      &      &      \\
  4 & 1   & 67    & 51     & 13     & 1    &      &      \\
  5 & 1   & 278   & 253    & 100    & 17   & 1    &      \\
  6 & 1   & 1272  & 1236   & 623    & 165  & 21   & 1   \\\hline
\end{array}
\end{eqnarray*}
Table 3.6. The first values of $r_{n, i}$.
\end{center}

\begin{theorem}\label{theom 3.4.1}
For any integers $n\geq i\geq 0$, there holds
\begin{eqnarray*}
r_{n, i}=\sum_{j=i}^{n}\frac{i}{2j-i}\binom{2j-i}{j}\sum_{k=0}^{n-j}\binom{j}{k}\binom{n+j-k}{n-j-k}.
\end{eqnarray*}
Moreover, $r_{n,i}$ is the $(n,i)$-entry of the
Riordan array
$$\Big(\frac{1}{1-x}, \frac{x(1+x)}{(1-x)^2}C\Big(\frac{x(1+x)}{(1-x)^2}\Big)\Big).$$
\end{theorem}
\pf For any G-Motzkin path $\mathbf{G}\in \mathcal{G}$ with $i$ return steps, $\mathbf{G}$ can be partitioned uniquely into
$$\mathbf{G}=\mathbf{h}^{k_0}\mathbf{u}\mathbf{G}_1\mathbf{z}_1\mathbf{h}^{k_1}\mathbf{u}\mathbf{G}_2\mathbf{z}_2\mathbf{h}^{k_2}\dots \mathbf{u}\mathbf{G}_{i}\mathbf{z}_i\mathbf{h}^{k_{i}}, $$
where $\mathbf{G}_1, \dots, \mathbf{G}_{i}\in \mathcal{G}$, $\mathbf{z}_1, \dots, \mathbf{z}_{i} \in \{\mathbf{d}, \mathbf{v}\}$ and $k_0, k_1, \dots, k_i\geq 0$. Similar to the proof of Theorem \ref{theom 3.1.1}, according to the length of $\mathbf{G}$, all $\mathbf{G}\in \mathcal{G}$ produce the generating function
$\frac{1}{(1-x)^{i+1}}x^{i}(1+x)^{i}G(x)^{i}$. Hence, the total number $r_{n,i}$ of return steps in all G-Motzkin paths $\mathbf{G}\in \mathcal{G}_{n}$ is the coefficient of $x^{n}$ in $\frac{1}{(1-x)^{i+1}}x^{i}(1+x)^{i}G(x)^{i}$, namely,
$$r_{n,i}=[x^{n}]\frac{1}{(1-x)^{i+1}}x^{i}(1+x)^{i}G(x)^{i}=[x^n]\frac{1}{1-x}\Big(\frac{x(1+x)}{1-x}G(x)\Big)^{i}. $$
By (\ref{eqn 2.3}), $r_{n,i}$ is the $(n,i)$-entry of the Riordan array
$$\Big(\frac{1}{1-x}, \frac{x(1+x)}{1-x}G(x)\Big)=\Big(\frac{1}{1-x}, \frac{x(1+x)}{(1-x)^2}C\Big(\frac{x(1+x)}{(1-x)^2}\Big)\Big). $$
By the relation $\frac{i}{2j-i}\binom{2j-i}{j}=[x^j](xC(x))^{i}$, we have
\begin{eqnarray*}
r_{n,i}\hskip-.22cm &=&\hskip-.22cm [x^n]\frac{1}{1-x}\Big(\frac{x(1+x)}{(1-x)^2}C\Big(\frac{x(1+x)}{(1-x)^2}\Big)\Big)^{i} \\
\hskip-.22cm &=&\hskip-.22cm [x^n]\frac{1}{1-x}\sum_{j=i}^{\infty}\frac{i}{2j-i}\binom{2j-i}{j}\frac{x^j(1+x)^j}{(1-x)^{2j}} \\
\hskip-.22cm &=&\hskip-.22cm \sum_{j=i}^{n}\frac{i}{2j-i}\binom{2j-i}{j}\sum_{k=0}^{n-j}\binom{j}{k}\binom{n+j-k}{n-j-k}.
\end{eqnarray*}
This completes the proof of Theorem \ref{theom 3.4.1}.    \qed\vskip0.2cm

By Theorem \ref{theom 3.4.1}, similar to the proof of Theorem \ref{theom 3.1.2}, we have
\begin{theorem}\label{theom 3.4.2}
For any integers $n, m\geq 0$, there holds
\begin{eqnarray}\label{eqn 3.4.1}
\sum_{i=0}^{n}(-1)^{n-i}\binom{n}{i}r_{n+m+i,m+i}=4^{n}.
\end{eqnarray}
\end{theorem}
One can be asked for a combinatorial proof of this identity (\ref{eqn 3.4.1}).

\section{The statistics ``number of $\mathbf{z}_1\mathbf{z}_2$-steps" in G-Motzkin paths}

In this section, we discuss the statistics ``number of $\mathbf{z}_1\mathbf{z}_2$-steps" in G-Motzkin paths for $\mathbf{z}_1, \mathbf{z}_2\in \{\mathbf{u}, \mathbf{h}, \mathbf{v}, \mathbf{d}\}$. Despite that there are $16$ cases to be considered, but in fact it only needs to study $10$ cases in the set $\{\mathbf{ud}, \mathbf{uh}, \mathbf{uu}, \mathbf{hh}, \mathbf{hd}, \mathbf{vu}, \mathbf{vv}, \mathbf{du}, \mathbf{dd}, \mathbf{dv}\}.$
Let $L_{n, i}^{\mathbf{z}_1\mathbf{z}_2}$ denote the number of G-Motzkin paths of length $n$ with $i$ $\mathbf{z}_1\mathbf{z}_2$ steps, it is not difficult to verify that

$$L_{n, i}^{\mathbf{uv}}=H_{n,i}, \ L_{n, i}^{\mathbf{hv}}=L_{n, i}^{\mathbf{vh}}=D_{n,i}, \ L_{n, i}^{\mathbf{hu}}=L_{n, i}^{\mathbf{uh}}, \ L_{n, i}^{\mathbf{hd}}=L_{n, i}^{\mathbf{dh}},   \ L_{n, i}^{\mathbf{dv}}=L_{n, i}^{\mathbf{vd}}.   $$

The results, $L_{n, i}^{\mathbf{uv}}=H_{n,i}$, $L_{n, i}^{\mathbf{hv}}=L_{n, i}^{\mathbf{vh}}=D_{n,i}$, can be easily obtained respectively by replacing all the $\mathbf{uv}$-steps ($\mathbf{hv}$-steps, $\mathbf{vh}$-steps) by new $\mathbf{h}$-steps ($\mathbf{d}$-steps) and replacing the old $\mathbf{h}$-steps ($\mathbf{d}$-steps) by $\mathbf{uv}$-steps ($\mathbf{hv}$-steps, $\mathbf{vh}$-steps) in each G-Motzkin paths of length $n$ and vice versa.
The results, $L_{n, i}^{\mathbf{hu}}=L_{n, i}^{\mathbf{uh}}$, $L_{n, i}^{\mathbf{hd}}=L_{n, i}^{\mathbf{dh}}$, $L_{n, i}^{\mathbf{dv}}=L_{n, i}^{\mathbf{vd}}$, can be easily obtained respectively by replacing all the $\mathbf{hu}$-steps ($\mathbf{hd}$-steps, $\mathbf{dv}$-steps) by new $\mathbf{uh}$-steps ($\mathbf{dh}$-steps, $\mathbf{vd}$-steps) and replacing the old $\mathbf{uh}$-steps ($\mathbf{dh}$-steps, $\mathbf{vd}$-steps) (if exist) by $\mathbf{hu}$-steps ($\mathbf{hd}$-steps, $\mathbf{dv}$-steps) in each G-Motzkin paths of length $n$ and vice versa.

\subsection{ The statistics ``number of $\mathbf{ud}$-peaks" }

Let $L^{\mathbf{ud}}(x, y)=\sum_{n=0}^{\infty}\sum_{i=0}^{\infty}L_{n, i}^{\mathbf{ud}}x^ny^i$ be the generating function of $L_{n, i}^{\mathbf{ud}}$, the number of G-Motzkin paths of length $n$ with $i$ $\mathbf{ud}$-peaks. According to the
first return decomposition, a G-Motzkin path $\mathbf{P}$ can be decomposed as one of the following five forms:
$$\mathbf{P}=\varepsilon, \ \mathbf{P}=\mathbf{h}\mathbf{Q}_1,\ \mathbf{P}=\mathbf{u}\mathbf{Q}_1\mathbf{v}\mathbf{Q}_2, \ \mathbf{P}=\mathbf{u}\mathbf{d}\mathbf{Q}_2, \ \mbox{or}\  \mathbf{P}=\mathbf{u}\mathbf{P}_1\mathbf{d}\mathbf{Q}_2. $$
where $\mathbf{Q}_1$ and $\mathbf{Q}_2 $ are (possibly empty) G-Motzkin paths and $\mathbf{P}_1$ are nonempty. Then we get the relation
\begin{eqnarray*}
L^{\mathbf{ud}}(x, y)\hskip-.22cm &=&\hskip-.22cm 1+xL^{\mathbf{ud}}(x, y)+xL^{\mathbf{ud}}(x, y)^2+x^2yL^{\mathbf{ud}}(x, y)+x^2(L^{\mathbf{ud}}(x, y)-1)L^{\mathbf{ud}}(x, y) \\
\hskip-.22cm &=&\hskip-.22cm 1+x(1-x+xy)L^{\mathbf{ud}}(x, y)+x(1+x)L^{\mathbf{ud}}(x, y)^2.
\end{eqnarray*}
Solve this, we have
\begin{eqnarray}\label{eqn 4.1}
L^{\mathbf{ud}}(x, y)\hskip-.22cm &=&\hskip-.22cm   \frac{1-x+x^2-x^2y-\sqrt{(1-x+x^2-x^2y)^2-4x(1+x)}}{2x(1+x)}  \nonumber\\
                     \hskip-.22cm &=&\hskip-.22cm   \frac{1}{1-x+x^2-x^2y}C\Big(\frac{x(1+x)}{(1-x+x^2-x^2y)^2}\Big).
\end{eqnarray}
By (\ref{eqn 1.1}) and (\ref{eqn 4.1}), taking the coefficient of $x^ny^i$ in $L^{\mathbf{ud}}(x, y)$, we derive that
\begin{theorem}
For any integers $n\geq i\geq 0$, there holds
\begin{eqnarray*}
L_{n, i}^{\mathbf{ud}} \hskip-.22cm &=&\hskip-.22cm \sum_{k=0}^{n-2i}\sum_{j=0}^{[\frac{n-k-2i}{3}]}(-1)^j\binom{2k+i}{i}\binom{2k+i+j}{j}\binom{3k+i+1}{n-k-2i-3j}C_k.
\end{eqnarray*}
\end{theorem}

The first values of $L_{n, i}^{\mathbf{ud}}$ are illustrated in Table 4.1.
\begin{center}
\begin{eqnarray*}
\begin{array}{c|ccccccc}\hline
n/i & 0        & 1        & 2       & 3     & 4    & 5        \\\hline
  0 & 1        &          &         &       &      &         \\
  1 & 2        &          &         &       &      &         \\
  2 & 6        & 1        &         &       &      &          \\
  3 & 24       & 5        &         &       &      &          \\
  4 & 106      & 26       & 1       &       &      &          \\
  5 & 498      & 143      & 9       &       &      &         \\
  6 & 2444     & 805      & 69      &  1    &      &         \\
  7 & 12382    & 4604     & 498     &  14   &      &         \\
  8 & 64270    & 26637    & 3471    &  146  &  1   &         \\\hline
\end{array}
\end{eqnarray*}
Table 4.1. The first values of $L_{n, i}^{\mathbf{ud}}$.
\end{center}

\begin{theorem}\label{theom 4.1.2ud}
For any integers $n, m\geq 0$, there holds
\begin{eqnarray*}
\sum_{i=0}^{2n}(-1)^{i}\binom{2n}{i}L_{n+2m+2i, m+i}^{\mathbf{ud}}=C_n.
\end{eqnarray*}
\end{theorem}
\pf Note that
\begin{eqnarray*}
L_{n+2m+2i, m+i}^{\mathbf{ud}} \hskip-.22cm &=&\hskip-.22cm \sum_{k=0}^{n}\sum_{j=0}^{[\frac{n-k}{3}]}(-1)^j\binom{2k+m+i}{2k}\binom{2k+m+i+j}{j}\binom{3k+m+i+1}{n-k-3j}C_k,
\end{eqnarray*}
each inner term in $L_{n+2m+2i, m+i}^{\mathbf{ud}}$, denoted by $h_{n, m, k, j}(i)$, is a polynomial on $i$ with degree
$$\partial h_{n, m, k, j}(i)=2k+j+(n-k-3j)=n+k-2j\leq 2n $$
such that $\partial h_{n, m, k, j}(i) = 2n$ if and only if $k=n$ and $j=0$. Clearly, the leading term in
\begin{eqnarray*}
h_{n, m, n, 0}(i) \hskip-.22cm &=&\hskip-.22cm \binom{2n+m+i}{2n}C_n
\end{eqnarray*}
is $\frac{C_n}{(2n)!}i^{2n}$. So $L_{n+2m+2i, m+i}^{\mathbf{ud}}$ is also a polynomial on $i$ with degree $2n$ such that the leading term is $\frac{C_n}{(2n)!}i^{2n}$. By the
Euler difference identity
\begin{eqnarray*}
\sum_{i=0}^{n}(-1)^{i}\binom{n}{i}i^{r}=\left\{
\begin{array}{rl}
0,          &  \mbox{if}\ 0\leq r<n, \\[5pt]
(-1)^nn!,   &  \mbox{if}\ r=n,
\end{array}\right.
\end{eqnarray*}
we have
\begin{eqnarray*}
\sum_{i=0}^{2n}(-1)^{i}\binom{2n}{i}L_{n+2m+2i, m+i}^{\mathbf{ud}}=\sum_{i=0}^{2n}(-1)^{i}\binom{2n}{i}i^{2n}\frac{C_n}{(2n)!}=C_n.
\end{eqnarray*}
This completes the proof.  \qed\vskip0.2cm

By (\ref{eqn 1.1}) and (\ref{eqn 4.1}), taking the coefficient of $x^n$ in $L^{\mathbf{ud}}(x, -\frac{3}{x})$, we derive that
\begin{theorem}\label{theom 4.1.3ud}
For any integer $n\geq 0$, there holds
\begin{eqnarray*}
\sum_{i=0}^{n}(-3)^{i}L_{n+i, i}^{\mathbf{ud}}=\sum_{i=0}^{n}(-1)^{n-k}\binom{n+2k+1}{n-k}C_k.
\end{eqnarray*}
\end{theorem}

\subsection{ The statistics ``number of $\mathbf{uh}$-steps" }

Let $L^{\mathbf{uh}}(x, y)=\sum_{n=0}^{\infty}\sum_{i=0}^{\infty}L_{n, i}^{\mathbf{uh}}x^ny^i$ be the generating function of $L_{n, i}^{\mathbf{uh}}$, the number of G-Motzkin paths of length $n$ with $i$ $\mathbf{uh}$-steps. According to the
first return decomposition, a G-Motzkin path $\mathbf{P}$ can be decomposed as one of the following four forms:
$$\mathbf{P}=\varepsilon, \ \mathbf{P}=\mathbf{h}\mathbf{Q}_2,\ \mathbf{P}=\mathbf{u}\mathbf{P}_1\mathbf{v}\mathbf{Q}_2, \ \mbox{or}\  \mathbf{P}=\mathbf{u}\mathbf{P}_1\mathbf{d}\mathbf{Q}_2. $$
where $\mathbf{P}_1$ and $\mathbf{Q}_2$ are (possibly empty) G-Motzkin paths. If $\mathbf{P}_1$ begins with an $\mathbf{h}$-step, $\mathbf{u}\mathbf{P}_1$ contributes at least one $\mathbf{uh}$-step. Then we get the relation
\begin{eqnarray*}
L^{\mathbf{uh}}(x, y)\hskip-.22cm &=&\hskip-.22cm 1+xL^{\mathbf{uh}}(x, y)+x(1-x)L^{\mathbf{uh}}(x, y)^2   \\
\hskip-.22cm & &\hskip-.22cm \  \  +\ x^2yL^{\mathbf{uh}}(x, y)^2+x^2(1-x)L^{\mathbf{uh}}(x, y)^2+x^3yL^{\mathbf{uh}}(x, y)^2 \\
\hskip-.22cm &=&\hskip-.22cm 1+xL^{\mathbf{uh}}(x, y)+x(1+x)(1-x+xy)L^{\mathbf{uh}}(x, y)^2.
\end{eqnarray*}
Solve this, we have
\begin{eqnarray}\label{eqn 4.2}
L^{\mathbf{uh}}(x, y)\hskip-.22cm &=&\hskip-.22cm   \frac{1-x-\sqrt{(1-x)^2-4x(1+x)(1-x+xy)}}{2x(1+x)(1-x+xy)}  \nonumber\\
                     \hskip-.22cm &=&\hskip-.22cm   \frac{1}{1-x}C\Big(\frac{x(1+x)(1-x+xy)}{(1-x)^2}\Big).
\end{eqnarray}
By (\ref{eqn 1.1}) and (\ref{eqn 4.2}), taking the coefficient of $x^ny^i$ in $L^{\mathbf{uh}}(x, y)$, we derive that
\begin{theorem}
For any integers $n\geq i\geq 0$, there holds
\begin{eqnarray*}
L_{n, i}^{\mathbf{uh}} \hskip-.22cm &=&\hskip-.22cm \sum_{k=i}^{n-i}\sum_{j=0}^{n-k-i}\binom{k}{i}\binom{k}{j}\binom{n-j}{n-k-i-j}C_k.
\end{eqnarray*}
\end{theorem}

The first values of $L_{n, i}^{\mathbf{uh}}$ are illustrated in Table 4.2. Note that $L_{2n, n}^{\mathbf{uh}}=C_n$, since any Dyck paths of length $2n$ can be obtained
by replacing each $\mathbf{uh}$-step in any G-Motzkin paths of length $2n$ with $n$ $\mathbf{uh}$-steps (no $\mathbf{d}$-steps implied) by a $\mathbf{u}$-step and each $\mathbf{v}$-step by a $\mathbf{d}$-step, and vice versa.

\begin{center}
\begin{eqnarray*}
\begin{array}{c|ccccccc}\hline
n/i & 0     & 1     & 2      & 3     & 4    & 5        \\\hline
  0 & 1     &       &        &       &      &         \\
  1 & 2     &       &        &       &      &         \\
  2 & 6     & 1     &        &       &      &          \\
  3 & 21    & 8     &        &       &      &          \\
  4 & 83    & 48    &  2     &       &      &          \\
  5 & 353   & 268   &  29    &       &      &         \\
  6 & 1577  & 1466  &  271   &  5    &      &         \\
  7 & 7294  & 7984  &  2114  &  106  &      &         \\
  8 & 34622 & 43509 &  15028 &  1352 &  14  &         \\\hline
\end{array}
\end{eqnarray*}
Table 4.2. The first values of $L_{n, i}^{\mathbf{uh}}$.
\end{center}

\subsection{ The statistics ``number of $\mathbf{uu}$-steps" }

Let $L^{\mathbf{uu}}(x, y)=\sum_{n=0}^{\infty}\sum_{i=0}^{\infty}L_{n, i}^{\mathbf{uu}}x^ny^i$ be the generating function of $L_{n, i}^{\mathbf{uu}}$, the number of G-Motzkin paths of length $n$ with $i$ $\mathbf{uu}$-steps. According to the
first return decomposition, a G-Motzkin path $\mathbf{P}$ can be decomposed as one of the following four forms:
$$\mathbf{P}=\varepsilon, \ \mathbf{P}=\mathbf{h}\mathbf{P}_0,\ \mathbf{P}=\mathbf{u}^k\mathbf{h}\mathbf{P}_0\mathbf{z}_1\mathbf{P}_{1}\mathbf{z}_2\mathbf{P}_{2}\dots \mathbf{z}_{k}\mathbf{P}_{k}, \ \mbox{or}\  \mathbf{P}=\mathbf{u}^k\mathbf{z}_1\mathbf{P}_{1}\mathbf{z}_2\mathbf{P}_{2}\dots \mathbf{z}_{k}\mathbf{P}_{k}. $$
where $\mathbf{P}_0, \mathbf{P}_1, \dots, \mathbf{P}_k$ are (possibly empty) G-Motzkin paths for certain $k\geq 1$ and
$\mathbf{z}_1,\dots, \mathbf{z}_k\in \{\mathbf{d}, \mathbf{v}\}$. Note that the last two cases contribute at least $k-1$ $\mathbf{uu}$-steps. Then we get the relation
\begin{eqnarray*}\small
L^{\mathbf{uu}}(x, y)\hskip-.22cm &=&\hskip-.22cm 1+xL^{\mathbf{uu}}(x, y)+\sum_{k=1}^{\infty}x^{k+1}(1+x)^ky^{k-1}
L^{\mathbf{uu}}(x, y)^{k+1}+\sum_{k=1}^{\infty}x^{k}(1+x)^ky^{k-1}L^{\mathbf{uu}}(x, y)^{k} \\
\hskip-.22cm &=&\hskip-.22cm 1+xL^{\mathbf{uu}}(x, y)+\frac{x^2(1+x)
L^{\mathbf{uu}}(x, y)^2}{1-x(1+x)yL^{\mathbf{uu}}(x, y)}+
\frac{x(1+x)L^{\mathbf{uu}}(x, y)}{1-x(1+x)yL^{\mathbf{uu}}(x, y)}  \\
\hskip-.22cm &=&\hskip-.22cm \big(1+xL^{\mathbf{uu}}(x, y)\big
)\Big(1+\frac{x(1+x)L^{\mathbf{uu}}(x, y)}{1-x(1+x)yL^{\mathbf{uu}}(x, y)}\Big).
\end{eqnarray*}
Solve this, we have
\begin{eqnarray}\label{eqn 4.3}
L^{\mathbf{uu}}(x, y)\hskip-.22cm &=&\hskip-.22cm   \frac{1-2x-x^2+x(1+x)y-\sqrt{(1-2x-x^2+x(1+x)y)^2-4x(1+x)(x+y-xy)}}{2x(1+x)(x+y-xy)}  \nonumber\\
                     \hskip-.22cm &=&\hskip-.22cm   \frac{1}{1-2x-x^2+x(1+x)y}C\Big(\frac{x(1+x)(x+y-xy)}{(1-2x-x^2+x(1+x)y)^2}\Big).
\end{eqnarray}

By (\ref{eqn 1.1}) and (\ref{eqn 4.3}), taking the coefficient of $x^ny^{i}$ in $L^{\mathbf{uu}}(x, y)$ in two different ways, we have
\begin{theorem}
For any integers $n\geq i\geq 0$, there holds
\begin{eqnarray*}
L_{n, i}^{\mathbf{uu}} \hskip-.22cm &=&\hskip-.22cm \sum_{k=0}^{n}\sum_{j=0 }^{k}\sum_{r=0}^{n-k-j}\sum_{\ell=0}^{k+r}(-1)^{i+j-k}
\binom{k}{j}\binom{2k+r}{r}\binom{j+r}{i+j-k}\binom{k+r}{\ell}\binom{n+k-j-\ell}{n-k-j-r-\ell}C_k, \\
\hskip-.22cm &=&\hskip-.22cm \sum_{k=0}^{n}\sum_{j=0 }^{k}\sum_{r=0}^{n-k-j}\sum_{\ell=0}^{k+r}(-1)^{i+j-k}
\binom{k}{j}\binom{2k+r}{r}\binom{r}{i+j-k}\binom{k+r}{\ell}\binom{n-\ell}{n-k-j-r-\ell}C_k.
\end{eqnarray*}
\end{theorem}

The first values of $L_{n, i}^{\mathbf{uu}}$ are illustrated in Table 4.3.
\begin{center}
\begin{eqnarray*}
\begin{array}{c|ccccccc}\hline
n/i & 0      & 1      & 2    & 3     & 4    & 5        \\\hline
  0 & 1      &        &      &       &      &         \\
  1 & 2      &        &      &       &      &         \\
  2 & 6      & 1      &      &       &      &          \\
  3 & 19     & 9      & 1    &       &      &          \\
  4 & 64     & 54     & 14   &  1    &      &          \\
  5 & 225    & 282    & 122  &  20   &  1   &         \\
  6 & 828    & 1386   & 847  &  230  &  27  &  1       \\\hline
\end{array}
\end{eqnarray*}
Table 4.3. The first values of $L_{n, i}^{\mathbf{uu}}$.
\end{center}

Specially, when $i=n\geq 1$, it implies that $j=0, r=n-k$ and $\ell =0$, so $L_{0, 0}^{\mathbf{uu}}=1$ and $L_{n, n}^{\mathbf{uu}}=0$ for $n\geq 1$ produce that
\begin{eqnarray*}
\sum_{k=0}^{n} (-1)^{n-k}\binom{n+k}{n-k}C_k=\left\{
\begin{array}{cl}
0,   &  \mbox{if}\ n\geq 1, \\[5pt]
1,   &  \mbox{if}\ n=0,
\end{array}\right.
\end{eqnarray*}
which has been derived by Chen and Pang \cite{ChenPang} and is a special case of the identity \cite{ChenPang, ManSun}
\begin{eqnarray*}
\sum_{k=0}^{n}\binom{n+k}{n-k}C_k(x-1)^{n-k} = \sum_{k=0}^{n}\frac{1}{n}\binom{n}{k-1}\binom{n}{k}x^{k}.
\end{eqnarray*}


\subsection{ The statistics ``number of $\mathbf{hh}$-steps" }

Let $L^{\mathbf{hh}}(x, y)=\sum_{n=0}^{\infty}\sum_{i=0}^{\infty}L_{n, i}^{\mathbf{hh}}x^ny^i$ be the generating function of $L_{n, i}^{\mathbf{hh}}$, the number of G-Motzkin paths of length $n$ with $i$ $\mathbf{hh}$-steps. According to the
first return decomposition, a G-Motzkin path $\mathbf{P}$ can be decomposed as one of the following six forms:
$$\mathbf{P}=\varepsilon, \ \mathbf{P}=\mathbf{h}^{k},\  \ \mathbf{P}=\mathbf{h}^{k}\mathbf{u}\mathbf{P}_1\mathbf{v}\mathbf{P}_2,\ \mathbf{P}=\mathbf{h}^{k}\mathbf{u}\mathbf{P}_1\mathbf{d}\mathbf{P}_2, \ \mathbf{P}=\mathbf{u}\mathbf{P}_1\mathbf{v}\mathbf{P}_2,\ \mbox{or}\  \mathbf{P}=\mathbf{u}\mathbf{P}_1\mathbf{d}\mathbf{P}_2 $$
for certain $k\geq 1$, where $\mathbf{P}_1$ and $\mathbf{P}_2$ are (possibly empty) G-Motzkin paths. Note that the $\mathbf{h}^{k}\mathbf{u}$ part contributes $k-1$ $\mathbf{hh}$-steps. Then we get the relation
\begin{eqnarray*}
L^{\mathbf{hh}}(x, y)\hskip-.22cm &=&\hskip-.22cm 1+\sum_{k=1}^{\infty}x^ky^{k-1}+\sum_{k=1}^{\infty}x^{k+1}y^{k-1}L^{\mathbf{uu}}(x, y)^{2}+\sum_{k=1}^{\infty}x^{k+2}y^{k-1}L^{\mathbf{uu}}(x, y)^{2} \\
\hskip-.22cm & &\hskip-.22cm \ \  +\ xL^{\mathbf{hh}}(x, y)^2+x^2L^{\mathbf{hh}}(x, y)^2 \\
\hskip-.22cm &=&\hskip-.22cm \Big(1+\frac{x}{1-xy}\Big)\big(1+x(1+x)L^{\mathbf{hh}}(x, y)^2\big).
\end{eqnarray*}
Solve this, we have
\begin{eqnarray}\label{eqn 4.4}
L^{\mathbf{hh}}(x, y)\hskip-.22cm &=&\hskip-.22cm   \frac{1-xy-\sqrt{(1-xy)^2-4x(1+x)(1+x-xy)^2}}{2x(1+x)(1+x-xy)}  \nonumber\\
                     \hskip-.22cm &=&\hskip-.22cm   \frac{1+x-xy}{1-xy}C\Big(\frac{x(1+x)(1+x-xy)^2}{(1-xy)^2}\Big)  \nonumber \\
                     \hskip-.22cm &=&\hskip-.22cm   \big(1+\frac{x}{1-xy}\big)C\Big(x(1+x)\big(1+\frac{x}{1-xy}\big)^2\Big).
\end{eqnarray}
By (\ref{eqn 1.1}) and (\ref{eqn 4.4}), taking the coefficient of $x^ny^i$ in $L^{\mathbf{hh}}(x, y)$, we derive that
\begin{theorem}
For any integers $n\geq i\geq 0$, there holds
\begin{eqnarray*}
L_{n, i}^{\mathbf{hh}} \hskip-.22cm &=&\hskip-.22cm \sum_{k=0}^{n-i}\sum_{j=0}^{n-k-i}\binom{2k+1}{j}\binom{i+j-1}{i}\binom{k}{n-k-i-j}C_k.
\end{eqnarray*}
\end{theorem}

The first values of $L_{n, i}^{\mathbf{hh}}$ are illustrated in Table 4.4.
\begin{center}
\begin{eqnarray*}
\begin{array}{c|ccccccc}\hline
n/i & 0       & 1       & 2        & 3     & 4    & 5        \\\hline
  0 & 1       &         &          &       &      &         \\
  1 & 2       &         &          &       &      &         \\
  2 & 6       & 1       &          &       &      &          \\
  3 & 25      & 3       & 1        &       &      &          \\
  4 & 110     & 19      & 3        &  1    &      &          \\
  5 & 520     & 104     & 22       &  3    &  1   &         \\
  6 & 2566    & 594     & 130      &  25   &  3   &  1      \\\hline
\end{array}
\end{eqnarray*}
Table 4.4. The first values of $L_{n, i}^{\mathbf{hh}}$.
\end{center}

\subsection{ The statistics ``number of $\mathbf{hd}$-steps" }

Let $L^{\mathbf{hd}}(x, y)=\sum_{n=0}^{\infty}\sum_{i=0}^{\infty}L_{n, i}^{\mathbf{hd}}x^ny^i$ be the generating function of $L_{n, i}^{\mathbf{hd}}$, the number of G-Motzkin paths of length $n$ with $i$ $\mathbf{hd}$-steps. According to the
first return decomposition, a G-Motzkin path $\mathbf{P}$ can be decomposed as one of the following four forms:
$$\mathbf{P}=\varepsilon, \ \mathbf{P}=\mathbf{h}\mathbf{P}_2,\ \mathbf{P}=\mathbf{u}\mathbf{P}_1\mathbf{v}\mathbf{P}_2, \ \mbox{or}\  \mathbf{P}=\mathbf{u}\mathbf{P}_1\mathbf{d}\mathbf{P}_2. $$
where $\mathbf{P}_1$ and $\mathbf{P}_2$ are (possibly empty) G-Motzkin paths. If $\mathbf{P}_1$ ends with an $\mathbf{h}$-step, $\mathbf{P}_1\mathbf{d}$ contributes at least one $\mathbf{hd}$-step. Then we get the relation
\begin{eqnarray*}
L^{\mathbf{hd}}(x, y)\hskip-.22cm &=&\hskip-.22cm 1+xL^{\mathbf{hd}}(x, y)+xL^{\mathbf{hd}}(x, y)^2+x^2(1-x)L^{\mathbf{hd}}(x, y)^2+x^3yL^{\mathbf{hd}}(x, y)^2 \\
\hskip-.22cm &=&\hskip-.22cm 1+xL^{\mathbf{hd}}(x, y)+x(1+x-x^2+x^2y)L^{\mathbf{hd}}(x, y)^2.
\end{eqnarray*}
Solve this, we have
\begin{eqnarray}\label{eqn 4.5}
L^{\mathbf{hd}}(x, y)\hskip-.22cm &=&\hskip-.22cm   \frac{1-x-\sqrt{(1-x)^2-4x(1+x-x^2+x^2y)}}{2x(1+x-x^2+x^2y)}  \nonumber\\
                     \hskip-.22cm &=&\hskip-.22cm   \frac{1}{1-x}C\Big(\frac{x(1+x-x^2+x^2y)}{(1-x)^2}\Big).
\end{eqnarray}
By (\ref{eqn 1.1}) and (\ref{eqn 4.5}), taking the coefficient of $x^ny^i$ in $L^{\mathbf{hd}}(x, y)$, we derive that
\begin{theorem}
For any integers $n\geq i\geq 0$, there holds
\begin{eqnarray*}
L_{n, i}^{\mathbf{hd}} \hskip-.22cm &=&\hskip-.22cm \sum_{k=i}^{n-2i}\sum_{j=0}^{k-i}\binom{k}{i}\binom{k-i}{j}\binom{n+k-2i-2j}{n-k-2i-j}C_k.
\end{eqnarray*}
\end{theorem}

The first values of $L_{n, i}^{\mathbf{hd}}$ are illustrated in Table 4.5. Note that $L_{3n, n}^{\mathbf{hd}}=C_n$, since any Dyck paths of length $2n$ can be obtained
by replacing each $\mathbf{hd}$-step in any G-Motzkin paths of length $3n$ with $n$ $\mathbf{hd}$-steps ($n$ $\mathbf{u}$-steps and no $\mathbf{v}$-steps implied) by a $\mathbf{d}$-step, and vice versa.
\begin{center}
\begin{eqnarray*}
\begin{array}{c|ccccccc}\hline
n/i & 0      & 1      & 2       & 3     & 4    & 5        \\\hline
  0 & 1      &        &         &       &      &         \\
  1 & 2      &        &         &       &      &         \\
  2 & 7      &        &         &       &      &          \\
  3 & 28     & 1      &         &       &      &          \\
  4 & 126    & 7      &         &       &      &          \\
  5 & 605    & 45     &         &       &      &         \\
  6 & 3040   & 277    & 2       &       &      &         \\
  7 & 15781  & 1692   & 25      &       &      &         \\
  8 & 83971  & 10320  & 234     &       &      &         \\
  9 & 455553 & 63026  & 1924    & 5     &      &         \\\hline
\end{array}
\end{eqnarray*}
Table 4.5. The first values of $L_{n, i}^{\mathbf{hd}}$.
\end{center}

\subsection{ The statistics ``number of $\mathbf{vu}$-valleys" }

Let $L^{\mathbf{vu}}(x, y)=\sum_{n=0}^{\infty}\sum_{i=0}^{\infty}L_{n, i}^{\mathbf{vu}}x^ny^i$ be the generating function of $L_{n, i}^{\mathbf{vu}}$, the number of G-Motzkin paths of length $n$ with $i$ $\mathbf{vu}$-valleys. According to the
first return decomposition, a G-Motzkin path $\mathbf{P}$ can be decomposed as one of the following four forms:
$$\mathbf{P}=\varepsilon, \ \mathbf{P}=\mathbf{h}\mathbf{P}_1,\ \mathbf{P}=\mathbf{u}\mathbf{P}_1\mathbf{v}\dots \mathbf{u}\mathbf{P}_k\mathbf{v}\mathbf{u}\mathbf{P}_{k+1}\mathbf{d}\mathbf{P}_{k+2}, \ \mbox{or}\  \mathbf{P}=\mathbf{u}\mathbf{P}_1\mathbf{v}\dots \mathbf{u}\mathbf{P}_k\mathbf{v}\mathbf{u}\mathbf{P}_{k+1}\mathbf{v}\mathbf{Q}. $$
for certain $k\geq 0$, where $\mathbf{P}_1, \dots, \mathbf{P}_{k+1}$ are (possibly empty) G-Motzkin paths and $\mathbf{Q}$ is empty or begins with an $\mathbf{h}$-step. Note that the last two cases contribute at least $k$ $\mathbf{vu}$-valleys. Then we get the relation
\begin{eqnarray}
L^{\mathbf{vu}}(x, y)\hskip-.22cm &=&\hskip-.22cm 1+xL^{\mathbf{vu}}(x, y)+\sum_{k=0}^{\infty}x^{k+2}y^{k}L^{\mathbf{vu}}(x, y)^{k+2} \nonumber \\
\hskip-.22cm & &\hskip-.22cm \ +\  \sum_{k=0}^{\infty}x^{k+1}y^{k}L^{\mathbf{vu}}(x, y)^{k+1}(1+xL^{\mathbf{vu}}(x, y)) \nonumber \\
\hskip-.22cm &=&\hskip-.22cm 1+xL^{\mathbf{vu}}(x, y)+\frac{x^2L^{\mathbf{vu}}(x, y)^2}{1-xyL^{\mathbf{vu}}(x, y)}+\frac{xL^{\mathbf{vu}}(x, y)(1+xL^{\mathbf{vu}}(x, y))}{1-xyL^{\mathbf{vu}}(x, y)}. \label{eqn 4.6vu}
\end{eqnarray}
Solve this, we have
\begin{eqnarray}
L^{\mathbf{vu}}(x, y)\hskip-.22cm &=&\hskip-.22cm   \frac{1-2x+xy-\sqrt{(1-2x+xy)^2-4x(2x+y-xy)}}{2x(2x+y-xy)}  \nonumber\\
                     \hskip-.22cm &=&\hskip-.22cm   \frac{1}{1-2x+xy}C\Big(\frac{x(2x+y-xy)}{(1-2x+xy)^2}\Big). \label{eqn 4.7vu}
\end{eqnarray}
Let $T=xL^{\mathbf{vu}}(x, y)$, using the Lagrange inversion formula in (\ref{eqn 4.6vu}), taking the coefficient of $x^{n+1}y^i$ in $xL^{\mathbf{vu}}(x, y)$, we derive that
\begin{eqnarray}
L_{n, i}^{\mathbf{vu}} \hskip-.22cm &=&\hskip-.22cm  [x^{n+1}y^i]T=[y^i][x^{n+1}]T=[y^i]\frac{1}{n+1}[T^{n}]\Big(1+T+\frac{T(1+2T)}{1-yT}\Big)^{n+1}  \nonumber\\
                       \hskip-.22cm &=&\hskip-.22cm  \frac{1}{n+1}[T^{n}][y^i]\Big(1+T+\frac{T(1+2T)}{1-yT}\Big)^{n+1}   \nonumber\\
                       \hskip-.22cm &=&\hskip-.22cm  \frac{1}{n+1}[T^{n}][y^i]\sum_{k=0}^{n+1}\binom{n+1}{k}\frac{1}{(1-yT)^{k}}T^k(1+2T)^k(1+T)^{n+1-k}   \nonumber\\
                       \hskip-.22cm &=&\hskip-.22cm  \frac{1}{n+1}[T^{n}]\sum_{k=0}^{n+1}\binom{n+1}{k}\binom{k+i-1}{i}T^{k+i}(1+2T)^k(1+T)^{n+1-k}     \label{eqn 4.7.1vu}    \\
                       \hskip-.22cm &=&\hskip-.22cm  \frac{1}{n+1}\sum_{k=0}^{n-i}\sum_{j=0}^{k}\binom{n+1}{k}\binom{k+i-1}{i}\binom{k}{j}\binom{n+1-k}{n-k-i-j}2^j.  \nonumber
\end{eqnarray}
Hence, we obtain the following result.
\begin{theorem}
For any integers $n\geq i\geq 0$, there holds
\begin{eqnarray*}
L_{n, i}^{\mathbf{vu}} \hskip-.22cm &=&\hskip-.22cm \frac{1}{n+1}\sum_{k=0}^{n-i}\sum_{j=0}^{k}\binom{n+1}{k}\binom{k+i-1}{i}\binom{k}{j}\binom{n+1-k}{n-k-i-j}2^j.
\end{eqnarray*}
\end{theorem}

\begin{remark} If expanding $(1+2T)^k=\sum_{j=0}^{k}\binom{k}{j}T^{j}(1+T)^{k-j}=\sum_{j=0}^{k}(-1)^{j}\binom{k}{j}2^{k-j}(1+T)^{k-j}$ in (\ref{eqn 4.7.1vu}),
we have another two formulas for $L_{n, i}^{\mathbf{vu}}$, i.e.,
\begin{eqnarray*}
L_{n, i}^{\mathbf{vu}} \hskip-.22cm &=&\hskip-.22cm \frac{1}{n+1}\sum_{k=0}^{n-i}\sum_{j=0}^{k}\binom{n+1}{k}\binom{k+i-1}{i}\binom{k}{j}\binom{n+1-j}{n-k-i-j}  \\
                       \hskip-.22cm &=&\hskip-.22cm \frac{1}{n+1}\sum_{k=0}^{n-i}\sum_{j=0}^{k}(-1)^{j}\binom{n+1}{k}\binom{k+i-1}{i}\binom{k}{j}\binom{n+1-j}{n-k-i}2^{k-j}.
\end{eqnarray*}
\end{remark}

The first values of $L_{n, i}^{\mathbf{vu}}$ are illustrated in Table 4.6.
\begin{center}
\begin{eqnarray*}
\begin{array}{c|ccccccc}\hline
n/i & 0     & 1     & 2      & 3     & 4    & 5        \\\hline
  0 & 1     &       &        &       &      &         \\
  1 & 2     &       &        &       &      &         \\
  2 & 6     &  1    &        &       &      &          \\
  3 & 20    &  8    &  1     &       &      &          \\
  4 & 72    & 48    &  12    & 1     &      &          \\
  5 & 272   & 260   &  100   & 17    &  1   &         \\
  6 & 1064  & 1340  &  706   & 185   &  23  & 1      \\\hline
\end{array}
\end{eqnarray*}
Table 4.6. The first values of $L_{n, i}^{\mathbf{vu}}$.
\end{center}

By (\ref{eqn 1.1}) and (\ref{eqn 4.7vu}), taking the coefficient of $x^{n}y^i$ in $L^{\mathbf{vu}}(x, 2y)$, we derive the following identity and provide a combinatorial proof for it.
\begin{theorem}\label{theom 4.6.1vu}
For any integer $n\geq 0$, there holds
\begin{eqnarray*}
 \sum_{i=0}^{n}2^iL_{n, i}^{\mathbf{vu}}      \hskip-.22cm &=&\hskip-.22cm  2^nC_n.
\end{eqnarray*}
\end{theorem}

\pf Let $\mathcal{L}_{n, i}^{\mathbf{vu}}$ denote the set of weighted G-Motzkin paths of length $n$ with $i$ number of $\mathbf{vu}$-steps such that each $\mathbf{vu}$-step is weighted by $2$ (regarded as $\mathbf{v}_1\mathbf{u}$ and $\mathbf{v}_2\mathbf{u}$ for convenience) and other steps are weighted by $1$. Let $\mathcal{L}_{n, i}^{\mathbf{h}}$ denote the set of weighted G-Motzkin paths of length $n$ with $i$ number of $\mathbf{h}$-steps at high level (level not less than $1$) such that each $\mathbf{h}$-step at hight level is weighted by $2$ (regarded as $\mathbf{h}_1$ and $\mathbf{h}_2$ for convenience) and other steps are weighted by $1$. For convenience, the $\mathbf{h}$-steps at level zero of the   weighted G-Motzkin paths in $\mathcal{L}_{n, i}^{\mathbf{h}}$ are also written as $\mathbf{h}_1$. Let $\mathcal{C}_{n}^{\mathbf{d}}$ denote the weighted Dyck paths of length $2n$ such that each $\mathbf{u}$-step is weighted by $1$ and each $\mathbf{d}$-step is weighted by $2$ (regarded as $\mathbf{d}_1$ and $\mathbf{d}_2$ for convenience). Set $\mathcal{L}_{n}^{\mathbf{vu}}=\bigcup_{i=0}^{n}\mathcal{L}_{n, i}^{\mathbf{vu}}$, $\mathcal{L}_{n}^{\mathbf{h}}=\bigcup_{i=0}^{n}\mathcal{L}_{n, i}^{\mathbf{h}}$ and $\mathcal{L}^{\mathbf{h}}=\bigcup_{n=0}^{\infty}\mathcal{L}_{n}^{\mathbf{h}}$.

Clearly,
$$w(\mathcal{L}_{n}^{\mathbf{vu}})= \sum_{i=0}^{n}w(\mathcal{L}_{n, i}^{\mathbf{vu}})= \sum_{i=0}^{n}2^{i}L_{n, i}^{\mathbf{vu}}\ \ \mbox{and}\ \ w(\mathcal{C}_{n}^{\mathbf{d}})=2^{n}C_{n}. $$

Firstly, there exists a simple bijection $\chi_1$ between $\mathcal{L}_{n, i}^{\mathbf{h}}$ and $\mathcal{L}_{n, i}^{\mathbf{vu}}$ for $n\geq i\geq 0$. For any $\mathbf{P}\in \mathcal{L}_{n, i}^{\mathbf{h}}$, $\chi_1$ can be easily obtained by replacing all the $\mathbf{h}$-steps (at high level) weighted by $2$ in $\mathbf{P}$ by new $\mathbf{vu}$-steps with weight $2$ (regarded as $\mathbf{h}_1\leftrightarrow\mathbf{v}_1\mathbf{u}$ and $\mathbf{h}_2\leftrightarrow\mathbf{v}_2\mathbf{u}$ for convenience) and replacing the old $\mathbf{vu}$-steps in $\mathbf{P}$ by $\mathbf{h}$-steps (at high level), then $\chi_1(\mathbf{P})\in \mathcal{L}_{n, i}^{\mathbf{vu}}$. The inverse procedure can be handled similarly. In fact, $\chi_1$ is an involution between $\mathcal{L}_{n, i}^{\mathbf{h}}$ and $\mathcal{L}_{n, i}^{\mathbf{vu}}$.

Secondly, there exists a recursive bijection $\chi_2$ between $\mathcal{C}_{n}^{\mathbf{d}}$ and $\mathcal{L}_{n}^{\mathbf{h}}$. For $n=0, 1$, we define
$$\chi_2(\varepsilon)=\varepsilon, \ \chi_2(\mathbf{ud}_1)=\mathbf{h}_1,\ \chi_2(\mathbf{ud}_2)=\mathbf{uv}. $$
For any $\mathbf{Q}\in \mathcal{C}_{n}^{\mathbf{d}}$, $\mathbf{Q}$ can be uniquely partitioned into $\mathbf{Q}=\mathbf{Q}_1\mathbf{Q}_2\dots \mathbf{Q}_k$ for $n\geq 2$ and certain $k\geq 1$, where $\mathbf{Q}_1, \mathbf{Q}_2, \dots, \mathbf{Q}_k$ are primitive. Then $\chi_2$ can be recursively defined by
$$ \chi_2(\mathbf{Q})=\chi_2(\mathbf{Q}_1)\chi_2(\mathbf{Q}_2)\dots\chi_2(\mathbf{Q}_k). $$
So it suffices to discuss the cases when $\mathbf{Q}\in \mathcal{C}_{n}^{\mathbf{d}}$ are primitive. There are four cases for such $\mathbf{Q}$ to be considered.

\subsection*{Case 1.} When $\mathbf{Q}=\mathbf{u}\mathbf{Q}'\mathbf{d}_2$, we define $\chi_2(\mathbf{Q})=\mathbf{u}\chi_2(\mathbf{Q}')\mathbf{v}$.

\subsection*{Case 2.} When $\mathbf{Q}=\mathbf{u}\mathbf{Q}'\mathbf{d}_1$, where $\mathbf{Q}'=\mathbf{u}\mathbf{Q}_1\mathbf{d}_1\mathbf{u}\mathbf{Q}_2\mathbf{d}_1\dots\mathbf{u}\mathbf{Q}_k\mathbf{d}_1$ for certain $k\geq 1$, namely, each return step of $\mathbf{Q}'$ is a $\mathbf{d}_1$-step, we define
$$\chi_2(\mathbf{Q})=\mathbf{u}\chi_2(\mathbf{Q}_1)\mathbf{h}_2\chi_2(\mathbf{Q}_2)\dots\mathbf{h}_2\chi_2(\mathbf{Q}_k)\mathbf{d}. $$

\subsection*{Case 3.} When $\mathbf{Q}=\mathbf{u}\mathbf{Q}''\mathbf{Q}'\mathbf{d}_1$, where the last return step of $\mathbf{Q}''$ is a $\mathbf{d}_2$-step and $\mathbf{Q}'=\mathbf{u}\mathbf{Q}_1\mathbf{d}_1\mathbf{u}\mathbf{Q}_2\mathbf{d}_1\dots\mathbf{u}\mathbf{Q}_k\mathbf{d}_1$ for certain $k\geq 1$, i.e., each return step of $\mathbf{Q}'$ is a $\mathbf{d}_1$-step, we define
$$\chi_2(\mathbf{Q})=\mathbf{u}\chi_2(\mathbf{Q}'')\mathbf{h}_2\chi_2(\mathbf{Q}_1)\mathbf{h}_2\chi_2(\mathbf{Q}_2)\dots\mathbf{h}_2\chi_2(\mathbf{Q}_k)\mathbf{v}. $$

\subsection*{Case 4.} When $\mathbf{Q}=\mathbf{u}\mathbf{Q}''\mathbf{Q}'\mathbf{d}_1$, where $\mathbf{Q}''$ is empty or the last return step of $\mathbf{Q}''$ is a $\mathbf{d}_1$-step, and $\mathbf{Q}'=\mathbf{u}\mathbf{Q}_1\mathbf{d}_2\mathbf{u}\mathbf{Q}_2\mathbf{d}_2\dots\mathbf{u}\mathbf{Q}_k\mathbf{d}_2$ for certain $k\geq 1$, i.e., each return step of $\mathbf{Q}'$ is a $\mathbf{d}_2$-step, we define
$$\chi_2(\mathbf{Q})=\mathbf{u}\chi_2(\mathbf{Q}'')\mathbf{h}_2\chi_2(\mathbf{Q}_1)\mathbf{h}_2\chi_2(\mathbf{Q}_2)\dots\mathbf{h}_2\chi_2(\mathbf{Q}_k)\mathbf{v}. $$

Note that according to the definition of $\chi_2$, $\chi_2(\mathbf{Q})$ has no $\mathbf{h}_2$-steps at level zero for any $\mathbf{Q}\in \mathcal{C}_{n}^{\mathbf{d}}$ with $n\geq 0$, and $\chi_2(\mathbf{Q})$ in the above four cases are always primitive for $\mathbf{Q}=\mathbf{ud}_2$ or $\mathbf{Q}$ being primitive with length at least 4. Moreover, in the first case $\chi_2(\mathbf{Q})$ has no $\mathbf{h}_2$-steps at level one; In the third case $\chi_2(\mathbf{Q}'')$ must proceed exactly to the leftmost $\mathbf{h}_2$-step at level one in $\chi_2(\mathbf{Q})$ and end with a $\mathbf{v}$-step, and the last primitive part of $\chi_2(\mathbf{Q}'')$ has no $\mathbf{h}_2$-steps at level one by the first case. In the forth case $\chi_2(\mathbf{Q}'')$ must proceed exactly to the leftmost $\mathbf{h}_2$-step at level one in $\chi_2(\mathbf{Q})$, once $\chi_2(\mathbf{Q}'')$ ends with a $\mathbf{v}$-step, there must exist $\mathbf{h}_2$-steps at level one in the last primitive part of $\chi_2(\mathbf{Q}'')$ recursively by the third or the fourth cases; In the second case $\chi_2(\mathbf{Q})$ is obviously distinguished from the other three cases in the last step.

From the above observation, one can handle the inverse procedure as follows. For any $\mathbf{P}\in \mathcal{L}_{n}^{\mathbf{h}}$, since $\mathbf{P}$ has no $\mathbf{h}_2$-steps at level zero, $\mathbf{P}$ can be uniquely partitioned into $\mathbf{P}=\mathbf{h}_1^{j_1}\mathbf{P}_1\mathbf{h}_1^{j_2}\mathbf{P}_2\dots \mathbf{h}_1^{j_k}\mathbf{P}_k\mathbf{h}_1^{j_{k+1}}$ for $n\geq k\geq 1$ and $j_1, \dots, j_{k+1}\geq 0$, where $\mathbf{P}_1, \mathbf{P}_2, \dots, \mathbf{P}_k\in \mathcal{L}^{\mathbf{h}}$ are primitive. Then $\chi_2^{-1}$ can be recursively defined by
$$ \chi_2^{-1}(\mathbf{P})=(\mathbf{ud}_1)^{j_1}\chi_2^{-1}(\mathbf{P}_1)(\mathbf{ud}_1)^{j_2}\chi_2^{-1}(\mathbf{P}_2)\dots(\mathbf{ud}_1)^{j_k}\chi_2^{-1}(\mathbf{P}_k)(\mathbf{ud}_1)^{j_{k+1}} $$
such that $\chi_2^{-1}(\varepsilon)=\varepsilon$, $\chi_2^{-1}(\mathbf{h}_1)=\mathbf{ud}_1$ and $\chi_2^{-1}(\mathbf{uv})=\mathbf{ud}_2$.
Naturally, it suffices to consider the cases when $\mathbf{P}\in \mathcal{L}_{n}^{\mathbf{h}}$ are primitive. There are four cases for such $\mathbf{P}$ to be considered.

\subsection*{Case 1.} When $\mathbf{P}=\mathbf{u}\mathbf{P}'\mathbf{v}$ and $\mathbf{P}$ has no $\mathbf{h}_2$-steps at level one, i.e., $\mathbf{P}'\in \mathcal{L}_{n-1}^{\mathbf{h}}$ has no $\mathbf{h}_2$-steps at level zero, in this case we define $\chi_2^{-1}(\mathbf{P})=\mathbf{u}\chi_2^{-1}(\mathbf{P}')\mathbf{d}_2$.

\subsection*{Case 2.} When $\mathbf{P}=\mathbf{u}\mathbf{P}'\mathbf{d}$ and $\mathbf{P}$ has exactly $k-1$ $\mathbf{h}_2$-steps at level one, where $\mathbf{P}'=\mathbf{P}_1\mathbf{h}_2\mathbf{P}_2\mathbf{h}_2\dots\mathbf{h}_2\mathbf{P}_k$ with $\mathbf{P}_1, \dots, \mathbf{P}_k\in \mathcal{L}^{\mathbf{h}}$ for $k\geq 1$, in this case we define
$$\chi_2^{-1}(\mathbf{P})=\mathbf{u}\mathbf{u}\chi_2^{-1}(\mathbf{P}_1)\mathbf{d}_1\mathbf{u}\chi_2^{-1}(\mathbf{P}_2)\mathbf{d}_1\dots\mathbf{u}\chi_2^{-1}(\mathbf{P}_k)\mathbf{d}_1\mathbf{d}_1. $$

\subsection*{Case 3.} When $\mathbf{P}=\mathbf{u}\mathbf{P}''\mathbf{P}'\mathbf{v}$ and $\mathbf{P}$ has exactly $k$ $\mathbf{h}_2$-steps at level one for $k\geq 1$, where
$\mathbf{P}'=\mathbf{h}_2\mathbf{P}_1\mathbf{h}_2\mathbf{P}_2\dots\mathbf{h}_2\mathbf{P}_k$ with $\mathbf{P}_1, \dots, \mathbf{P}_k\in \mathcal{L}^{\mathbf{h}}$, and $\mathbf{P}''\in \mathcal{L}^{\mathbf{h}}$ is nonempty and ends with a $\mathbf{v}$-step such that the last primitive part of $\mathbf{P}''$ also has no $\mathbf{h}_2$-steps at level one, in this case we define
$$\chi_2^{-1}(\mathbf{P})=\mathbf{u}\chi_2^{-1}(\mathbf{P}'')\mathbf{u}\chi_2^{-1}(\mathbf{P}_1)\mathbf{d}_1\mathbf{u}\chi_2^{-1}(\mathbf{P}_2)\mathbf{d}_1\dots \mathbf{u}\chi_2^{-1}(\mathbf{P}_k)\mathbf{d}_1\mathbf{d}_1. $$

\subsection*{Case 4.} When $\mathbf{P}=\mathbf{u}\mathbf{P}''\mathbf{P}'\mathbf{v}$ and $\mathbf{P}$ has exactly $k$ $\mathbf{h}_2$-steps at level one for $k\geq 1$, where
$\mathbf{P}'=\mathbf{h}_2\mathbf{P}_1\mathbf{h}_2\mathbf{P}_2\dots\mathbf{h}_2\mathbf{P}_k$ with $\mathbf{P}_1, \dots, \mathbf{P}_k\in \mathcal{L}^{\mathbf{h}}$, $\mathbf{P}''\in \mathcal{L}^{\mathbf{h}}$ is possibly empty or ends with a $\mathbf{z}$-step for $\mathbf{z}\in \{\mathbf{h}_1, \mathbf{h}_2, \mathbf{d}, \mathbf{v}\}$ such that once $\mathbf{P}''$ ends with a $\mathbf{v}$-step, then there exist $\mathbf{h}_2$-steps at level one in the last primitive part of $\mathbf{P}''$, in this case we define
$$\chi_2^{-1}(\mathbf{P})=\mathbf{u}\chi_2^{-1}(\mathbf{P}'')\mathbf{u}\chi_2^{-1}(\mathbf{P}_1)\mathbf{d}_2\mathbf{u}\chi_2^{-1}(\mathbf{P}_2)\mathbf{d}_2\dots \mathbf{u}\chi_2^{-1}(\mathbf{P}_k)\mathbf{d}_2\mathbf{d}_1. $$

The above two procedures verify that $\chi_2$ is indeed a bijection between $\mathcal{C}_{n}^{\mathbf{d}}$ and $\mathcal{L}_{n}^{\mathbf{h}}$. Hence $\chi_1\chi_2$ is a desired bijection between $\mathcal{C}_{n}^{\mathbf{d}}$ and $\mathcal{L}_{n}^{\mathbf{vu}}$.  \qed\vskip0.2cm

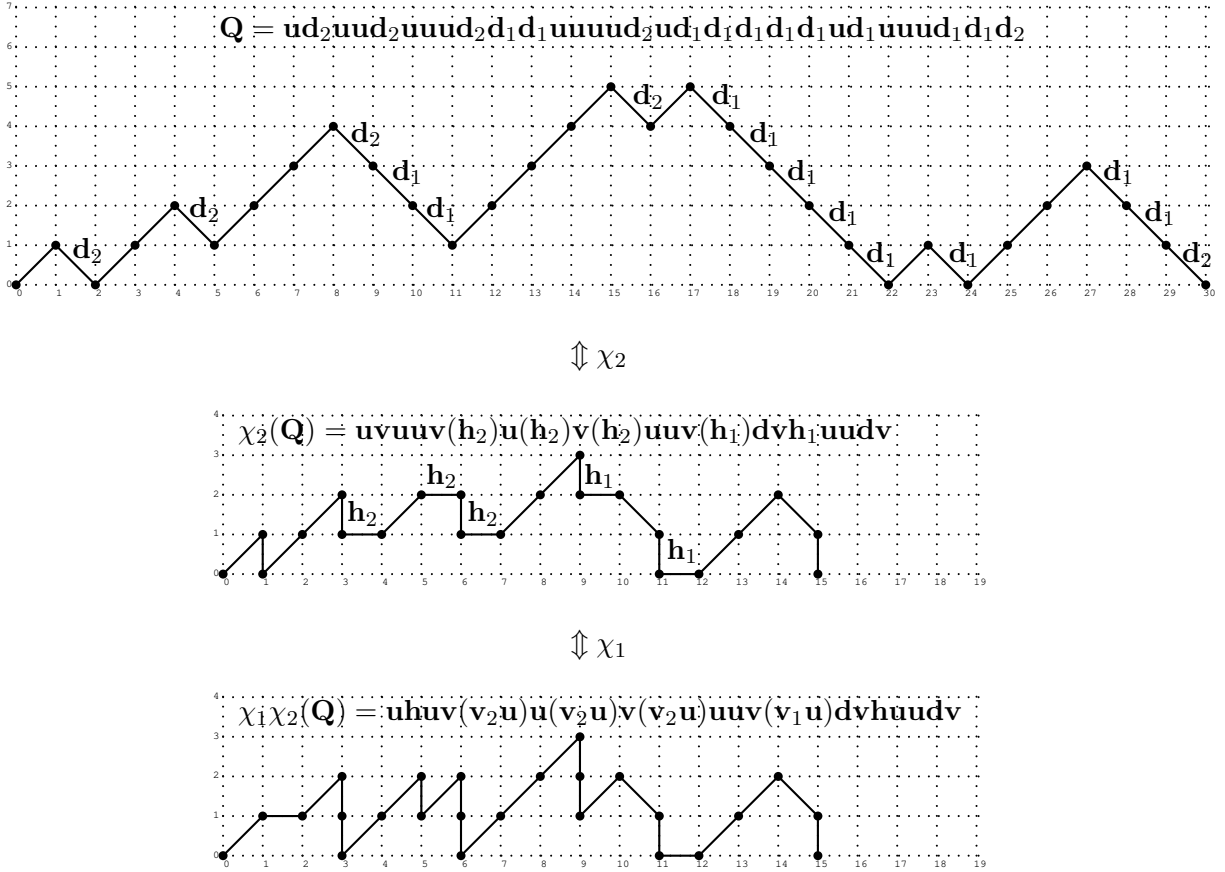
\begin{figure}[h] \setlength{\unitlength}{0.5mm}

\begin{center}
\begin{pspicture}(18,3.5)
\psset{xunit=15pt,yunit=15pt}\psgrid[subgriddiv=1,griddots=4,
gridlabels=4pt](0,0)(30,7)

\psline(0,0)(1,1)(2,0)(4,2)(5,1)(8,4)(11,1)(14,4)(15,5)(16,4)(17,5)(22,0)(23,1)(24,0)(27,3)(30,0)

\pscircle*(0,0){0.06}\pscircle*(1,1){0.06}\pscircle*(2,0){0.06}
\pscircle*(3,1){0.06}\pscircle*(4,2){0.06}\pscircle*(5,1){0.06}
\pscircle*(6,2){0.06}\pscircle*(7,3){0.06}\pscircle*(8,4){0.06}
\pscircle*(9,3){0.06}\pscircle*(10,2){0.06}\pscircle*(11,1){0.06}
\pscircle*(12,2){0.06}\pscircle*(13,3){0.06}\pscircle*(14,4){0.06}\pscircle*(15,5){0.06}
\pscircle*(16,4){0.06}\pscircle*(17,5){0.06}\pscircle*(18,4){0.06}
\pscircle*(19,3){0.06}\pscircle*(20,2){0.06}\pscircle*(21,1){0.06}\pscircle*(22,0){0.06}
\pscircle*(23,1){0.06}\pscircle*(24,0){0.06}\pscircle*(25,1){0.06}
\pscircle*(26,2){0.06}\pscircle*(27,3){0.06}\pscircle*(28,2){0.06}\pscircle*(29,1){0.06}\pscircle*(30,0){0.06}

\put(2.7,3.3){$\mathbf{Q}=\mathbf{u}\mathbf{d}_2\mathbf{uu}\mathbf{d}_2
\mathbf{uuu}\mathbf{d}_2\mathbf{d}_1\mathbf{d}_1\mathbf{uuuu}\mathbf{d}_2
\mathbf{u}\mathbf{d}_1\mathbf{d}_1\mathbf{d}_1\mathbf{d}_1\mathbf{d}_1\mathbf{u}\mathbf{d}_1\mathbf{uuu}\mathbf{d}_1\mathbf{d}_1\mathbf{d}_2$}

\put(.75,.35){$\mathbf{d}_2$}\put(2.3,.9){$\mathbf{d}_2$}\put(4.45,1.9){$\mathbf{d}_2$}\put(5,1.4){$\mathbf{d}_1$}\put(5.45,.9){$\mathbf{d}_1$}
\put(8.2,2.4){$\mathbf{d}_2$}\put(9.25,2.4){$\mathbf{d}_1$}\put(9.75,1.9){$\mathbf{d}_1$}\put(10.25,1.4){$\mathbf{d}_1$}\put(10.8,.85){$\mathbf{d}_1$}\put(11.3,.3){$\mathbf{d}_1$}
\put(12.4,.3){$\mathbf{d}_1$}\put(14.45,1.4){$\mathbf{d}_1$}\put(15,.85){$\mathbf{d}_1$}\put(15.5,.3){$\mathbf{d}_2$}

\end{pspicture}
\end{center}\vskip0.5cm

$\Updownarrow \chi_2$

\begin{center}
\begin{pspicture}(10,2.6)
\psset{xunit=15pt,yunit=15pt}\psgrid[subgriddiv=1,griddots=4,
gridlabels=4pt](0,0)(19,4)

\psline(0,0)(1,1)(1,0)(2,1)(3,2)(3,1)(4,1)(5,2)(6,2)(6,1)(7,1)(9,3)(9,2)(10,2)(11,1)(11,0)(12,0)(14,2)(15,1)(15,0)

\pscircle*(0,0){0.06}\pscircle*(1,0){0.06}\pscircle*(1,1){0.06}\pscircle*(2,1){0.06}
\pscircle*(3,1){0.06}\pscircle*(3,2){0.06}\pscircle*(4,1){0.06}\pscircle*(5,2){0.06}
\pscircle*(6,2){0.06}\pscircle*(6,1){0.06}

\pscircle*(7,1){0.06}\pscircle*(8,2){0.06}\pscircle*(9,2){0.06}
\pscircle*(9,3){0.06}\pscircle*(10,2){0.06}\pscircle*(11,1){0.06}\pscircle*(11,0){0.06}
\pscircle*(12,0){0.06}\pscircle*(13,1){0.06}\pscircle*(14,2){0.06}
\pscircle*(15,0){0.06}\pscircle*(15,1){0.06}

\put(2.7,1.2){$\mathbf{h}_2$}\put(1.65,.7){$\mathbf{h}_2$}\put(3.25,.7){$\mathbf{h}_2$}\put(4.8,1.2){$\mathbf{h}_1$}
\put(5.9, .2){$\mathbf{h}_1$}

\put(.2,1.8){$\chi_2(\mathbf{Q})=\mathbf{u}\mathbf{v}\mathbf{uu}\mathbf{v}(\mathbf{h}_2)
\mathbf{u}(\mathbf{h}_2)\mathbf{v}(\mathbf{h}_2)\mathbf{uu}\mathbf{v}(\mathbf{h}_1)\mathbf{d}\mathbf{v}\mathbf{h}_1\mathbf{uu}\mathbf{dv}$}

\end{pspicture}
\end{center}\vskip0.5cm
$\Updownarrow \chi_1$

\begin{center}
\begin{pspicture}(10,2.5)
\psset{xunit=15pt,yunit=15pt}\psgrid[subgriddiv=1,griddots=4,
gridlabels=4pt](0,0)(19,4)

\psline(0,0)(1,1)(2,1)(3,2)(3,0)(4,1)(5,2)(5,1)(6,2)(6,0)(7,1)(9,3)(9,1)(10,2)(11,1)(11,0)(12,0)(14,2)(15,1)(15,0)

\pscircle*(0,0){0.06}\pscircle*(1,1){0.06}\pscircle*(2,1){0.06}\pscircle*(3,0){0.06}
\pscircle*(3,1){0.06}\pscircle*(3,2){0.06}\pscircle*(4,1){0.06}\pscircle*(5,1){0.06}\pscircle*(5,2){0.06}
\pscircle*(6,2){0.06}\pscircle*(6,0){0.06}\pscircle*(6,1){0.06}

\pscircle*(7,1){0.06}\pscircle*(8,2){0.06}\pscircle*(9,1){0.06}\pscircle*(9,2){0.06}
\pscircle*(9,3){0.06}\pscircle*(10,2){0.06}\pscircle*(11,1){0.06}\pscircle*(11,0){0.06}
\pscircle*(12,0){0.06}\pscircle*(13,1){0.06}\pscircle*(14,2){0.06}
\pscircle*(15,0){0.06}\pscircle*(15,1){0.06}

\put(.2,1.8){$\chi_1\chi_2(\mathbf{Q})=\mathbf{u}\mathbf{h}\mathbf{u}\mathbf{v}(\mathbf{v}_2\mathbf{u})
\mathbf{u}(\mathbf{v}_2\mathbf{u})\mathbf{v}(\mathbf{v}_2\mathbf{u})
\mathbf{uu}\mathbf{v}(\mathbf{v}_1\mathbf{u})\mathbf{d}\mathbf{v}\mathbf{h}
\mathbf{uu}\mathbf{dv}$}

\end{pspicture}
\end{center}\vskip0.2cm

\caption{\small An example of the bijections $\chi_1$ and $\chi_2$ described in the proof of Theorem \ref{theom 4.6.1vu}. }

\end{figure}

In order to give a more intuitive view on the bijection $\chi_1$ and $\chi_2$, a pictorial description of
$\chi_1$ and $\chi_2$ is presented for $\mathbf{Q}=\mathbf{u}\mathbf{d}_2\mathbf{uu}\mathbf{d}_2
\mathbf{uuu}\mathbf{d}_2\mathbf{d}_1\mathbf{d}_1\mathbf{uuuu}\mathbf{d}_2
\mathbf{u}\mathbf{d}_1\mathbf{d}_1\mathbf{d}_1\mathbf{d}_1\mathbf{d}_1\mathbf{u}\mathbf{d}_1\mathbf{uuu}\mathbf{d}_1\mathbf{d}_1\mathbf{d}_2 \in \mathcal{C}_{30}^{\mathbf{d}}$, we have
$$\chi_2(\mathbf{Q})=\mathbf{u}\mathbf{v}\mathbf{uu}\mathbf{v}(\mathbf{h}_2)\mathbf{u}(\mathbf{h}_2)\mathbf{v}(\mathbf{h}_2)
\mathbf{uu}\mathbf{v}(\mathbf{h}_1)\mathbf{d}\mathbf{v}\mathbf{h}_1\mathbf{uu}\mathbf{dv}\in \mathcal{L}_{15}^{\mathbf{h}}\ \ \mbox{and} $$
$$\chi_1\chi_2(\mathbf{Q})=\mathbf{u}\mathbf{h}\mathbf{u}\mathbf{v}(\mathbf{v}_2\mathbf{u})\mathbf{u}(\mathbf{v}_2\mathbf{u})\mathbf{v}(\mathbf{v}_2\mathbf{u})
\mathbf{uu}\mathbf{v}(\mathbf{v}_1\mathbf{u})\mathbf{d}\mathbf{v}\mathbf{h}\mathbf{uu}\mathbf{dv}\in \mathcal{L}_{15}^{\mathbf{vu}}.$$
See Figure 3 for detailed illustrations.

By (\ref{eqn 1.1}) and (\ref{eqn 4.7vu}), taking the coefficient of $x^{n}$ in $L^{\mathbf{vu}}(x, -1)$, we derive that
\begin{theorem}\label{theom 4.6.2vu}
For any integer $n\geq 0$, there holds
\begin{eqnarray*}
 \sum_{i=0}^{n}(-1)^iL_{n, i}^{\mathbf{vu}}   \hskip-.22cm &=&\hskip-.22cm \sum_{i=0}^{n}(-1)^i\binom{n}{i}3^{n-i}C_i.
\end{eqnarray*}
\end{theorem}

By (\ref{eqn 4.7vu}), taking the coefficient of $x^{n}$ in $L^{\mathbf{vu}}(x, -2)=\frac{1}{1-2x}$, we derive the following identity whose combinatorial proof is also provided.
\begin{theorem}\label{theom 4.6.3vu}
For any integer $n\geq 0$, there holds
\begin{eqnarray*}
 \sum_{i=0}^{n}(-2)^iL_{n, i}^{\mathbf{vu}}   \hskip-.22cm &=&\hskip-.22cm 2^{n}.
\end{eqnarray*}
\end{theorem}
\pf Set $\mathcal{L}_{n}^{\mathbf{h}, e}=\bigcup_{i=0}^{[\frac{n}{2}]}\mathcal{L}_{n, 2i}^{\mathbf{h}}$ and $\mathcal{L}_{n}^{\mathbf{h}, o}=\bigcup_{i=0}^{[\frac{n-1}{2}]}\mathcal{L}_{n, 2i+1}^{\mathbf{h}}$. Note that $\chi_1$ in Theorem \ref{theom 4.6.1vu} is a bijection between $\mathcal{L}_{n, i}^{\mathbf{h}}$ and $\mathcal{L}_{n, i}^{\mathbf{vu}}$ for $n\geq i\geq 0$, together with the definition of $\mathcal{D}_n^{*}$ in the proof of (\ref{eqn D2.1}), in order to prove Theorem \ref{theom 4.6.3vu}, it is sufficient to establish a bijection between $\mathcal{L}_{n}^{\mathbf{h}, e}/\mathcal{D}_n^{*}$ and $\mathcal{L}_{n}^{\mathbf{h}, o}$ for $n\geq 0$.

It is trivial for $n=0, 1$. For $n\geq 2$, any $\mathbf{P}\in \mathcal{L}_{n}^{\mathbf{h}, e}/\mathcal{D}_n^{*}$ has at least one of the four subpaths, $\mathbf{h}_1\mathbf{v}$, $\mathbf{h}_2\mathbf{v}$, $\mathbf{d}$ and $\mathbf{uvv}$, find the last one, say $\mathbf{z}$, $\mathbf{P}$ can be partitioned uniquely into $\mathbf{P}=\mathbf{P}_1\mathbf{z}\mathbf{v}^r\mathbf{P}_2$, where $\mathbf{P}_2\in\mathcal{D}_k^{*}$ for certain $0\leq r, k\leq n-2$. Then define $\varrho(\mathbf{P})=\mathbf{P}_1\mathbf{z}'\mathbf{v}^r\mathbf{P}_2$, where
\begin{eqnarray*}
\mathbf{z}'=\left\{
\begin{array}{rl}
 \mathbf{uvv},           &  \mbox{if}\ \mathbf{z}=\mathbf{h}_1\mathbf{v}, \\
 \mathbf{d},             &  \mbox{if}\ \mathbf{z}=\mathbf{h}_2\mathbf{v}, \\
 \mathbf{h}_2\mathbf{v}, &  \mbox{if}\ \mathbf{z}=\mathbf{d},             \\
 \mathbf{h}_1\mathbf{v}, &  \mbox{if}\ \mathbf{z}=\mathbf{uvv}.
\end{array}\right.
\end{eqnarray*}
This way ensures that the number of $\mathbf{h}$-steps at hight level in $\varrho(\mathbf{P})$ is one more or less than that in $\mathbf{P}\in \mathcal{L}_{n}^{\mathbf{h}, e}/\mathcal{D}_n^{*}$, so $\varrho(\mathbf{P})\in \mathcal{L}_{n}^{\mathbf{h}, o}$ and
$\mathbf{z}'$ in $\varrho(\mathbf{P})$ is also the last one of the four subpaths, $\mathbf{h}_1\mathbf{v}$, $\mathbf{h}_2\mathbf{v}$, $\mathbf{d}$ and $\mathbf{uvv}$.
The reverse procedure can be handled similarly. Hence, $\varrho$ is a bijection between $\mathcal{L}_{n}^{\mathbf{h}, e}/\mathcal{D}_n^{*}$ and $\mathcal{L}_{n}^{\mathbf{h}, o}$. This completes the proof. \qed\vskip0.2cm

\begin{theorem}\label{theom 4.6.4vu}
For any integers $n, m\geq 0$, there holds
\begin{eqnarray*}
\sum_{i=0}^{2n}(-1)^{i}\binom{2n}{i}L_{n+m+i+2, m+i+1}^{\mathbf{vu}}=C_n.
\end{eqnarray*}
\end{theorem}
\pf Set $N=n+m+i+3$, one has
\begin{eqnarray*}
L_{n+m+i+2, m+i+1}^{\mathbf{vu}} \hskip-.22cm &=&\hskip-.22cm \sum_{k=0}^{n+1}\sum_{j=0}^{k}\frac{1}{N}\binom{N}{k}\binom{k+m+i}{k-1}\binom{k}{j}\binom{N-k-j}{n+1-k-j}2^j,
\end{eqnarray*}
each inner term in $L_{n+m+i+2, m+i+1}^{\mathbf{vu}}$, denoted by $g_{n, m, k, j}(i)$, is a polynomial on $i$ with degree
$$\partial g_{n, m, k, j}(i)=-1+k+(k-1)+(n+1-k-j)=n+k-1-j\leq 2n $$
such that $\partial g_{n, m, k, j}(i) = 2n$ if and only if $k=n+1$ and $j=0$. Clearly, the leading term in
\begin{eqnarray*}
g_{n, m, n+1, 0}(i) \hskip-.22cm &=&\hskip-.22cm \frac{1}{N}\binom{N}{n+1}\binom{n+m+i+1}{n}=\frac{1}{n+1}\binom{n+m+i+2}{n}\binom{n+m+i+1}{n}
\end{eqnarray*}
is $\frac{i^{2n}}{(n+1)!n!}$. So $L_{n+m+i+2, m+i+1}^{\mathbf{vu}}$ is also a polynomial on $i$ with degree $2n$ such that the leading term is $\frac{i^{2n}}{(n+1)!n!}$. Similar to the proof of Theorem \ref{theom 4.1.2ud}, one can have the result.  \qed\vskip0.2cm

\subsection{ The statistics ``number of $\mathbf{vv}$-steps" }

Let $L^{\mathbf{vv}}(x, y)=\sum_{n=0}^{\infty}\sum_{i=0}^{\infty}L_{n, i}^{\mathbf{vv}}x^ny^i$ be the generating function of $L_{n, i}^{\mathbf{vv}}$, the number of G-Motzkin paths of length $n$ with $i$ $\mathbf{vv}$-steps. According to the
first return decomposition, a G-Motzkin path $\mathbf{P}$ can be decomposed as one of the following six forms:
$$\mathbf{P}=\varepsilon, \ \mathbf{P}=\mathbf{h}\mathbf{P}_1,\ \mathbf{P}=\mathbf{u}\mathbf{P}_1\mathbf{d}\mathbf{P}_{2}, \    \mathbf{P}=\mathbf{u}\mathbf{P}_1\mathbf{u}\mathbf{P}_2\dots \mathbf{u}\mathbf{P}_k\mathbf{u}\mathbf{P}_{k+1}\mathbf{d}\mathbf{v}^{k}\mathbf{P}_{k+2}, $$
$$ \mathbf{P}=\mathbf{u}\mathbf{P}_1\mathbf{u}\mathbf{P}_2\dots \mathbf{u}\mathbf{P}_k\mathbf{h}\mathbf{v}^{k}\mathbf{P}_{k+1}\  \ \mbox{or}\  \mathbf{P}=\mathbf{u}\mathbf{P}_1\mathbf{u}\mathbf{P}_2\dots \mathbf{u}\mathbf{P}_{k-1}\mathbf{u}\mathbf{v}^{k}\mathbf{P}_{k}$$
for certain $k\geq 1$, where $\mathbf{P}_1, \dots, \mathbf{P}_{k+2}$ are (possibly empty) G-Motzkin paths. Note that the last three cases contribute at least $k-1$ $\mathbf{vv}$-steps. Then we get the relation
\begin{eqnarray}
L^{\mathbf{vv}}(x, y)\hskip-.22cm &=&\hskip-.22cm 1+xL^{\mathbf{vv}}(x, y)+x^2L^{\mathbf{vv}}(x, y)^2+\sum_{k=1}^{\infty}x^{k+2}y^{k-1}L^{\mathbf{vv}}(x, y)^{k+2} \nonumber \\
\hskip-.22cm & &\hskip-.22cm \ +\  \sum_{k=1}^{\infty}x^{k+1}y^{k-1}L^{\mathbf{vv}}(x, y)^{k+1}+\sum_{k=1}^{\infty}x^{k}y^{k-1}L^{\mathbf{vv}}(x, y)^{k} \nonumber \\
\hskip-.22cm &=&\hskip-.22cm \big(1+xL^{\mathbf{vv}}(x, y)+x^2L^{\mathbf{vv}}(x, y)^2\big)\Big(1+\frac{xL^{\mathbf{vv}}(x, y)}{1-xyL^{\mathbf{vv}}(x, y)}\Big). \label{eqn 4.8vv}
\end{eqnarray}

Let $T=xL^{\mathbf{vv}}(x, y)$, using the Lagrange inversion formula in (\ref{eqn 4.8vv}), taking the coefficient of $x^{n+1}y^i$ in $xL^{\mathbf{vv}}(x, y)$, we derive that
\begin{theorem}
For any integers $n\geq i\geq 0$, there holds
\begin{eqnarray*}
L_{n, i}^{\mathbf{vv}} \hskip-.22cm &=&\hskip-.22cm \frac{1}{n+1}\sum_{k=0}^{n-i}\sum_{j=0}^{[\frac{n-k-i}{2}]}\binom{n+1}{k}\binom{n+1}{j}\binom{k+i-1}{i}\binom{n-j+1}{n-k-i-2j}.
\end{eqnarray*}
\end{theorem}

The first values of $L_{n, i}^{\mathbf{vv}}$ are illustrated in Table 4.7.
\begin{center}
\begin{eqnarray*}
\begin{array}{c|ccccccc}\hline
n/i & 0     & 1     & 2      & 3     & 4    & 5        \\\hline
  0 & 1     &       &        &       &      &         \\
  1 & 2     &       &        &       &      &         \\
  2 & 6     &  1    &        &       &      &          \\
  3 & 21    &  7    &  1     &       &      &          \\
  4 & 80    &  41   &  11    & 1     &      &          \\
  5 & 322   & 225   &  86    & 16    &  1   &         \\
  6 & 1347  & 1198  &  589   & 162   &  22  &  1       \\\hline
\end{array}
\end{eqnarray*}
Table 4.7. The first values of $L_{n, i}^{\mathbf{vv}}$.
\end{center}

\begin{theorem}\label{theom 4.7.3vv}
For any integers $n, m\geq 0$, there holds
\begin{eqnarray*}
\sum_{i=0}^{2n}(-1)^{i}\binom{2n}{i}L_{n+m+i+2, m+i+1}^{\mathbf{vv}}=C_n.
\end{eqnarray*}
\end{theorem}
\pf Set $N=n+m+i+3$, one has
\begin{eqnarray*}
L_{n+m+i+2, m+i+1}^{\mathbf{vv}} \hskip-.22cm &=&\hskip-.22cm \sum_{k=0}^{n+1}\sum_{j=0}^{[\frac{n-k+1}{2}]}\frac{1}{N}\binom{N}{k}\binom{N}{j}\binom{k+m+i}{k-1}\binom{N-j}{n+1-k-2j},
\end{eqnarray*}
each inner term in $L_{n+m+i+2, m+i+1}^{\mathbf{vv}}$, denoted by $f_{n, m, k, j}(i)$, is a polynomial on $i$ with degree
$$\partial f_{n, m, k, j}(i)=-1+k+j+(k-1)+(n+1-k-2j)=n+k-1-j\leq 2n $$
such that $\partial f_{n, m, k, j}(i) = 2n$ if and only if $k=n+1$ and $j=0$. Clearly, the leading term in
\begin{eqnarray*}
f_{n, m, n+1, 0}(i) \hskip-.22cm &=&\hskip-.22cm \frac{1}{N}\binom{N}{n+1}\binom{n+m+i+1}{n}=\frac{1}{n+1}\binom{n+m+i+2}{n}\binom{n+m+i+1}{n}
\end{eqnarray*}
is $\frac{i^{2n}}{(n+1)!n!}$. So $L_{n+m+i+2, m+i+1}^{\mathbf{vv}}$ is also a polynomial on $i$ with degree $2n$ such that the leading term is $\frac{i^{2n}}{(n+1)!n!}$. Similar to the proof of Theorem \ref{theom 4.1.2ud}, one can have the result.  \qed\vskip0.2cm

\subsection{ The statistics ``number of $\mathbf{du}$-valleys" }

Let $L^{\mathbf{du}}(x, y)=\sum_{n=0}^{\infty}\sum_{i=0}^{\infty}L_{n, i}^{\mathbf{du}}x^ny^i$ be the generating function of $L_{n, i}^{\mathbf{du}}$, the number of G-Motzkin paths of length $n$ with $i$ $\mathbf{du}$-valleys. According to the
first return decomposition, a G-Motzkin path $\mathbf{P}$ can be decomposed as one of the following four forms:
$$\mathbf{P}=\varepsilon, \ \mathbf{P}=\mathbf{h}\mathbf{P}_1,\ \mathbf{P}=\mathbf{u}\mathbf{P}_1\mathbf{d}\mathbf{u}\mathbf{P}_2\mathbf{d}\dots \mathbf{u}\mathbf{P}_k\mathbf{d}\mathbf{u}\mathbf{P}_{k+1}\mathbf{v}\mathbf{P}_{k+2}, \  \mathbf{P}=\mathbf{u}\mathbf{P}_1\mathbf{d}\mathbf{u}\mathbf{P}_2\mathbf{d}\dots \mathbf{u}\mathbf{P}_k\mathbf{d}\mathbf{u}\mathbf{P}_{k+1}\mathbf{d}\mathbf{Q} $$
for certain $k\geq 0$, where $\mathbf{P}_1, \dots, \mathbf{P}_{k+2}$ are (possibly empty) G-Motzkin paths and $\mathbf{Q}$ is empty or begins with an $\mathbf{h}$-step. Note that the last two cases contribute at least $k$ $\mathbf{du}$-valleys. Then we get the relation
\begin{eqnarray*}
L^{\mathbf{du}}(x, y)\hskip-.22cm &=&\hskip-.22cm 1+xL^{\mathbf{du}}(x, y)+\sum_{k=0}^{\infty}x^{2k+1}y^{k}L^{\mathbf{du}}(x, y)^{k+2} \\
\hskip-.22cm & &\hskip-.22cm + \sum_{k=0}^{\infty}x^{2k+2}y^{k}L^{\mathbf{du}}(x, y)^{k+1}(1+xL^{\mathbf{du}}(x, y)) \nonumber \\
\hskip-.22cm &=&\hskip-.22cm 1+xL^{\mathbf{du}}(x, y)+\frac{xL^{\mathbf{du}}(x, y)^2}{1-x^2yL^{\mathbf{du}}(x, y)}+x^2L^{\mathbf{du}}(x, y)\Big(\frac{1+xL^{\mathbf{du}}(x, y)}{1-x^2yL^{\mathbf{du}}(x, y)}\Big).
\end{eqnarray*}
Solve this, we have
\begin{eqnarray}
L^{\mathbf{du}}(x, y)\hskip-.22cm &=&\hskip-.22cm   \frac{1-x-x^2+x^2y-\sqrt{(1-x-x^2+x^2y)^2-4x(1+x^2+xy-x^2y)}}{2x(1+x^2+xy-x^2y)} \nonumber   \\
                     \hskip-.22cm &=&\hskip-.22cm   \frac{1}{1-x-x^2+x^2y}C\Big(\frac{x(1+x^2+xy-x^2y)}{(1-x-x^2+x^2y)^2}\Big).   \label{eqn 4.9.1du}
\end{eqnarray}

By (\ref{eqn 1.1}) and (\ref{eqn 4.9.1du}), taking the coefficient of $x^ny^i$ in $L^{\mathbf{du}}(x, y)$, we derive that
\begin{theorem}
For any integers $n\geq i\geq 0$, there holds
\begin{eqnarray*}
L_{n, i}^{\mathbf{du}} \hskip-.22cm &=&\hskip-.22cm \sum_{k=0}^{n}\sum_{r=0}^{i}\sum_{j=0}^{k-r}\sum_{\ell=0}^{[\frac{n+r-k}{2}]-i-j}(-1)^{i-r}\binom{2k+i-r}{i-r}\binom{k}{r}\binom{k-r}{j}\cdot \\
\hskip-.22cm & &\hskip-.22cm \hskip3cm  \binom{2k+i+\ell-r}{\ell}\binom{n+k-r-i-2j-\ell}{n+r-k-2i-2j-2\ell} C_k.
\end{eqnarray*}
Specially,
\begin{eqnarray*}
L_{n, 0}^{\mathbf{du}} \hskip-.22cm &=&\hskip-.22cm \sum_{k=0}^{n}\sum_{j=0}^{k}\sum_{\ell=0}^{[\frac{n-k}{2}]-j}\binom{2k+\ell}{\ell}\binom{k}{j}\binom{n+k-2j-\ell}{n-k-2j-2\ell} C_k.
\end{eqnarray*}
\end{theorem}

The first values of $L_{n, i}^{\mathbf{du}}$ are illustrated in Table 4.8.
\begin{center}
\begin{eqnarray*}
\begin{array}{c|ccccccc}\hline
n/i & 0      & 1      & 2       & 3     & 4    & 5        \\\hline
  0 & 1      &        &         &       &      &         \\
  1 & 2      &        &         &       &      &         \\
  2 & 7      &        &         &       &      &          \\
  3 & 28     &  1     &         &       &      &          \\
  4 & 123    & 10     &         &       &      &          \\
  5 & 576    & 73     &  1      &       &      &         \\
  6 & 2819   & 485    &  15     &       &      &         \\
  7 & 14250  & 3093   &  154    &  1    &      &         \\
  8 & 73833  & 19325  &  1346   &  21   &      &         \\
  9 & 390048 & 119418 &  10758  &  283  & 1    &         \\\hline
\end{array}
\end{eqnarray*}
Table 4.8. The first values of $L_{n, i}^{\mathbf{du}}$.
\end{center}

\subsection{ The statistics ``number of $\mathbf{dd}$-steps" }

Let $L^{\mathbf{dd}}(x, y)=\sum_{n=0}^{\infty}\sum_{i=0}^{\infty}L_{n, i}^{\mathbf{dd}}x^ny^i$ be the generating function of $L_{n, i}^{\mathbf{dd}}$, the number of G-Motzkin paths of length $n$ with $i$ $\mathbf{dd}$-steps. According to the
first return decomposition, a G-Motzkin path $\mathbf{P}$ can be decomposed as one of the following six forms:
$$\mathbf{P}=\varepsilon, \ \mathbf{P}=\mathbf{h}\mathbf{P}_1,\ \mathbf{P}=\mathbf{u}\mathbf{P}_1\mathbf{v}\mathbf{P}_{2}, \    \mathbf{P}=\mathbf{u}\mathbf{P}_1\mathbf{u}\mathbf{P}_2\dots \mathbf{u}\mathbf{P}_k\mathbf{u}\mathbf{P}_{k+1}\mathbf{v}\mathbf{d}^{k}\mathbf{P}_{k+2}, $$
$$ \mathbf{P}=\mathbf{u}\mathbf{P}_1\mathbf{u}\mathbf{P}_2\dots \mathbf{u}\mathbf{P}_k\mathbf{h}\mathbf{d}^{k}\mathbf{P}_{k+1}\   \mbox{or}\  \mathbf{P}=\mathbf{u}\mathbf{P}_1\mathbf{u}\mathbf{P}_2\dots \mathbf{u}\mathbf{P}_{k-1}\mathbf{u}\mathbf{d}^{k}\mathbf{P}_{k}$$
for certain $k\geq 1$, where $\mathbf{P}_1, \dots, \mathbf{P}_{k+2}$ are (possibly empty) G-Motzkin paths. Note that the last three cases contribute at least $k-1$ $\mathbf{dd}$-steps. Then we get the relation
\begin{eqnarray*}
L^{\mathbf{dd}}(x, y)\hskip-.22cm &=&\hskip-.22cm 1+xL^{\mathbf{dd}}(x, y)+xL^{\mathbf{dd}}(x, y)^2+\sum_{k=1}^{\infty}x^{2k+1}y^{k-1}L^{\mathbf{dd}}(x, y)^{k+2}  \nonumber \\
\hskip-.22cm & &\hskip-.22cm \  +\ \sum_{k=1}^{\infty}x^{2k+1}y^{k-1}L^{\mathbf{dd}}(x, y)^{k+1} + \sum_{k=1}^{\infty}x^{2k}y^{k-1}L^{\mathbf{dd}}(x, y)^{k}   \\
\hskip-.22cm &=&\hskip-.22cm 1+xL^{\mathbf{dd}}(x, y)+xL^{\mathbf{dd}}(x, y)^2+x^2 L^{\mathbf{dd}}(x, y)\Big(\frac{1+xL^{\mathbf{dd}}(x, y)+xL^{\mathbf{dd}}(x, y)^2}{1-x^2yL^{\mathbf{dd}}(x, y)}\Big).
\end{eqnarray*}

However, the exact formula for $L_{n, i}^{\mathbf{dd}}$ is still unknown. Here we give the array
$L_{n, i}^{\mathbf{dd}}$ for $0 \leq n \leq 9$ and $0 \leq i \leq 5$, see Table 4.9.
\begin{center}
\begin{eqnarray*}
\begin{array}{c|ccccccc}\hline
n/i & 0       & 1       & 2      & 3     & 4    & 5        \\\hline
  0 & 1       &         &        &       &      &         \\
  1 & 2       &         &        &       &      &         \\
  2 & 7       &         &        &       &      &          \\
  3 & 29      &         &        &       &      &          \\
  4 & 132     &  1      &        &       &      &          \\
  5 & 641     &  9      &        &       &      &         \\
  6 & 3254    &  64     &  1     &       &      &         \\
  7 & 17060   &  427    &  11    &       &      &         \\
  8 & 91663   &  2770   &  91    &  1    &      &         \\
  9 & 499569  &  20219  &  707   &  13   &      &         \\\hline
\end{array}
\end{eqnarray*}
Table 4.9. The first values of $L_{n, i}^{\mathbf{dd}}$.
\end{center}

\subsection{ The statistics ``number of $\mathbf{dv}$-steps" }

Let $L^{\mathbf{dv}}(x, y)=\sum_{n=0}^{\infty}\sum_{i=0}^{\infty}L_{n, i}^{\mathbf{dv}}x^ny^i$ be the generating function of $L_{n, i}^{\mathbf{dv}}$, the number of G-Motzkin paths of length $n$ with $i$ $\mathbf{dv}$-steps. According to the
first return decomposition, a G-Motzkin path $\mathbf{P}$ can be decomposed as one of the following six forms:
$$\mathbf{P}=\varepsilon, \ \mathbf{P}=\mathbf{h}\mathbf{P}_1,\ \mathbf{P}=\mathbf{u}\mathbf{P}_1\mathbf{d}\mathbf{P}_{2},\   \mathbf{P}=\mathbf{u}\mathbf{P}_1\mathbf{u}\mathbf{P}_2\dots \mathbf{u}\mathbf{P}_{k-1}\mathbf{u}\mathbf{v}^{k}\mathbf{P}_{k}, $$
$$ \mathbf{P}=\mathbf{u}\mathbf{P}_1\mathbf{u}\mathbf{P}_2\dots \mathbf{u}\mathbf{P}_k\mathbf{h}\mathbf{v}^{k}\mathbf{P}_{k+1}\  \mbox{or}\   \mathbf{P}=\mathbf{u}\mathbf{P}_1\mathbf{u}\mathbf{P}_2\dots \mathbf{u}\mathbf{P}_k\mathbf{u}\mathbf{P}_{k+1}\mathbf{d}\mathbf{v}^{k}\mathbf{P}_{k+2}$$
for certain $k\geq 1$, where $\mathbf{P}_1, \dots, \mathbf{P}_{k+2}$ are (possibly empty) G-Motzkin paths. Note that the last case contributes at least one $\mathbf{dv}$-step. Then we get the relation
\begin{eqnarray}
L^{\mathbf{dv}}(x, y)\hskip-.22cm &=&\hskip-.22cm 1+xL^{\mathbf{dv}}(x, y)+x^2L^{\mathbf{dv}}(x, y)^2+\sum_{k=1}^{\infty}x^{k}L^{\mathbf{dv}}(x, y)^{k}  \nonumber \\
\hskip-.22cm & &\hskip-.22cm \  +\ \sum_{k=1}^{\infty}x^{k+1}L^{\mathbf{dv}}(x, y)^{k+1} + y\sum_{k=1}^{\infty}x^{k+2}L^{\mathbf{dv}}(x, y)^{k+2}  \nonumber \\
\hskip-.22cm &=&\hskip-.22cm 1+xL^{\mathbf{dv}}(x, y)+x^2L^{\mathbf{dv}}(x, y)^2+\frac{xL^{\mathbf{dv}}(x, y)\big(1+x L^{\mathbf{dv}}(x, y)+yx^2 L^{\mathbf{dv}}(x, y)^2\big)}{1-xL^{\mathbf{dv}}(x, y)}. \label{eqn 4.10dv}
\end{eqnarray}

Let $T=xL^{\mathbf{dv}}(x, y)$, using the Lagrange inversion formula in (\ref{eqn 4.10dv}), taking the coefficient of $x^{n+1}y^i$ in $xL^{\mathbf{dv}}(x, y)$ in two different ways, we derive that
\begin{theorem}
For any integers $n\geq i\geq 0$, there holds
\begin{eqnarray*}
L_{n, i}^{\mathbf{dv}}
\hskip-.22cm &=&\hskip-.22cm \frac{1}{n+1}\binom{n+1}{i}\sum_{k=0}^{[\frac{n-3i}{2}]}\sum_{j=0}^{n-3i-2k}\binom{n+1-i}{k}\binom{n+1-i-k}{j}\binom{2n-3i-3k-j}{n-3i-2k-j}\\
\hskip-.22cm &=&\hskip-.22cm \frac{1}{n+1}\binom{n+1}{i}\sum_{k=0}^{[\frac{n}{3}]-i}\sum_{j=0}^{[\frac{n}{3}]-i-k}(-1)^{k+j}\binom{n+1-i}{k}\binom{n+1-i-k}{j}\binom{3n-4i-4k-3j+1}{n-3i-3k-3j}
\end{eqnarray*}
\end{theorem}

The first values of $L_{n, i}^{\mathbf{dv}}$ are illustrated in Table 4.10. Clearly, $L_{3n, n}^{\mathbf{dv}}=\frac{1}{3n+1}\binom{3n+1}{n}$ is
the Fuss-Catalan numbers $C_k(n)=\frac{1}{kn+1}\binom{kn+1}{n}$ of the third order $(k=3)$ which also counts the ternary trees with $n$ internal
vertices and counts the number of lattice paths from $(0,0)$ to $(2n,0)$ in the first quadrant of the XOY-plane using an up-steps $\mathbf{u}=(1,1)$ and a down-steps $\mathbf{d}_2=(0,-2)$ \cite[A001764]{Sloane}. This kind of lattice paths of length $2n$ can be easily corresponded bijectively to G-Motzkin paths of length $3n$ with $n$ $\mathbf{dv}$-steps by replacing $\mathbf{d}_2$-steps by $\mathbf{dv}$-steps and vice versa.
\begin{center}
\begin{eqnarray*}
\begin{array}{c|ccccccc}\hline
n/i &    0         & 1      & 2    & 3     & 4    & 5        \\\hline
  0 &    1         &        &      &       &      &         \\
  1 &    2         &        &      &       &      &         \\
  2 &    7         &        &      &       &      &          \\
  3 &    28        &  1     &      &       &      &          \\
  4 &    124       &  9     &      &       &      &          \\
  5 &    584       & 66     &      &       &      &         \\
  6 &    2873      & 443    &  3   &       &      &         \\
  7 &    14592     & 2857   & 49   &       &      &         \\
  8 &    75944     & 18037  & 544  &       &      &         \\
  9 &   402928     & 112510 & 5058 &  12    &      &         \\\hline
\end{array}
\end{eqnarray*}
Table 4.10. The first values of $L_{n, i}^{\mathbf{dv}}$.
\end{center}

\vskip0.5cm
\section{Concluding remarks and further works}

The main objective of this paper has been achieved in Section 2, 3 and 4, where we enumerate ``number of $\mathbf{z}$-steps", ``number of $\mathbf{z}$-steps" at given level, ``number of $\mathbf{z}_1\mathbf{z}_2$-steps" in G-Motzkin paths for $\mathbf{z}, \mathbf{z}_1, \mathbf{z}_2\in \{\mathbf{u}, \mathbf{h}, \mathbf{v}, \mathbf{d}\}$. Some explicit formulas are obtained by bijective and algebraic methods, including generating functions and the Lagrange inversion formula.

Despite that several identities are proved by algebraic methods, we naturally expect their combinatorial proofs, especially for Theorem \ref{theom 2.2.3}, Theorem \ref{theom 3.1.2}, Theorem \ref{theom 3.1.3}, Theorem \ref{theom 3.2.3}, Theorem \ref{theom 3.3.3}, Theorem \ref{theom 3.4.2}, Theorem \ref{theom 4.1.2ud}, Theorem \ref{theom 4.6.2vu}, Theorem \ref{theom 4.6.3vu} and Theorem \ref{theom 4.7.3vv}.

In two forthcoming papers, we further consider the statistics ``number of occurrences of $\Gamma$'' for an arbitrary string $\Gamma$ with at least three steps and the statistics ``number of occurrences of $\Gamma$'' at given height, including even or odd hight, just as done in \cite{ManSapTasTsi, SapTasTsi, Sun}. Moreover, we also discuss the generalized Shr\"{o}der paths with steps in $\{\mathbf{H}, \mathbf{u}, \mathbf{v}_i (i\geq 1), \mathbf{d} \}$ analogous to G-Motzkin paths studied in this paper, where $\mathbf{v}_i=(0,-i)$.


\vskip0.5cm
\section*{Declaration of competing interest}

The authors declare that they have no known competing financial interests or personal relationships that could have
appeared to influence the work reported in this paper.

\section*{Acknowledgements} {The authors are grateful to the referees for
the helpful suggestions and comments. The Project is sponsored by ``Liaoning
BaiQianWan Talents Program". }

\vskip.2cm


\end{document}